\documentclass[12pt]{article}

\usepackage{amsmath, amsthm, amssymb}
\usepackage{setspace}
\usepackage{graphicx}
\usepackage{color}
\def\R{{\mathbb R}}
\def\N{{\mathbb N}}

\def\bpsi{\boldsymbol\psi}
\def\bphi{\boldsymbol\varphi}
\def\bdeta{\boldsymbol\eta}
\def\btau{\boldsymbol\tau}

\def\bn{{\bf n}}
\def\bx{{\bf x}}
\def\by{{\bf y}}
\def\bu{{\bf u}}
\def\bq{{\bf q}}
\def\bv{{\bf v}}
\def\bt{{\bf t}}
\def\bnu{{\bf n}}

\def\bsigma{\boldsymbol\sigma}

\def\dt{{\Delta t}}
\newtheorem {definition} {Definition}[section]

\newtheorem {remark}{Remark}[section]

  \newtheorem {theorem} {Theorem}[section]
  \newtheorem {lemma}   {Lemma}[section]
  
  \newtheorem {corollary}  {Corollary}[section]
  \newtheorem {proposition}{Proposition}[section]

\numberwithin{equation}{section}
\begin{document}

\title{\sc{A Generalization of the Aubin-Lions-Simon Compactness Lemma for Problems on Moving Domains}}
\author{ Boris Muha\thanks{Department of Mathematics,
    University of Zagreb, 10000 Zagreb, Croatia,
    borism@math.hr} \and Sun\v{c}ica \v{C}ani\'{c}\thanks{Department of Mathematics,
    University of Houston, Houston, Texas 77204-3476,
    canic@math.uh.edu} }

\date{}

\maketitle

\begin{abstract}
This work addresses an extension of the Aubin-Lions-Simon compactness result to generalized Bochner spaces  $L^2(0,T;H(t))$, where $H(t)$ is a family of Hilbert spaces,
parameterized by $t$. A compactness result of this type is needed in the study of the existence of weak solutions to nonlinear evolution problems governed by 
partial differential equations defined on moving domains. We identify the conditions on the regularity of the domain motion in time under which our extension of the Aubin-Lions-Simon compactness result holds. Concrete examples of the application
of the compactness theorem are presented, including a classical problem  for the incompressible,
Navier-Stokes equations defined on a {\sl given} non-cylindrical domain, and a class of fluid-structure interaction problems for the incompressible, Navier-Stokes equations,
coupled to the elastodynamics of a Koiter shell. 
The compactness result presented in this manuscript is crucial in obtaining constructive existence proofs to nonlinear, moving boundary problems, using Rothe's method.
\end{abstract}

\section{Introduction}\label{sec:introduction}

We are interested in studying compactness of sequences  in  generalized Bochner spaces $L^2(0,T;H(t))$, 
where $H(t)$ is a family of Hilbert spaces which depend on time.
Problems of this type arise, for example, in studying evolution problems modeled by partial differential equations defined on domains that depend on time.
Examples include general moving-boundary problems such as 
fluid flows problems in time-dependent domains that may either be given {\sl{a priori}}, or in fluid domains that are not known
{\sl{a priori}} but 
depend on the solution of a fluid-structure interaction problem.
In the latter case the elastodynamics of a compliant (elastic, or viscoelastic) structure determines the fluid domain.
Thus, the spatial domain depends on time through 
the unknowns of the problem, giving rise to a strong geometric nonlinearity. 

Problems of this type can be described in general as evolution problems:
\begin{equation}\label{ProblemA}
\begin{array}{rcl}
\displaystyle{\frac{d \bu}{dt}} &=& A^t\bu, t\in(0,T),\\
\bu(0) &=& \bu_0,
\end{array}
\end{equation}
where $A^t:V(t) \to W(t) $  is a family of  (nonlinear) spatial differential operators that depend on time as a parameter. 
For example, $\displaystyle{\frac{d \bu}{dt}} = A^t\bu$ may correspond to the Navier-Stokes equations for an incompressible, viscous fluid 
defined on a moving domain $\Omega(t)$.
In this case, $A^t$ is a spatial differential operator that associates to each $\bu$ the function $\nabla \cdot \boldsymbol\sigma-\bu\cdot\nabla\bu$, 
where $\boldsymbol\sigma$ is the fluid Cauchy stress tensor,
and $\bu$ is divergence free, 
satisfying certain  boundary conditions on  $\Omega(t)$. 

A way to ``solve'' this class of problems is to semi-discretize the problem in time by
sub-dividing the time interval $(0,T)$ into $N$ sub-intervals of width $\dt=T/N$,
and introducing the piecewise constant approximate functions
\begin{equation}\label{FuncDef}
\bu_\dt=\bu_\dt^n\quad {\rm for}\; t \in ((n-1)\dt,n\dt],\;n=1,\dots,N,
\end{equation}
which satisfy, e.g., a backward Euler approximation of the problem on $(t^n,t^{n+1})$:
$$
\frac{\bu_\dt^{n+1} - \bu_\dt^n}{\dt} = A^{t^{n+1}} \bu_\dt^{n+1} \quad {\rm or}\quad 
\frac{\bu_\dt^{n+1} - \bu_\dt^n}{\dt} = A^{t^{n}} \bu_\dt^{n+1}, 
$$
where the choice of $A^{t^{n+1}}$ or $A^{t^{n}}$ depends on the problem at hand.
For example, if the motion of the domain $\Omega(t)$ is specified {\sl{a priori}}, 
$A^{t^{n+1}}$ is typically used, where $A^{t^{n+1}}$ describes an approximation of the spatial differential operator defined on the "current" domain $\Omega(t^{n+1})$.
If the motion of the domain $\Omega(t)$ is not know {\sl{a priori}}, but it depends on the solution of the underlying problem, 
then $A^{t^{n}}$ is typically used, 
where $A^{t^{n}}$ describes an approximation of the spatial differential operator defined on the "previous"  domain $\Omega(t^{n})$.

Functions $\bu_\dt$ are defined for all $t \in(0,T)$ and they are piecewise constant on the interval $((n-1)\dt,n\dt]$,
where the constant is defined by its value at $n\dt$, as stated in \eqref{FuncDef}. This approach to solving the evolution problem 
\eqref{ProblemA}  is sometimes called the Rothe's method.

We are interested in compactness arguments that need to be employed when studying existence of weak solutions
to  \eqref{ProblemA} using Rothe's method. 
Rothe's method provides a constructive proof which uses semi-discretization of the continuous problem with respect to time 
to  design approximate solutions $\{\bu_\dt\}$ where $\dt = T/N$, for every $N \in \N$.
The aim is to prove the existence of  a sub-sequence of $\{\bu_\dt\}$ which converges to a weak solution of \eqref{ProblemA} as $\dt \to 0$,
 or equivalently, as $N \to \infty$. Since the problem is nonlinear, weak convergence is not sufficient to show that 
 the limit is a weak solution of the underlying problem. This is why compactness arguments need to be employed to conclude that
there exists a sub-sequence  $\{\bu_\dt\}$, which is precompact in a certain generalized Bochner space $L^2(0,T;H(t))$. 
 This will allow passage to the limit in nonlinear terms in order to show that the limit, as $\dt \to 0$, 
 of approximate weak solutions satisfies the weak formulation of the continuous problem.

 Employing this strategy to prove the existence of weak solutions to this class of problems
 is highly nontrivial, and is at the center of the current research in this area.
 The main source of difficulties is associated with the fact that for every $N\in\N$ and $n\in\{1,\dots,N\}$, the approximate weak solutions $\bu_\dt^{n}$,
 which are functions of the spatial variables, belong to different solution spaces $V_\dt^{n}$,
 which are associated with the operators $A^{t^n} : V_\dt^n \to W_\dt^n$, and are defined on different domains $\Omega(t^n_\dt)$, thus $V_\dt^n = V(\Omega(t^n_\dt))$.
  We would like to find the conditions under which  $\{ \bu_\dt \}$ is precompact in some $L^2(0,T; H(\Omega_\dt(t)))$,
where the definition of $L^2(0,T; H(\Omega_\dt(t)))$ needs to be made precise.
Namely, we want to find the conditions under which there exists a sub-sequence, also denoted by $\{ \bu_\dt \}$, which converges in $L^2(0,T; H(\Omega_\dt(t)))$ to a function in $L^2(0,T; H(\Omega(t)))$,
as $\dt \to 0$.

  There are two ways how to make the notion of convergence  in $L^2(0,T; H(\Omega_\dt(t)))$ precise.
One way is to introduce a family of {\it mappings}, which map the domains $\Omega(t^n_\dt)$
onto a fixed domain $\Omega$, and work in the space $L^2(0,T;H(\Omega))$.
The other approach is to {\it extend} the functions $\bu_\dt^n$ onto a larger, fixed domain $\Omega^M$, 
and work in the space $L^2(0,T;H(\Omega^M))$.
In both cases, certain conditions describing the regularity in time of the domain motion need to be satisfied, 
in order for a compactness argument to hold.
In this manuscript we identify those conditions, and state a generalization of the Aubin-Lions-Simon compactness result,
which can be used in both approaches, mentioned above.
Namely, we provide a generalization of the Aubin-Lions-Simon compactness result holding in classical Bochner spaces $L^2(0,T;H)$ \cite{aubin1963theoreme,lions1969quelques,Simon}.

We present several concrete examples, which are at the center of today's research in mathematical fluid dynamics, 
where this compactness result is useful.
They include a classical fluid flow problem on a given moving domain, and two
fluid-elastic structure interaction problems, one incorporating the no-slip condition
at the fluid-structure interface, and the other incorporating the Navier slip boundary condition.

\if 1 = 0
--------------------

 give rise to the associated test spaces $Q_\dt^n$ in the weak formulation:
 find $\bu_\dt^{n+1} \in V_\dt^{n+1}$ such that
 \begin{align}
\displaystyle{\left<   \frac{\bu_\dt^{n+1} - \bu_\dt^{n}}{\dt},\bq \right> }&= \left<A^{t^{n}} \bu_\dt^{n+1},\bq \right>, \ \forall \bq \in   Q_\dt^n.
\end{align}
The meaning of the left hand-side in the above equation is a problem, since $\bu_\dt^{n+1}$ and $\bu_\dt^{n}$ are defined on different domains,
and thus belong to different solution spaces.
Additional assumptions need to be made to make 
sense of the problem, and to be able to study convergence of $\bu_\dt$, as $\dt \to 0$.

In this manuscript we identify the conditions on the solution spaces and test spaces that need to be satisfied
in order for this problem to be well-defined, and identify the conditions on the approximate sequences 
that need to satisfied in order for a compactness result to hold.
More precisely, we would like to obtain a compactness result by which one can show that the sequence $\{ \bu_\dt \}$
is precompact in some $L^2(0,T; H(t))$.

Functions $\bu_\dt^n$, $n = 1,\dots,N$, which define $\bu_\dt$, are typically associated with solutions of approximate PDEs, defined on domains $\Omega(t_\dt^n)$, where $t_\dt^n = n\dt$.
Therefore, functions $\bu_\dt^n$ belong to spaces $V(\Omega(t^n_\dt))$, which depend on discrete time. 
We would like to find the conditions under which  $\{ \bu_\dt \}$ is precompact in some $L^2(0,T; H(\Omega_\dt(t)))$,
where the definition of $L^2(0,T; H(\Omega_\dt(t)))$ needs to be made precise.
Namely, we want to find the conditions under which there exists a sub-sequence, also denoted by $\{ \bu_\dt \}$, which converges in $L^2(0,T; H(\Omega_\dt(t)))$ to a function in $L^2(0,T; H(\Omega(t)))$,
as $\dt \to 0$.

We would like to obtain a compactness result by which one can show that a sequence $\{ \bu_i \}_{i=1}^\infty$
is precompact in $L^2(0,T; H(t))$.
As mentioned above, the sequence $\{ \bu_i \}$ may be generated by the Rothe's method, 
in which case we can think of the functions $\bu_i = \bu_\dt$ being generated via the approximations $\bu_\dt^n$, $n = 1,\dots,N$.
Functions $\bu_\dt^n$ are typically associated with solutions of approximate PDEs, defined on domains $\Omega(t_\dt^n)$, where $t_\dt^n = n\dt$.
Therefore, functions $\bu_\dt^n$ belong to spaces $V(\Omega(t^n_\dt))$, which depend on discrete time. 
We would like to find the conditions under which  $\{ \bu_\dt \}$ is precompact in some $L^2(0,T; H(\Omega_\dt(t)))$,
where the definition of $L^2(0,T; H(\Omega_\dt(t)))$ needs to be made precise.
Namely, we want to find the conditions under which there exists a sub-sequence, also denoted by $\{ \bu_\dt \}$, which converges in $L^2(0,T; H(\Omega_\dt(t)))$ to a function in $L^2(0,T; H(\Omega(t)))$,
as $\dt \to 0$.

  There are two ways how to make the notion of convergence  in $L^2(0,T; H(\Omega_\dt(t)))$ precise.
One way is to introduce a family of {\it mappings}, which map the domains $\Omega(t^n_\dt)$
onto a fixed domain $\Omega$, and work in the space $L^2(0,T;H(\Omega))$.
The other approach is to {\it extend} the functions $\bu_\dt^n$ onto a larger, fixed domain $\Omega^M$, 
and work in the space $L^2(0,T;H(\Omega^M))$.
In both cases, certain conditions need to be satisfied, which describe the regularity in time of the domain motion,
in order for our compactness argument to hold.

Our result applies to both approaches, and, in fact holds, for general function spaces, which are not 
necessarily associated with sequence generated by Rothe's method. This is why we state or compactness
result in terms of general Hilbert spaces, and sequences $\{ \bu_\dt \}$ that are not necessarily
generated by Rothe's method.

The compactness result is a generalization of the Aubin-Lions-Simon Lemma for Bochner spaces $L^2(0,T;H)$ \cite{aubin1963theoreme,lions1969quelques,Simon},
to the case where $H=H(t)$  depends on time.
We state and prove the generalized compactness result, and present a series of concrete examples where this compactness
result provides existence of a weak solution. 

\fi

\section{Literature review}
To the best of our knowledge, there is no general compactness theory similar to Aubin-Lions-Simon lemma for spaces $L^2(0,T;H(t))$, where $H$ depends on time. 
There are several compactness results for particular, specific problems for which the spatial domain depends on time, but they were proved using assumptions
that hold for that particular problem at hand.
The first result of this type
can be found in \cite{fujita1970existence}, where the authors studied the Navier-Stokes equations in a given, non-cylindrical domain,
namely, in a domain which depends on time, and whose motion is given {\sl a priori}. 
Similarly, in  \cite{conca2000existence}, the authors used a compactness result to study a fluid-rigid body interaction problem,
where an assumption on high regularity of the domain motion had to be used to obtain an existence result (see also \cite{feireisl2008motion,Starovoitov}). 
In the case of fluid-elastic structure interaction problems, the assumption on high regularity of the interface is typically not satisfied, 
thus different approaches need to be employed. 
 
 In the context of fluid-elastic structure interaction problems, 
 we mention \cite{DEGT,CG} where the authors considered a fluid-elastic structure interaction problem between the flow of a viscous, incompressible fluid and an elastic/viscoelastic plate, in which a compactness argument based on Simon's theorem was used to show
$L^2$-strong convergence of approximate solutions.
 We also mention \cite{hundertmark2016existence} where this approach was used in the case of a non-Newtonian fluid. 
 A similar problem, but in a more general geometrical setting, was studied in \cite{LenRuz}, where compactness of a set of approximate weak solutions,
 based on a particular linearization and regularization of the problem, was proved by using a modification of the ideas from the proof of  Aubin-Lions lemma. 
 Both approaches used the fluid viscosity and kinematic coupling condition to control high frequency oscillations of the structure velocity. 
 Recently,  a version of Aubin-Lions lemma for a moving domain problem was proved in the context of {\it compressible fluids}, see \cite{breit2017compressible},
 where an existence of a solution to an FSI problem between a {compressible fluid} and a linearly elastic shell was obtained.
 The lack of the fluid incompressibility constraint simplifies the compactness argument
for the velocity field.

Compactness results in more general frameworks were studied 
in \cite{nagele2017functional,NageleDiss}, where the authors developed a functional framework based on the flow method and the Piola transform
for problems in \textit{smoothly} moving domains, where the flow causing domain motion was given {\sl a priori}. 
In those works
a version of the Aubin-Lions lemma was obtained within this framework. A different version of Aubin-Lions lemma, in a more general form,
 was also considered in \cite{Moussa2016}. 
The approach in  \cite{Moussa2016} was based on negative Sobolev space-type estimates, defined on non-cylindrical, i.e., time-dependent domains.
The latter approach did not require high degree of smoothness of the domain motion.

We also mention the results obtained in  \cite{SuliKompaktnost,ChenJungelLiu14,Moussa2016},
where generalizations of the Aubin-Lions-Simon lemma in various types of nonlinear settings
were obtained, and the work in \cite{droniou2016gradient}
where  a version of Aubin-Lions-Simon result was obtained in the context of  finite element spaces.

Most of the works mentioned above were obtained for continuous time, i.e., the time variable was not discretized,
and most of them were tailored for a particular application in mind. 
Working with discretized time brings some additional difficulties in terms of the uniform bounds for the time-shifts (translations in time).
In the time-discretized case, namely, for the approaches based on Rothe's semi-discretization method, 
the uniform bounds on the time-shifts need to be somewhat stronger to guarantee compactness,
see Proposition 2 in \cite{dreher2012compact}. In particular, the work in \cite{dreher2012compact}
addresses a version of  Aubin-Lions-Simon result for piecewise constant functions in time, obtained using Rothe's method, but for a problem defined on 
a {\it fixed} Banach space. 

In contrast, our work presented in this manuscript concerns an extension of the Aubin-Lions-Simon result involving 
Hilbert spaces that are {\it solution dependent}, and not necessarily known {\it a priori}. 
To deal with this additional difficulty, we identify a new set of conditions (conditions C below),
which quantify the dependence of the Hilbert spaces on time so that an extension of Aubin-Lions-Simon result can be
applied to a sequence of approximate solutions constructed using Rothe's method.

\section{The general compactness results}\label{sec:compactness}



The goal is to specify the uniform bounds on the sequence $\{ \bu_\dt \}$,
and the precise dependence of the function spaces on time, under which a generalization of the Aubin-Lions-Simon 
compactness result will hold.

For an arbitrary function $\bf f$ denote by $\tau_h {\bf f}$ the time-shifts:
\begin{equation}\label{time_shift}
\tau_h {\bf f}(t,.)={\bf f}(t-h,.),\  h\in\R.
\end{equation}

\if 1 = 0 
\noindent
{\bf Basic assumptions on spaces:}
\begin{itemize}
\item Let $V$ and $H$ be Hilbert spaces such that $V$ is compactly embedded in $H$, $V\subset\subset H$. 
\item For every $N\in\N$ and $n\in\{1,\dots,N\}$, the functions 
 $\bu_\dt^n$ belong to the spaces  $V^n_\dt \hookrightarrow V$, which are continuously embedded in $V$,
 where the embedding is  uniformly continuous w.r.t. $\dt$ and $n$,
 i.e., the embedding constant does not depend on $\dt$ or $n$.
\item  For each  $V^n_\dt$ there is an associated ``test space'' $Q^n_\dt \hookrightarrow V$,
 which is also continuously embedded in $V$, where the embedding is uniformly continuous in $\dt$ and $n$.
 (Notice that we do not require $ Q^n_\dt \subset V^n_\dt$.)
\item 
Let $\{\bu_\dt\} \subset L^2(0,T;H)$.
\end{itemize}

\begin{remark} In the remainder of this manuscript we will be identifying the functions $\bu \in V^n_\dt$ with the corresponding ones in $V$.
More precisely, if ${\cal{E}}^n_\dt:V^n_\dt\to V$ are the embeddings given by the above assumptions,
since there is no possibility of confusion,  in the rest of the paper we will be identifying
$\bu\in V^n_\dt$ with ${\cal{E}}^n_\dt(\bu)\in V$. The same will hold for the test functions.
\end{remark}

\noindent
{\bf Main assumptions on $\{\bu_\dt \}$ and on the time-dependence of function spaces:}


\begin{itemize}
\item[(A)]\textbf{Strong bounds (Uniform Energy Estimates) on $\{ \bu_\dt \}$}. There exists a universal constant $C>0$ such that for every $\dt$ the following estimates hold:
\begin{enumerate}
\item[A1.] The $L^2(0,T;V)$ estimate:
$$
\|\bu_\dt\|^2_{L^2(0,T;V)}\leq C\sum_{n=1}^N\|\bu^n_\dt\|^2_{V^n_\dt}\dt\leq C,
$$
\item[A2.] The $L^\infty(0,T;H)$ estimate:
$$
\|\bu_\dt\|_{L^\infty(0,T;H)}=\max_{n=1,\dots,N}\|\bu^n_\dt\|_H\leq C,
$$
\item[A3.]  Numerical dissipation estimate:
$$
\|\tau_\dt\bu_\dt-\bu_\dt\|^2_{L^2(\dt,T;H)}=\sum_{n=1}^N\|\bu^{n}_\dt-\bu^{n-1}_\dt\|^2_H\dt\leq C\dt.
$$
\end{enumerate}

\item[(B)] \textbf{A uniform (weak) time-derivative bound on $\{ \bu_\dt \}$.}
The exists an absolute constant $C>0$ such that
$$
\|P^{n}_\dt\frac{\bu^{n+1}_\dt-\bu^{n}_\dt}{\dt}\|_{(Q^n_\dt)'}\leq C(\|\bu^{n+1}_\dt\|_{V^{n+1}_\dt}+1),\; n=0,\dots,N-1,
$$
where $P^n_\dt$ is the orthogonal projector onto the space $H^n_\dt=\overline{Q^n_\dt}^H$, which is the closure of ${Q^n_\dt}$ in $H$.

{{\sl Remark:} This property can be slightly modified depending on the problem at hand. 
Namely, on the left hand side one can take $P^{n+1}_\dt$ and $(Q^{n+1}_\dt)'$ instead of $P^{n}_\dt$ and $(Q^{n}_\dt)'$, or/and $\|\bu^n_\dt\|_{V^n_\dt}$. The compactness proof stays the same.}


\item[(C)]\label{FuncSp} \textbf{Smooth dependence of (approximate) function spaces on time}:
\begin{enumerate}
\item[C1.] \textbf{Existence of a common test space.}
For every sub-division $N\in\N$ of the time interval $(0,T)$, and for every $l\in \{1,\dots,N\}$, and $n\in\{1,\dots,N-l\}$,  
there exists a space $Q^{n,l}_\dt\subset V$ which is ``common''' for all the sub-spaces $Q^{n+i}_\dt$, $i=0,1,\dots,l$,
in the sense that there exist the operators 
$${{J}}^i_{\Delta t,l,n}:Q^{n,l}_\dt\to Q^{n+i}_\dt, \quad i=0,1,\dots,l,$$ 
and a universal constant $C>0$ such that for every $\bq\in Q^{n,l}_\dt$  the following holds:
\begin{itemize}
\item{\bf{Continuity of the embeddings, uniform in $\Delta t, n$ and $l$}}:
$$\|J^i_{\Delta t,l,n}\bq\|_{Q^{n+i}_\dt}\leq C\|\bq\|_{Q^{n,l}_\dt},$$ 
\item {\bf{Smooth dependence of the  test spaces $Q^{n+i}_\dt$ on time}}:
\begin{equation}\label{Ji1}
\|J^i_{\Delta t,l,n}\bq-J^j_{\Delta t,l,n}\bq\|_H\leq C |i-j|\dt\|\bq\|_{Q^{n,l}_\dt},\quad i,j\in \{0,\dots,l\},
\end{equation}
\item {\bf Approximation property of $Q^{n,l}_\dt$}:
\begin{equation}\label{Ji2}
\|J^i_{\Delta t,l,n}\bq-\bq\|_H\leq C l\dt\|\bq\|_{Q^{n,l}_\dt},\quad i\in \{0,\dots,l\}.
\end{equation}
Notice that the norms on the left hand-side are all in $H$, i.e., the smooth dependence on (discrete) time
is with respect to the $H$-norm.
\end{itemize}


{\sl Remark:} This condition is related to making sense of taking the time derivatives on different domains.

\item[C2.] \textbf{Approximation property of solution spaces in time.}
Let $V^{n,l}_\dt$ be the closure of  $Q^{n,l}_\dt$ in $V$. 
The approximation property requires that every function in $V^{n+i}_\dt$ can be approximated by a function in $V^{n,l}_\dt$. 
More precisely, we require that there exist the operators 
$$I^i_{\dt,l,n}:V^{n+i}_\dt\to V^{n,l}_\dt, \  i=0,1,\dots,l,$$ 
and a universal constant $C>0$, such that for every $\bv\in V^{n+i}_\dt$ the following holds: 

\begin{itemize}
\item {\bf{Continuity of the embeddings, uniform in $\Delta t, n$ and $l$}}
\begin{equation}\label{C21}
\|I^i_{\dt,l,n} \bv\|_{V^{n,l}_\dt}\leq C\|\bv\|_{V^{n+i}_\dt}, \quad i\in \{0,\dots,l\},
\end{equation}
\item {\bf{Approximation property of  $V^{n,l}_\dt$}}:
\begin{equation}\label{C22}
\|I^i_{\dt,l,n} \bv-\bv\|_H\leq  g(l\dt)\|\bv\|_{V^{n+i}_\dt},\quad i\in \{0,\dots,l\},
\end{equation}
where $g:\R_+\to\R_+$ is a universal, monotonically increasing function such that $g(h)\to 0$ as $h\to 0$.
Notice that the norm on the left hand-side in the second estimate above is again in the space $H$.
\end{itemize}
\item[C3.] \textbf{Uniform Ehrling Property.}
Let  $H^{n,l}_\dt$ be the closure of $Q^{n,l}_\dt$ in $H$. 
By identifying $H_\dt^{n,l}$ with its dual, we have 
$V_\dt^{n,l}\hookrightarrow H^{n,l}_\dt \simeq (H^{n,l}_\dt)' \hookrightarrow (Q^{n,l}_\dt)'$ where the first embedding is compact and the last is injective. 
We require that for all the functions $\bv\in V^{n,l}_\dt$ the following uniform Ehrling property holds:

{\it For every $\delta>0$ there exists a constant $C(\delta)$ independent of $n,l$ and $\dt$, such that
\begin{equation}\label{Ehrling}
\|\bv\|_{H}\leq \delta\|\bv\|_{V^{n,l}_\dt}+C(\delta)\|\bv\|_{(Q^{n,l}_\dt)'},\quad \bv\in V^{n,l}_\dt.
\end{equation}
}
\end{enumerate}

\end{itemize}

We will shown in Sec.~\ref{sec:compactness} that these assumptions provide relative compactness of 
the sequence $\{ \bu_\dt \}$ in $L^2(0,T;H)$. Before we do that, we pause to state a few remarks about the assumptions 
presented above.

\noindent
{\bf{Remarks:}}
\begin{enumerate}
\item Conditions (A) and (B) correspond the ``classical'' conditions in the Aubin-Lions lemma.
They correspond directly to the conditions introduced
for problems on {\emph{fixed domains}} in \cite{dreher2012compact}, with  $V^n_\dt=V$, $Q^{n,l}_\dt=Q^n_\dt=Q$. 
Namely, in that case condition (C) is trivially satisfied.
\item { Condition (A3) can be replaced by (B), as will be shown in Theorem \ref{A3Bez}, below. 
Thus, for compactness, this condition does not have to be satisfied.
However, we leave this condition in because it is usually satisfied as a byproduct of Rothe's method, and it nicely shows the role of numerical dissipation
in the compactness proof, as presented in this manuscript.
 }
\item Conditions (C) are new, and are crucial for compactness arguments on moving domains,
in which case the function spaces depend on domain motion. 
Conditions (C1) and (C2) are ``local'' in the sense that they require smooth dependence of the test spaces and solution spaces on time, locally in time,
namely, for the time shifts from $n\dt$ to time $(n+l)\dt$.
Condition (C3), i.e., the Uniform Ehrling Property, provides a ``global'' estimate.
More precisely, while the classical Ehrling
 lemma holds for any set of function spaces $\{V^{n,l}_\dt,H^{n,l}_\dt,(Q^{n,l}_\dt)'\}$ parameterized by $l,n$ and $\dt$,
where the constants $C(\delta)$ {\it a priori} depend on $l,n$ and $\dt$,
the Uniform Ehrling Property, i.e., condition (C3), requires that those constants {\it do not depend} on the parameters $n,l$ and $\dt$,
 which is a condition on the entire family of function spaces. Thus, the Uniform Ehrling Property requires that the classical Ehrling lemma holds
{\it  uniformly} for the entire family of function spaces $\{V^{n,l}_\dt,H^{n,l}_\dt,(Q^{n,l}_\dt)'\}$ parameterized by $l,n$ and $\dt$. {A version of the uniform Ehrling property was also used in \cite{NageleDiss}, and in \cite{droniou2016gradient} (see Def. 4.15 and Lemma 4.16) in a slightly different context.}
 
\end{enumerate}

\fi 

\begin{theorem}\label{Compactness} {\rm{\bf{(Main Result)}}}
Let $V$ and $H$ be Hilbert spaces such that  $V\subset\subset H$.
Suppose that 
$\{\bu_\dt\} \subset L^2(0,T;H)$ is a sequence such that
$
\bu_\dt (t,\cdot)=\bu_\dt^n(\cdot) \ {\rm on}\; ((n-1)\dt,n\dt],\;n=1,\dots,N,
$
with $N\Delta t = T$.
Let $V_\dt^n$ and $Q_\dt^n$ be Hilbert spaces such that $(V_\dt^n,Q_\dt^n)\hookrightarrow V\times V$,  
where the embeddings are uniformly continuous w.r.t. $\dt$ and $n$,
and $V_\dt^n \subset \subset \overline{Q_\dt^n}^H \hookrightarrow (Q_\dt^n)'$.
Let $\bu_\dt^n \in V_\dt^n$, $n = 1,\dots,N$.
If the following is true: 
\begin{itemize}
\item[(A)] There exists a universal constant $C>0$ such that for every $\dt$ 
\begin{enumerate}
\item[A1.] 
$
\sum_{n=1}^N\|\bu^n_\dt\|_{V^n_\dt}^2\dt\leq C,
$
\item[A2.] 
$
\|\bu_\dt\|_{L^\infty(0,T;H)}\leq C,
$
\item[A3.]
$
\|\tau_\dt\bu_\dt-\bu_\dt\|^2_{L^2(\dt,T;H)}\leq C\dt.
$

\end{enumerate}

\item[(B)] 
There exists a universal constant $C>0$ such that
$$
\|P^{n}_\dt\frac{\bu^{n+1}_\dt-\bu^{n}_\dt}{\dt}\|_{(Q^n_\dt)'}\leq C(\|\bu^{n+1}_\dt\|_{V^{n+1}_\dt}+1),\; n=0,\dots,N-1,
$$
where $P^n_\dt$ is the orthogonal projector onto  $\overline{Q^n_\dt}^H$.

\item[(C)]\label{FuncSp} The function spaces $Q_\dt^n$ and $V_\dt^n$ depend smoothly on time in the following sense:
\begin{enumerate}
\item[C1.] 
For every $\dt > 0$, and for every $l\in \{1,\dots,N\}$ and $n\in\{1,\dots,N-l\}$,  
there exists a space $Q^{n,l}_\dt\subset V$ and the operators 
${{J}}^i_{\Delta t,l,n}:Q^{n,l}_\dt\to Q^{n+i}_\dt, i=0,1,\dots,l,$ such that
$\|J^i_{\Delta t,l,n}\bq\|_{Q^{n+i}_\dt}\leq C\|\bq\|_{Q^{n,l}_\dt}, \ \forall \bq\in Q^{n,l}_\dt$, and
\begin{equation}\label{Ji1}
\Big ((J^{j+1}_{\Delta t,l,n}\bq-J^j_{\Delta t,l,n}\bq),\bu^{n+j+1}_\dt\Big )_H\leq C \dt\|\bq\|_{Q^{n,l}_\dt}\|\bu^{n+j+1}_\dt\|_{V^{n+j+1}_\dt},\quad j\in \{0,\dots,l-1\},
\end{equation}
\begin{equation}\label{Ji2}
\|J^i_{\Delta t,l,n}\bq-\bq\|_H\leq C \sqrt{l\dt}\|\bq\|_{Q^{n,l}_\dt},\quad i\in \{0,\dots,l\},
\end{equation}
where $C > 0$ is independent of $\dt,n$ and $l$.

\item[C2.] 
Let $V^{n,l}_\dt = \overline{Q^{n,l}_\dt}^V$. 
There exist the functions $I^i_{\dt,l,n}:V^{n+i}_\dt\to V^{n,l}_\dt, \  i=0,1,\dots,l,$
and a universal constant $C>0$, such that for every $\bv\in V^{n+i}_\dt$ 
\begin{equation}\label{C21}
\|I^i_{\dt,l,n} \bv\|_{V^{n,l}_\dt}\leq C\|\bv\|_{V^{n+i}_\dt}, \quad i\in \{0,\dots,l\},
\end{equation}
\begin{equation}\label{C22}
\|I^i_{\dt,l,n} \bv-\bv\|_H\leq  g(l\dt)\|\bv\|_{V^{n+i}_\dt},\quad i\in \{0,\dots,l\},
\end{equation}
where $g:\R_+\to\R_+$ is a universal, monotonically increasing function such that $g(h)\to 0$ as $h\to 0$.

\item[C3.] {\rm{Uniform Ehrling property:}}
%
For every $\delta>0$ there exists a constant $C(\delta)$ independent of $n,l$ and $\dt$, such that
\begin{equation}\label{Ehrling}
\|\bv\|_{H}\leq \delta\|\bv\|_{V^{n,l}_\dt}+C(\delta)\|\bv\|_{(Q^{n,l}_\dt)'},\quad \bv\in V^{n,l}_\dt.
\end{equation}
\end{enumerate}
\end{itemize}
then $\{\bu_{\Delta t}\}$ is relatively compact in $L^2(0,T;H)$,
\end{theorem}

\noindent
\begin{remark}
{\rm{The following are general remarks about the statement of the theorem:}}
\begin{enumerate}
{\rm{
\item Conditions (A) and (B) correspond the ``classical'' conditions in the Aubin-Lions lemma.
They correspond directly to the conditions introduced
for problems on {\emph{fixed domains}} in \cite{dreher2012compact}, with  $V^n_\dt=V$, $Q^{n,l}_\dt=Q^n_\dt=Q$. 
Namely, in that case condition (C) is trivially satisfied.
\item { Condition (A3) is not necessary, as will be shown in Theorem \ref{A3Bez}, below. 
However,  this condition is usually satisfied as a by-product of Rothe's method, it simplifies the proof, and it nicely shows the role of numerical dissipation
in compactness arguments. These are the reasons why we left it in the statement and proof of the main theorem, and then showed, in Theorem \ref{A3Bez},
that the result of the main theorem holds true even without (A3). }
\item Conditions (C) are new, and are crucial for compactness arguments on moving domains.
They are directly related to making sense of taking the time derivative on different domains. 
Conditions (C1) and (C2) are ``local'' in the sense that they require smooth dependence of the test spaces and solution spaces on time, locally in time,
namely, for the time shifts from $n\dt$ to time $(n+l)\dt$.
Condition (C3), i.e., the Uniform Ehrling Property, provides a ``global'' estimate.
More precisely, while the classical Ehrling
 lemma holds for any set of function spaces parameterized by $l,n$ and $\dt$, i.e., $\{V^{n,l}_\dt,H^{n,l}_\dt,(Q^{n,l}_\dt)'\}=\{V^{n,l}_\dt,\overline{Q^{n,l}_\dt}^H,(Q^{n,l}_\dt)'\}$,
$V_\dt^{n,l}\hookrightarrow H^{n,l}_\dt \simeq (H^{n,l}_\dt)' \hookrightarrow (Q^{n,l}_\dt)'$ where the first embedding is compact and the last is injective,
with the constants $C(\delta)$ depending on $l,n$ and $\dt$,
the Uniform Ehrling Property, i.e., (C3) requires that those constants {\it do not depend} on the parameters $n,l$ and $\dt$,
 which is a condition on the entire family of function spaces. Thus, the Uniform Ehrling Property requires that the classical Ehrling lemma holds
{\it  uniformly} for the entire family of function spaces $\{V^{n,l}_\dt,H^{n,l}_\dt,(Q^{n,l}_\dt)'\}$ parameterized by $l,n$ and $\dt$. {A version of the uniform Ehrling property was also used in \cite{NageleDiss}, and in \cite{droniou2016gradient} (see Def. 4.15 and Lemma 4.16) in a slightly different context.}
}} 
\end{enumerate}
\end{remark}

\begin{remark} 
{\rm{In the rest of the manuscript we will be identifying the functions $\bu \in V^n_\dt$ with the corresponding ones in $V$.
More precisely, if ${\cal{E}}^n_\dt:V^n_\dt\to V$ are the embeddings given by the above assumptions,
we will be using  $\bu\in V^n_\dt$ to denote ${\cal{E}}^n_\dt(\bu)\in V$. The same will be true for the test functions.}}
\end{remark}


\proof
To prove Theorem~\ref{Compactness} we aim at using Simon's Lemma \cite{Simon}.
Thus, we need to show that  sequence $\{\bu_\dt\}$  satisfies the following
two properties: (1) a spatial compactness property, namely, that the sequence $\{\int_{t_1}^{t_2} \bu_\dt (t) dt\}$ is relatively compact in $H$, for $0 < t_1<t_2<T$,
and (2) a uniform integral equicontinuity property for the time-shifts, namely, $\| \tau_h \bu_\dt - \bu_\dt \|_{L^2(h,T;H)} \to 0$ as $h \to 0$, uniformly in $\dt$.

The spatial compactness property follows from the fact that
the sequence $\{\bu_\dt\}$ is bounded in $L^2(0,T;V)$ and $V$ is compactly embedded in $H$. 

Therefore, it remains to prove the uniform equicontinuity property of the time-shifts in $L^2(0,T;H)$.
More precisely, we will show that for every $\epsilon > 0$, there exists an $h_0 > 0$, such that for every sub-division $N \in \N$ of the time interval $(0,T)$,
giving rise to $\dt = T/N$, the following holds:
\begin{equation}\label{compactness_estimate}
\|\tau_h\bu_\dt-\bu_\dt\|_{L^2(h,T;H)}\leq \varepsilon, \; 0\leq h<h_0.
\end{equation}
The main difficulty in showing this inequality lies in the fact that the functions $\bu_\dt$ are defined via $\bu_\dt^n$ 
on different spaces $V_\dt^n$. This is why properties (C ) will be important.

We start by estimating 
$
\|\tau_h\bu_\dt-\bu_\dt\|_{L^2(h,T;H)},
$
where $h > 0$.
Notice that  for every $\dt > 0$ there exists an $l \in \N_0$ and an $0 \leq s < \dt$ such that $h$ can be written as $h = l \dt + s$.
Thus, we can write $\tau_h\bu_\dt-\bu_\dt = \tau_s\tau_{l\dt} \bu_\dt-\bu_\dt  \pm \tau_{\l\dt}\bu_\dt$,
and use the triangle inequality to estimate:
\begin{equation}\label{estimate6}
\|\tau_h\bu_\dt-\bu_\dt\|_{L^2(h,T;H)}\leq \|\tau_s\tau_{l\dt}\bu_\dt-\tau_{l\dt}\bu_\dt\|_{L^2(h,T;H)}+\|\tau_{l\dt}\bu_\dt-\bu_\dt\|_{L^2(h,T;H)}.
\end{equation}
The first term on the right hand-side is expressed as follows:
\begin{equation}\label{estimate6a}
\|\tau_s\tau_{l\dt}\bu_\dt-\tau_{l\dt}\bu_\dt\|^2_{L^2(h,T;H)} = \sum_{n=l+1}^{N}\int_{n\dt - s}^{n\dt}  \| \tau_{l\dt} \bu_\dt - \tau_{(l+1)\dt} \bu_\dt^{n+1} \|^2_{H}dt
=  s \sum_{n=0}^{N-1-l} \|  \bu_\dt^{l+n} -  \bu_\dt^{l+n+1} \|^2_{H}.
\end{equation}
Now, from property (A3) the sum on the right hand-side is bounded by $C$, and so we have:
\begin{equation}\label{estimate7}
\|\tau_s\tau_{l\dt}\bu_\dt-\tau_{l\dt}\bu_\dt\|_{L^2(h,T;H)} \le C \sqrt{s}.
\end{equation}
To estimate the second term on the right hand-side of \eqref{estimate6} we 
need the following two lemmas.


\begin{lemma}\label{DualEstimate}
Under the assumption that properties (A1), (A2), (B) and (C1) hold,
there exists a uniform constant $C>0$ such that for
every $N\in\N$ large enough, $l\in \{1,\dots,N\}$, $n\in\{1,\dots,N-l\}$ the following estimate holds:
\begin{equation}\label{eq:DualEstimate}
\|P^{n,l}_\dt(\bu_{\Delta t}^{n+l}-\bu^n_\dt)\|_{(Q_\dt^{n,l})'}\leq C\sqrt{l\dt},
\end{equation}
where $P^{n,l}_\dt$ is the orthogonal projector onto $H^{n,l}_\dt$.
\end{lemma}
\proof
Consider a fixed, arbitrary $N\in\N$, $l\in \{1,\dots,N\}$, $n\in\{1,\dots,N-l\}$ and let $\bphi\in Q^{n,l}_\dt$ be a test function in $Q^{n,l}_\dt$. 
For each $i=0,\dots, l$ consider the mappings $J^i_{\dt,l,n} : Q^{n,l}_\dt \to Q^{n+i}_\dt$ that associate to each $\bphi \in  Q^{n,l}_\dt$
the test functions $\bphi^i=J^i_{\dt,l,n} \bphi\in Q^{n+i}_\dt$.
Since there is no chance of confusion, we drop the sub-scripts $\dt,l,n$ in the rest of the proof, and replace $J^i_{\dt,l,n}$ by $J^i$.
 To show the dual estimate \eqref{eq:DualEstimate} we consider the products $(P^{n,l}_\dt(\bu^{n+l}- \bu^n),\bphi)_{H}$.
 Since $P^{n,l}_\dt$ is the orthogonal projector onto $H^{n,l}_\dt$, these are the same as $(\bu^{n+l}- \bu^n,\bphi)_{H}$. 
 We calculate:
$$
(P^{n,l}_\dt(\bu^{n+l}- \bu^n),\bphi)_{H}=(\bu^{n+l}- \bu^n,\bphi)_{H}=\sum_{i=0}^{l-1}(\bu^{n+i+1}-\bu^{n+i},\bphi)_H
=\sum_{i=0}^{l-1}(\bu^{n+i+1}-\bu^{n+i},\bphi^i-\bphi^i+\bphi)_H
$$
\begin{equation}\label{estimate}
=\dt\sum_{i=0}^{l-1}(P^{n+i}_\dt\frac{\bu^{n+i+1}-\bu^{n+i}}{\dt},\bphi^i)_{H}
+\sum_{i=0}^{l-1}(\bu^{n+i+1}-\bu^{n+i},\bphi-\bphi^i)_H.
\end{equation}
The first sum in \eqref{estimate} is estimated as follows: from (B) we first get
$$
\dt\sum_{i=0}^{l-1}(P^{n+i}_\dt\frac{\bu^{n+i+1}-\bu^{n+i}}{\dt},\bphi^i)_{H}
\leq C\dt \sum_{i=0}^{l-1}(\|\bu^{n+i+1}\|_{V^{n+i+1}_\dt}+1)\|\bphi^i\|_{Q^{n+i}_\dt}.
$$
Then, by property (C1), the Cauchy-Schwarz inequality, and property (A1) this is estimated as follows:
$$
\dt \sum_{i=0}^{l-1}(\|\bu^{n+i+1}\|_{V^{n+i+1}_\dt}+1)\|\bphi^i\|_{Q^{n+i}_\dt}
\leq C \dt  \sum_{i=0}^{l-1}(\|\bu^{n+i+1}\|_{V^{n+i+1}_\dt}+1) \|\bphi\|_{Q^{n,l}_\dt}
$$
\begin{equation}
\leq C \dt \left(  \sqrt{\frac{l}{{\dt}}} \sqrt{ \sum_{i=0}^{l-1} \| \bu^{n+i+1} \|_{V^{n+i+1}_\dt} \dt }  + l  \right) \|\bphi\|_{Q^{n,l}_\dt}
\leq {C} (\sqrt{l\dt}+l\dt)\|\bphi\|_{Q_\dt^{n,l}}.
\end{equation}
In the last inequality above, the constant ${C}$ includes the constant from property (A1) since:
$$
\sum_{i=0}^{l-1}\|{\bf u}^{n+i+1}\|_{V^{n+i+1}_{\dt}}
\leq\sqrt{\sum_{i=0}^{l-1}\Delta t\|\nabla {\bf u}^{n+i+1}\|^2_{V^{n+i+1}_{\dt}}}\sqrt{\sum_{i=0}^{l-1}\frac{1}{\Delta t}}\leq C\sqrt{\frac{l}{\Delta t}}.
$$
Therefore,  we have obtained the following estimate:
\begin{equation}\label{estimate2}
\dt\sum_{i=0}^{l-1}(P^{n+i}_\dt\frac{\bu^{n+i+1}-\bu^{n+i}}{\dt},\bphi^i)_{H}
\leq {C} (\sqrt{l\dt}+l\dt)\|\bphi\|_{Q_\dt^{n,l}}.
\end{equation}

To estimate the second term in \eqref{estimate} we use the summation by parts formula:
$$
\sum_{i =  0}^{m} f_i \Delta g_i = - \sum_{i=0}^{m} g_{i+1} \Delta f_i + [f_{m+1} g_{m+1} - f_0 g_0],
$$
where $\Delta g_i := g_{i+1} - g_i$ and $\Delta f_i := f_{i+1} - f_i$. We apply summation by parts to $f _i= \bphi - \bphi^i$ and $g_i = \bu^{n+i}$,
to obtain:
$$
\sum_{i=0}^{l-1}(\bu^{n+i+1}-\bu^{n+i},\bphi-\bphi^i)_H=
-\sum_{i=0}^{l-1}(\bu^{n+i+1},\bphi^{i}-\bphi^{i+1})_H+(\bu^{n+l},\bphi-\bphi^l)_H-(\bu^n,\bphi-\bphi^0)_H.
$$
Now, by using the Cauchy-Schwarz inequality, and by the properties (C1) and (A2), the right hand-side can be estimated as follows:
$$
\leq C\sum_{i=0}^{l-1}\|\bu^{n+i+1}\|_{V^{n+i+1}_\dt}\dt\|\bphi\|_{Q_\dt^{n,l}}+C\sqrt{l\dt}\|\bphi\|_{Q_\dt^{n,l}}(\|\bu^n\|_H+\|\bu^{n+l}\|_H)\leq C\sqrt{l\dt}\|\bphi\|_{Q_\dt^{n,l}}.
$$
Therefore,
\begin{equation}\label{estimate3}
\sum_{i=0}^{l-1}(\bu^{n+i+1}-\bu^{n+i},\bphi-\bphi^i)_H
\leq C\sqrt{l\dt}\|\bphi\|_{Q_\dt^{n,l}}.
\end{equation}
By combining \eqref{estimate}, \eqref{estimate2}, and \eqref{estimate3}, and by recalling that $N$ is large enough, i.e.,
$\dt$ is small, $\sqrt{l\dt}$ is dominant over $l\dt$, and so estimate \eqref{eq:DualEstimate} holds.
\qed

Notice how in this dual space estimate of the shift in $\bu_\dt$ by $l \Delta t$, the right hand-side of the estimate is given exactly in terms of the shift $l \Delta t$.
This will be important for the completion of the compactness proof. 
By using this lemma we can now prove the following estimate, which is crucial for the compactness argument.
\begin{lemma}\label{HEst}
Let $\delta > 0$. Then there exists a constant $C(\delta) >0$ such that for 
every $N\in\N$ large enough, $l\in \{1,\dots,N\}$, $n\in\{1,\dots,N-l\}$  the following estimate holds:
$$
\|\bu^{n+l}_\dt-\bu^n_\dt\|_H\leq C(\|\bu^{n+l}\|_{V^{n+l}_\dt}+\|\bu^{n}\|_{V^{n}_\dt})\big ( g(l\dt)C(\delta)+\delta\big )+C(\delta)\sqrt{l\dt}.
$$
\end{lemma}
\proof
We start by rewriting the left hand-side by adding and subtracting  $I^{l}_\dt\bu^{n+l}_\dt$ and $I^0\bu^{n}_\dt$, and using the triangle inequality
to obtain:
$$
\|\bu^{n+l}_\dt-\bu^n_\dt\|_H\leq \|I^{l}_\dt\bu^{n+l}_\dt-\bu_\dt^{n+l}\|_H+\|I^{l}_\dt\bu^{n+l}_\dt-I^0\bu^{n}_\dt\|_H+\|I^0\bu^n_\dt-\bu^n_{\dt}\|_H,
$$
where
$I^l\bu^{n+l}_\dt$ and $I^0\bu^{n}_\dt$ are the approximations introduced in property (C2).
Notice that we have dropped the sub-scripts $n$ and $l$ in the notation of $I^l_\dt$  to simplify notation 
within this proof.
By property (C2),
the first and the third terms above are bounded as follows:
$$
\|I^{l}_\dt\bu^{n+l}_\dt-\bu_\dt^{n+l}\|_H+\|I^0\bu^n_\dt-\bu^n_\dt\|_H\leq g(l\dt)(\|\bu^{n+l}_\dt\|_{V^{n+l}_\dt}+\|\bu^{n}_\dt\|_{V^{n}_\dt}).
$$
What remains is to estimate the middle term $\|I^{l}_\dt\bu^{n+l}_\dt-I^0\bu^{n}_\dt\|_H$.
For this purpose we use the Uniform Ehrling's Property (C3), which states that for every $\delta > 0$ there exists a uniform constant $C(\delta) > 0$ such
that:
\begin{equation}\label{ehrling1}
\|I^{l}_\dt\bu^{n+l}_\dt-I^0\bu^{n}_\dt\|_H\leq \delta \|I^{l}_\dt\bu^{n+l}_\dt-I^0\bu^{n}_\dt\|_{V^{n,l}_\dt}+C(\delta)\|I^{l}_\dt\bu^{n+l}_\dt-I^0\bu^{n}_\dt\|_{(Q^{n,l}_\dt)'}.
\end{equation}
The first term on the right hand-side of \eqref{ehrling1} is further estimated using the triangle inequality and property (C2):
$$
\|I^{l}_\dt\bu^{n+l}_\dt-I^0\bu^{n}_\dt\|_{V^{n,l}_\dt}
\leq C (\|\bu^{n+l}_\dt\|_{V^{n+l}_\dt}+\|\bu^{n}_\dt\|_{V^{n}_\dt}).
$$
The second term on the right hand-side of \eqref{ehrling1} is estimated as follows:
\begin{equation}\label{estimate4}
\|I^l\bu^{n+l}_\dt-I^0\bu^n_\dt\|_{(Q^{n,l}_\dt)'}
\leq 
\|I^l\bu^{n+l}_\dt-P^{n,l}_\dt\bu^{n+l}_\dt\|_{(Q^{n,l}_\dt)'}
+\|P^{n,l}_\dt\bu^{n+l}_\dt-P^{n,l}_\dt\bu^{n}_\dt\|_{(Q^{n,l}_\dt)'}
+\|P^{n,l}_\dt\bu^{n}_\dt-I^0\bu^{n}_\dt\|_{(Q^{n,l}_\dt)'}.
\end{equation}
Now, for the first and last terms above we use the following:
$$
\|I^i\bu^{n+i}_\dt-P^{n,l}_\dt\bu^{n+i}\|_{(Q^{n,l}_\dt)'} \leq C \|I^i\bu^{n+i}_\dt-P^{n,l}_\dt\bu^{n+i}\|_H
\leq C \left(  \|I^i\bu^{n+i}_\dt-\bu^{n+i}_\dt\|_H +\|\bu^{n+i}_\dt-P^{n,l}_\dt\bu^{n+i}_\dt\|_H \right),
$$
where $i = 0,l$.
Since orthogonal projection is the best approximation, we get that the two terms on the right hand-side can be estimated by:
$$
\|I^i\bu^{n+i}_\dt-\bu^{n+i}_\dt\|_H+\|\bu^{n+i}_\dt-P^{n,l}_\dt\bu^{n+i}\|_H\leq 2\|I^i\bu^{n+i}_\dt-\bu^{n+i}_\dt\|_H\leq g(l\dt)\|\bu^{n+i}_\dt\|_{V^{n+i}_\dt},
$$
where we used property (C2) in the last inequality. Therefore, the first and the last terms on the right hand-side in \eqref{estimate4} are estimated as follows:
\begin{equation}
\|I^l\bu^{n+l}_\dt-P^{n,l}_\dt\bu^{n+l}_\dt\|_{(Q^{n,l}_\dt)'} + \|P^{n,l}_\dt\bu^{n}_\dt-I^0\bu^{n}_\dt\|_{(Q^{n,l}_\dt)'}
\leq g(l\dt)  (\|\bu^{n+l}_\dt\|_{V^{n+l}_\dt}+\|\bu^{n}_\dt\|_{V^{n}_\dt}).
\end{equation}
Notice that we used the monotonicity of the function $g$ above.

What is left is to estimate the middle term in \eqref{estimate4}, i.e., $\|P^{n,l}_\dt\bu^{n+l}_\dt-P^{n,l}_\dt\bu^{n}_\dt\|_{(Q^{n,l}_\dt)'} $. 
Here we use Lemma~\ref{DualEstimate} to obtain:
\begin{equation}\label{estimate5}
\|P^{n,l}_\dt\bu^{n+l}_\dt-P^{n,l}_\dt\bu^{n}_\dt\|_{(Q^{n,l}_\dt)'} 
\leq 
C \sqrt{l \Delta t}.
\end{equation}

Finally, by putting all the estimates together we get:
$$
\|\bu^{n+l}_\dt-\bu^n_\dt\|_H\leq C(\|\bu^{n+l}_\dt\|_{V^{n+l}_\dt}+\|\bu^{n}_\dt\|_{V^{n}_\dt})\big ( (g(l\dt)(1+C(\delta))+\delta\big )+C(\delta)\sqrt{l\dt}.
$$
Since constant $C(\delta)$ comes from the uniform Ehrling property, this constant is usually large, and so the number $1$
on the right hand-side can be ``swallowed'' by the constant
$C(\delta)$, giving rise to the estimate in the statement of Lemma~\ref{HEst}.
\qed


Based on these two lemmas we can now continue the proof of the main compactness theorem
by estimating the second term on the right hand-side of \eqref{estimate6}.
We start by the estimate from Lemma \ref{HEst} and square both sides of the estimate in Lemma~\ref{HEst} to obtain:
\begin{equation}\label{EstimateL}
\|\bu^{n+l}_\dt-\bu^n_\dt\|^2_H\leq C\big (\|\bu^{n+l}_\dt\|^2_{V^{n+l}_\dt}+\|\bu^{n}_\dt\|^2_{V^{n}_\dt}\big ) (\delta+C(\delta)g(l\dt))^2+C(\delta)^2l\dt.
\end{equation}
By multiplying \eqref{EstimateL} by $\dt$ and summing w.r.t. $n$ (a discrete analogue of integration with respect to $t$) we obtain:
\begin{equation}\label{estimate10}
\|\tau_{l\dt}\bu_\dt-\bu_\dt\|^2_{L^2(l\dt,T;H)}\leq C\|\bu_\dt\|^2_{L^2(0,T;V)} (\delta+C(\delta)g(l\dt))^2+C(\delta)^2l\dt.
\end{equation}
Since all the terms on the right hand-side are positive, the right hand-side is less than or equal to the complete square of the quantities on the right hand-side.
By taking the square root on both sides we obtain:
\begin{equation}\label{EstimateLSum}
\|\tau_{l\dt}\bu_\dt-\bu_\dt\|_{L^2(l\dt,T;H)}
\leq 
{C\|\bu_\dt\|_{L^2(0,T;V)} (\delta+C(\delta)g(l\dt))+C(\delta)\sqrt{l\dt}}.
\end{equation}
Finally, by combining estimates \eqref{estimate6}, \eqref{estimate7}, and \eqref{EstimateLSum} we get:
$$
\|\tau_h\bu_\dt-\bu_\dt\|_{L^2(h,T;H)} \leq C\sqrt{s} + C\|\bu_\dt\|_{L^2(0,T;V)} (\delta+C(\delta)g(l\dt))+C(\delta)\sqrt{l\dt}.
$$
Now, for every $N \in \N_0$ there exist an $l$ and an $s$ such that 
$h=l\dt + s$. We will use this, together with the assumption that $g$ is a monotonically increasing function,
 to estimate from above the terms containing $l\dt$ and $s$, by the appropriate terms containing $h$. We obtain: 
\begin{equation}\label{estimate8}
\|\tau_h\bu_\dt-\bu_\dt\|_{L^2(h,T;H)} \leq C\sqrt{h} + C\|\bu_\dt\|_{L^2(0,T;V)} (\delta+C(\delta)g(h))+C(\delta) \sqrt{h}.
\end{equation}
This inequality provides the desired compactness estimate \eqref{compactness_estimate}. 
More precisely, from \eqref{estimate8} we see that for every $\varepsilon > 0$ we can find a $\delta > 0$ 
in the uniform Ehrling's property (C3) so that
$C\|\bu_\dt\|_{L^2(0,T;V)}\delta<\varepsilon/2$. This defines a constant $C(\delta) > 0$.
With the given $C(\delta)$ we can choose an $h_0 > 0$ such that 
the sum of the three terms on the right hand-side in \eqref{estimate8} containing $h$ 
is less than or equal to $\varepsilon/2$.
Since $g$ is monotonically increasing, any $h$ such that $0 \le h < h_0$ will give 
$$
\|\tau_h\bu_\dt-\bu_\dt\|_{L^2(h,T;B)}\leq \varepsilon.
$$
Therefore, we have just shown that
for every $\varepsilon > 0$, there exists an $h_0 > 0$ such that 
for every $N\in\N$
$$
\|\tau_h\bu_\dt-\bu_\dt\|_{L^2(h,T;B)}\leq \varepsilon, \; 0\leq h<h_0.
$$
This completes the proof of Theorem~\ref{Compactness}.
\qed

\begin{theorem}\label{A3Bez}
The conclusion of Theorem \ref{Compactness} is valid without condition (A3).
\end{theorem}
\proof
First, we note that the only place where we use condition (A3) in the proof of Theorem \ref{Compactness} is to prove estimate \eqref{estimate6}. 
We will now obtain an estimate of $\|\tau_h \bu_\dt-\bu_\dt\|_{L^2(h,T;H)}$ without assuming (A3). Using \eqref{estimate6a} we have:
$$
\|\tau_s\tau_{l\dt}\bu_\dt-\tau_{l\dt}\bu_\dt\|^2_{L^2(0,T-h;H)} = 
\frac{s}{\dt}\dt\sum_{n=0}^{N-l-1} \| \bu_\dt^{l+n} -  \bu_\dt^{l+n+1} \|^2_{H}\leq\frac{s}{\dt}\|\tau_\dt\bu_\dt-\bu_\dt\|^2_{L^2(\dt,T;H)}.
$$
Now, instead of \eqref{estimate8} we have:
\begin{equation}\label{new_estimate}
\|\tau_h\bu_\dt-\bu_\dt\|_{L^2(h,T;H)} \leq\frac{s}{\dt}\|\tau_\dt\bu_\dt-\bu_\dt\|_{L^2(0,T;H)} + C\|\bu_\dt\|_{L^2(0,T;V)} (\delta+C(\delta)g(l\dt))+C(\delta)\sqrt{l\dt}.
\end{equation}
We will show that for every $\varepsilon>0$, there exists an $h^* > 0$, such that for every $\dt > 0$
$$
\|\tau_h\bu_\dt-\bu_\dt\|_{L^2(h,T;B)}\leq \varepsilon, \; 0\leq h<h^*.
$$
In particular, we only need to estimate the first term on the right hand-side of \eqref{new_estimate}. Namely, we want to show that
for every $\varepsilon>0$, there exists an $\tilde{h}_0 > 0$, such that for every $\dt > 0$
$$
\frac{s}{\dt}\|\tau_\dt\bu_\dt-\bu_\dt\|_{L^2(0,T;H)} < \frac{\varepsilon}{2}.
$$

We have two possibilities:

\begin{enumerate}
\item $\dt \leq h$. In this case the proof is the same as before. Namely, by \eqref{estimate10} we have:
$$
\frac{s}{\dt}\|\tau_\dt\bu_\dt-\bu_\dt\|_{L^2(\dt,T;H)}
\leq \|\tau_\dt\bu_\dt-\bu_\dt\|_{L^2(\dt,T;H)}
$$
$$
\leq {C} (\delta+C(\delta)g(\dt))+C(\delta)\sqrt{\dt}
\leq {C} (\delta+C(\delta)g(h))+C(\delta)\sqrt{h}\leq \varepsilon.
$$
Recall, $l=1$. Thus, since $g$ is monotonically increasing, we see that again, for every $\varepsilon > 0$,
there exists an $h_1 > 0$, such that the above inequality holds for all $h < h_1$.
\item $\dt > h = s$.
From \eqref{conv_to_0}, we see that
there exists a $\dt'>0$ such that
$$
\|\tau_\dt\bu_\dt-\bu_\dt\|_{L^2(0,T;H)}<\varepsilon,\; \dt<\dt'.
$$
Again, we need to estimate the following expression uniformly in $\dt$:
$$
\frac{s}{\dt}\|\tau_\dt\bu_\dt-\bu_\dt\|_{L^2(0,T;H)}.
$$
Again there are two cases:
\begin{itemize}
\item[$2a$)]$\dt<\dt'$: Here, we immediately have 
$$
\frac{s}{\dt}\|\tau_\dt\bu_\dt-\bu_\dt\|_{L^2(0,T;H)}
\leq \|\tau_\dt\bu_\dt-\bu_\dt\|_{L^2(0,T;H)}<\varepsilon,
$$
for all $\dt < \dt'$, and all $h$ (in particular, here we have $h < \dt$).
\item[$2b$)] $\dt>\dt'$: 
Here we first notice that by taking $l=1$ in \eqref{EstimateLSum}, using the same $\varepsilon$, $\delta$ argument as before, we have
\begin{equation}\label{conv_to_0}
\|\tau_\dt\bu_\dt-\bu_\dt\|_{L^2(\dt,T;H)}\to 0,
\end{equation}
as $\dt \to 0$.
Therefore there exists an $M>0$ such that:
\begin{equation}\label{F1A3}
\|\tau_\dt\bu_\dt-\bu_\dt\|_{L^2(\Delta t,T;H)}\leq M,\; \dt>0.
\end{equation}

Now we can take $h_2$ to be such that $h_2<\varepsilon\frac{\dt'}{M}$ to obtain that
$$
\frac{s}{\dt}\|\tau_\dt\bu_\dt-\bu_\dt\|_{L^2(0,T;H)}
\leq \frac{h_2}{\dt'}M<\varepsilon,
$$
for all $h < h_2$.
\end{itemize}
\end{enumerate}
By taking  $\tilde{h}_0=\min\{h_1,h_2\}$ we have shown that for every $\varepsilon > 0$, there exists an $\tilde{h}_0>0$, such that
$$
\frac{s}{\dt}\|\tau_\dt\bu_\dt-\bu_\dt\|_{L^2(0,T;H)} < \frac{\varepsilon}{2},
$$
for all $h < \tilde{h}_0$, which completes the proof.
\qed

\begin{remark}{\rm{
Theorem~\ref{Compactness} provides the ``classical'' $L^2$-compactness of Kolmogorov type in the case of moving domains. 
More precisely, let $\Omega\subset\R^d$ be a bounded domain, and $\bdeta$ a Lipschitz function defined on $[0,T]\times\Omega$
such that for every $t\in [0,T]$, $\bdeta(t,.)$ 
is bi-Lipschitz. Introduce the space-time domain
$$
\Omega^T=\{(t,\bx)\in [0,T]\times\R^d:\bx\in\bdeta(t,\Omega)\}
$$
defined by the family of functions $\bdeta(t,\cdot)$.
Since  $\bdeta$ are continuous, there exists a cylindrical ``super''-domain $\tilde{\Omega}^T=(0,T)\times\Omega_M$ such that $\Omega^T\subset\tilde{\Omega}^T$
($\Omega_M$ contains all the spatial domains defined by the family $\bdeta(t,\cdot)$). 

Consider, for example, a family of functions $\{ \bu_\dt \} \subset L^2(\Omega^T)$ such that 
$\bu_\dt (t,\cdot)=\bu_\dt^n(\cdot) \ {\rm on}\; ((n-1)\dt,n\dt],\;n=1,\dots,N,$
with $N\Delta t = T$, where $\bu^n_\dt \in H^s(\eta(t\dt,\Omega))$, $s < 1/2$. 
Then, the Fr\' echet-Kolmogorov Theorem (Thrm.~4.24 \cite{brezis2010functional}) implies that $\{ \bu_\dt \}$ is precompact in $L^2(\Omega^T)$.

Indeed, the conditions for the Fr\' echet-Kolmogorov Theorem are that the family of functions $\{ \bu_\dt \} \subset L^2(\Omega^T)$ needs to be bounded,
and that the translations both in space and time satisfy $\lim_{h\to 0} \| \tau_h \bu_\dt -\bu_\dt \|_{L^2(\Omega^T)} \to 0$, uniformly with respect to $\dt$. 
This follows directly from Theorem~\ref{Compactness}. More precisely, if we 
set $H=L^2(\Omega_M)$, $V=H^s(\Omega_M)$, $s<1/2$, $V^n_\dt=H^s(\eta(t\dt,\Omega))$, $Q^n_\dt=H^m_0(\eta(t\dt,\Omega))$, $m\in\N$,
with $(Q^n_\dt)'=H^{-m}(\eta(t\dt,\Omega))$, since $s<1/2$ the embedding $(V^n_\dt,Q^n_\dt)\hookrightarrow V\times V$ is given by extending the 
functions by $0$ to $\Omega_M$. Now, the boundedness of $\bu_\dt$ and the uniform boundedness of space translations 
is a direct consequence of property (A). To prove uniform boundedness of time translations, one follows the main steps of the proof of Theorem~\ref{Compactness},
where interpolation and the boundedness of the time-derivatives imply the desired property for the time-translations.
}}
\end{remark}

\section{Examples}

\subsection{The Navier-Stokes equations defined on a given moving domain}\label{example:NS}

As a first example we consider the flow of a viscous, incompressible fluid in a moving, time-dependent domain, whose motion is given {\sl a priori}
by a function $\bdeta=\bdeta(t,\bx)$, which measures the displacement
 from a fixed, reference Lipschitz domain $\Omega \subset \R^d$.
This is a classical problem in fluid dynamics, studied by many authors.
In particular, the existence of a {local in time} strong solution by Rothe's method was obtained by Ladyzhenskaya in \cite{ladyzhenskaya1970initial},
and by a different method, using a suitable change of variables, by Inoue and Wakimoto \cite{inoue1977existence}. 
Weak solutions were studied by various methods, such as Rothe's method \cite{NeustupaRothe09}, the penalty method \cite{fujita1970existence}, 
a change of variables by which the problem was mapped onto a fixed domain \cite{MyTe}, and by elliptic regularization \cite{rodolfo1985existence}. 
Depending on the technique used, different assumptions on the smoothness of the domain motion were needed.
Our approach, based on the compactness result presented in Theorem~\ref{Compactness},
has several advantages 
over the methods discussed above.
It uses semi-disretization in time for the Navier-Stokes equations defined in a cylindrical domain, as presented
in \cite{Tem}, Ch III.4, which is then combined with 
an ALE mapping to deal with the motion of the fluid domain. 
This approach based on the ALE-type discretization of the time-derivative in  moving domains 
is popular in numerical computations,
 and the proof that we present below shows that the numerical schemes, based on this semi-discretization in time,  are convergent. 
 Moreover, our proof works in any dimension, see \cite{Tem}, and relaxes the regularity assumptions on $\bdeta$,
 used in other works. 
 This is especially important for generalizing the result to moving boundary problems, and to fluid-structure interaction problems
 where the regularity of the fluid domain is not known {\sl apriori}, as we discuss in Section~\ref{sec:FSI}.

\subsubsection{Problem definition}

We assume that domain motion, described by $\bdeta$, is such that:
\begin{enumerate}
\item $\bdeta\in W^{1,\infty}([0,T]\times\Omega)$,
\item $\bdeta(t,.):\Omega\to\Omega(t)$ is a volume preserving bijective mapping, i.e. det$\nabla \bdeta(t,.)=1$, $t\in [0,T],$ 
\item $\bdeta$ is Lipschitz in the sense that:
\begin{itemize}
\item there exist constants $c_L, C_L > 0$ such that $\bdeta$ is bi-Lipschitz:
\begin{equation}\label{Bi-Lip}
c_L|\bx-\by|\leq |\bdeta(t,\bx)-\bdeta(t,\by)|\leq C_L|\bx-\by|,
\end{equation}
\item  $\bdeta$ is bi-Lipschitz in both space and time $(t,\bx)$ with the same Lipschitz constant $C_L$:
\begin{equation}\label{Lip}
|\bdeta(s,\bx)-\bdeta(t,\by)|\leq C_L(|\bx-\by|+|s-t|).
\end{equation}
\end{itemize}
\end{enumerate}
We define the moving  domain $\Omega(t)$ at time $t$ through $\bdeta$ as $\Omega(t)=\bdeta(t,\Omega)$, and denote
$$
(0,T)\times\Omega(t)=\cup_{t\in (0,T)}\Omega(t).
$$
With a slight abuse of notation we denote by $\bdeta^{-1}$ the function defined by $t\mapsto \bdeta(t,.)^{-1}$.

We are interested in studying the flow of an incompressible, viscous fluid in the cylinder $(0,T)\times\Omega(t)$,
where the flow is modeled by the Navier-Stokes equations
for an incompressible, viscous fluid, with the no-slip boundary condition at the moving boundary:
$$\bu(t,.)=\partial_t\bdeta(t,.)\circ\bdeta^{-1}(t,.),$$
describing continuity of velocities between the fluid and the boundary.
To simplify matters, since this condition can be treated in the same way as a homogeneous Dirichlet condition on $\partial\Omega(t)$, 
with an additional mild regularity assumption on $\partial_{tt}\bdeta\in L^2(0,T;H^{-1}(\Omega(t))$, $\partial_t\bdeta\in L^2(0,T;H^1(\Omega(t))$,
we consider here, as in \cite{ladyzhenskaya1970initial}, the corresponding initial-boundary value problem with zero Dirichlet data on the boundary:

\begin{align}\label{NSMoving}
\partial_t\bu+\bu\cdot\nabla\bu-\Delta\bu+\nabla p&=0,\;
\nabla\cdot\bu=0\;{\rm in}\; (0,T)\times\Omega(t),\\
\label{NSMovingBC}
\bu(t,.)&=0\quad{\rm on}\;\partial\Omega(t),\\
\label{NSMovingIC}
\bu(0,.)&=\bu_0\quad{\rm in}\;\partial\Omega(0).
\end{align}

\begin{definition}
A function 
$\bu\in L^{\infty}(0,T;L^2(\Omega(t)))\cap L^2(0,T;H^1(\Omega(t)))$ is a weak solution to problem \eqref{NSMoving}-\eqref{NSMovingIC}
 if $\nabla\cdot\bu=0$ in the sense of distributions, and:
\begin{equation}
\int_0^T\int_{\Omega(t)}\Big (-\bu\cdot\partial_t\bq+(\bu\cdot\nabla)\bu\cdot\bq +\nabla\bu:\nabla\bq \Big )=
\int_{\Omega(0)}\bu_0\cdot\bq(t,.),
\end{equation}
for all $\bq\in C_c^{\infty}([0,T)\times\Omega(t))$ such that $\nabla\cdot\bq=0$.
\end{definition}

\subsubsection{The semi-discretized problem}\label{sec:semi}
We construct a sequence of approximate solutions by semi-discretizing the problem in time.
Let $\dt = T/N$, where $N$ corresponds to the sub-division of the time interval $(0,T)$ into $N$ sub-intervals.
Since the fluid domain is moving, at every time step we need to solve a problem defined on a different domain:
 $\Omega^n_\dt=\bdeta(n\dt,\Omega)$, $n \in \{1,\dots,N\}$.
 This introduces difficulties in defining the discretized time derivative, since each of the functions, defined by a finite difference approximation,
is defined on a different domain.
To get around this difficulty, we map the current domain onto the previous domain via an ALE mapping  $A^n_\dt$, defined by:
\begin{equation}
A^n_\dt:\Omega^{n+1}_\dt\to\Omega^{n}_\dt, \ {\rm such \ that}\ A^n_\dt=\bdeta(n\dt,.)\circ\bdeta((n+1)\dt,.)^{-1}.
\end{equation}
Notice that $A^n$ is volume preserving, i.e. det$\nabla A^n=1$.  Moreover, the ALE velocity, i.e., the domain velocity, is given by
$${\bf w}^n_\dt=\partial_t\bdeta(n\dt,.)\circ\bdeta(n\dt,.)^{-1}.$$
Approximate solutions $\bu^n_\dt$ are then defined by the following time-marching scheme:
\begin{equation}\label{Ex1AppSol}
\begin{array}{c}
\displaystyle{\int_{\Omega^{n+1}_\dt}\big (\frac{\bu^{n+1}_\dt-{\bu}^n_\dt\circ A^n_\dt}{\dt}}+(({\bu}^n_\dt-{\bf w}^{n+1}_\dt)\cdot\nabla)\bu^{n+1}_\dt\big )\cdot\bq+\int_{\Omega^{n+1}_\dt}\nabla\bu^{n+1}_\dt:\nabla\bq=0,\; \bq\in C^{\infty}_c(\Omega^{n+1}_\dt),\\
{\rm with} \ \nabla\cdot\bq=\nabla\cdot\bu^{n+1}_\dt=0\ {\rm in \ distribution \ sense}.
\end{array}
\end{equation}

\begin{remark}{\rm{
We will be assuming that the functions $\bu^n_\dt$ are extended by $0$ outside of $\Omega^n_\dt$. }}
\end{remark}
By this remark, the nonlinear advection term is well-defined, since, outside $\Omega^{n}_\dt$ the function $\bu^n_\dt$ is zero, and outside 
$\Omega^{n+1}_\dt$ the function $\bu^{n+1}_\dt$ is zero. In the limit as $\dt \to 0$, the ``error'' will converge to zero, because
the functions $\bdeta$, which determine the fluid domains, are continuous both in space and time, and so the 
characteristic functions of the difference $\Omega^{n+1}_\dt \setminus \Omega^n_\dt$ will converge to zero strongly, as $\dt \to 0$.
Therefore, passing to the limit as $\dt \to 0$ will give the correct nonlinear advection term defined on the current fluid domain $\Omega(t)$. 

\begin{remark}{\rm{
Instead of the ALE formulation \eqref{Ex1AppSol}, one can formulate the approximate solutions by taking $\bu^n_\dt$ instead of ${\bu}^n_\dt\circ A^n_\dt$, where $\bu^n_\dt$ is defined on $\Omega^{n+1}_\dt$ by using extension by $0$. In that case there is no ALE derivative term $\bf w$. This approach was used by Ladyzhenskaya in \cite{ladyzhenskaya1970initial} to prove the existence of a strong solution to problem \eqref{NSMovingIC}. However, this approach is not suitable for more complicated moving boundary problems, since extensions by $0$ cannot be used.
This, together with the fact that ALE is widely used in numerical computations, are the reasons why we prefer
to use the ALE formulation presented above.
}}
\end{remark}

The existence of the $\bu^n_\dt$'s satisfying \eqref{Ex1AppSol} can be shown in a classical way, see, e.g., Temam \cite{Tem}, Chapter III.4, or
\cite{BorSun} for the existence proof of a similar problem. 
Furthermore, it can be shown that such a solution satisfies the following energy estimate.

\begin{lemma}
There exists an absolute constant $C>0$ such that for each $\dt>0$ and $n\in \{1,\dots,N\}$
solution $\bu^n_\dt$ of \eqref{Ex1AppSol} satisfies the following energy estimate:
\begin{equation}\label{EnergyNS}
\|\bu^n_\dt\|_{L^2(\Omega^n_\dt)}^2+\sum_{k=0}^{n-1}\|\bu^{k+1}_\dt-({\bu}_\dt^k \circ A^k_\dt) \|^2_{L^2(\Omega^{k+1}_\dt)}+\sum_{k=1}^{n}\|\nabla\bu^n_\dt\|_{L^2(\Omega^k_\dt)}^2\dt\leq C\|\bu_0\|^2_{L^2(\Omega)}.
\end{equation}
\end{lemma}

For each fixed $\dt > 0$ approximations $ \bu^n_\dt$ define a function $\bu_\dt$, which is piecewise constant on each sub-interval 
$(n\dt, (n+1) \dt) \subset (0,T)$, as in \eqref{FuncDef}. The goal of this section is to show that the compactness result from 
Sec.~\ref{sec:compactness} can be applied to the sequence $ \{ \bu^n_\dt\}$, defined 
by \eqref{Ex1AppSol}. Namely, we want to show that the sequence $\{ \bu^n_\dt\}$  is precompact in $L^2(0,T;H)$,
where $H$ is a Hilbert space defined below.

\subsubsection{Compactness result}

To state the main compactness result of this section we have to define the relevant function spaces. 
For this purpose we first introduce the ``maximal'' fluid domain $\Omega_M$, which contains all the 
time-dependent domains $\Omega(t)$, i.e., 
$\Omega_M$ is a domain such that $\Omega(t)\subset\Omega_M$, $t\in [0,T]$. 
We will be assuming that $\Omega(t) \subset \Omega_M \subset \R^d$, where $d$ is the dimension of the physical space.
As a result, the fluid velocity is a vector function with $d$ components. 
We introduce the following function spaces:
\begin{itemize}
\item The overarching function spaces: 
\begin{equation}\label{Ex1FSa}
H=L^2(\Omega_M)^d,\; V=H^1_0(\Omega_M)^d,
\end{equation}
\item The approximation function spaces:
\begin{equation}\label{Ex1FSb}
 V^{n}_\dt=\{\bu\in H^1_0(\Omega^{n}_\dt)^d:\nabla\cdot\bu=0\},\; Q^n_\dt=H^s(\Omega^{n}_\dt)^d\cap V^{n}_\dt,\ s>n/2,
\end{equation}
\item The function spaces associated with the time shifts $\tau_h \bu(t,\cdot)$, where $h = l \dt$, for $l \in \{1,\dots,N\}$:
\begin{equation}\label{Ex1CTspace}
Q^{n,l}_\dt=\{\bq\in H^s(\Omega^{n,l}_\dt)^d\cap V^n_\dt:\nabla\cdot\bq=0\},
\end{equation}
where $\Omega^{n,l}_\dt$ is a common domain contained in all the time-shifts 
by $i \dt$, $i \in \{0,\dots,l\}$, constructed below.
\end{itemize}
\begin{remark}{\rm{
In the rest of this section we assume that all the functions are extended by $0$ to $\Omega_M$.}}
\end{remark}

To construct the domain $\Omega^{n,l}_\dt$ from \eqref{Ex1CTspace}, which is contained
in all the time-shifts, we consider  
$$
\Omega_{\gamma}=\{\bx\in\Omega:{\rm dist}(\bx,\partial\Omega)>\gamma\},
$$
where $\gamma > 0$. 
For $\Omega$ Lipschitz, and $\gamma$ small enough, $\Omega_\gamma$ is also a Lipschitz domain. 
We would like to find a $\gamma > 0$, which depends on $h= l \dt$, so that the image $\bdeta(n\dt,\Omega_\gamma)$ of $\Omega_\gamma$ is contained 
in all the domains $\Omega^{n,i}_\dt$, for $i\in \{0,\dots,l\}$.
For this purpose we prove the following simple geometric lemma (see also \cite{Moussa2016} Propositions 4 and 5) . 
\begin{lemma}\label{LemGeom1}
Let $t\in[0,T]$. Then for $\gamma = 2\frac{c_L}{C_L}h$, where $c_L$ and $C_L$ are the Lipschitz constants given by
\eqref{Bi-Lip} and \eqref{Lip}, and $h = l\dt$,
we have
$\bdeta(t,\Omega_{\gamma})\subset\bdeta(s,\Omega)$, for all $s\in [0,T]$, such that $|s-t|<h$.
\end{lemma}
\proof
The proof is a consequence of the Lipschitz continuity of $\bdeta$ and the bounds on the Lipschitz constants in \eqref{Bi-Lip} and \eqref{Lip}.
More precisely, we will show that the boundary $\partial(\bdeta(s,\Omega))$ of the domain corresponding to time $s$,
cannot intersect the boundary $\bdeta(t,\Omega_\gamma)$, which is at least $\gamma$ units away from the boundary of the domain 
corresponding to the time $t$, as long as $| s - t | < h$, where $\gamma = 2\frac{c_L}{C_L}h$.
Take $\bx'\in \partial(\bdeta(t,\Omega))$. Since $\bdeta(t,.)$ is injective and continuous, there exists an $\bx\in\partial\Omega$ such that $\bx'=\bdeta(t,\bx)$. Moreover, let $\by'\in \bdeta(t,\Omega_\gamma)$, $\by'=\bdeta(t,\by)$, $\by\in \Omega_{\gamma}$. Then by \eqref{Bi-Lip} we have:
\begin{equation}\label{NSTmp1}
|\bx'-\by'|=|\bdeta(t,\bx)-\bdeta(t,\by)|\geq c_L|\bx-\by|\geq 2C_Lh.
\end{equation}
Now, let us take $\tilde{\bf x}=\bdeta(s,{\bx})\in\partial(\bdeta(s,\Omega))$, ${\bf x}\in\partial\Omega$. From \eqref{Lip} we have:
$$
|\tilde{\bf x}-\bx'|=|\bdeta(s,{\bx})-\bdeta(t,{\bx})|\leq C_L|s-t|\leq C_L h.
$$
Therefore, from \eqref{NSTmp1} and $\gamma = 2\frac{c_L}{C_L}h$, we conclude that $\tilde{\bx}\notin\bdeta(t,\Omega_{\gamma})$, which proves the lemma.
\qed

\begin{corollary}\label{corollary:shifts}
Domain $\Omega^{n,l}_\dt :=\bdeta(n\dt,\Omega_{2\frac{c_L}{C_L}l\dt})$ has the property that
$\Omega^{n,l}_\dt\subset\Omega^{n+i}_\dt$, $\forall i\in\{0,\dots,l\}$. 
\end{corollary}

Before we proceed, we mention that one can use several approaches inside our framework to prove the existence of a weak solution. 
One natural approach would be to use the change of variables introduced in [11] to transfer the test function between the fluid domains at different time. 
Namely, the change of variables from [11] preserves the divergence free condition provided that there exists a volume preserving diffeomorphism 
between different fluid domains, and is therefore also suitable for transforming the test functions. 
However, this approach requires higher regularity assumptions on $\bdeta$. 
We propose a different approach, which relaxes the regularity assumptions on $\bdeta$.

We continue by defining the velocity functions that depend on both time and space by introducing
$\{\bu_\dt\} \subset L^2(0,T;H)$  which are piecewise constant in $t$, i.e.
\begin{equation*}
\bu_\dt=\bu_\dt^n
\quad {\rm on}\; ((n-1)\dt,n\dt],\;n=1,\dots,N,
\end{equation*}
 as well as the corresponding time-shifts, denoted by $\tau_h$, defined by
$
\tau_h \bu_\dt (t,.)=\bu_\dt(t-h,.),\  h\in\R.
$

\begin{theorem}{\bf{(Compactness result for moving boundary problem)}}
The sequence $\{ \bu_\dt \}$, introduced in Sec.~\ref{sec:semi}, satisfying the weak formulation \eqref{Ex1AppSol},
and energy estimate \eqref{EnergyNS}, is relatively compact in $L^2(0,T;(L^2(\Omega_M))^d)$.
\end{theorem}

\proof
We would like to show that the assumptions (A)-(C) from Theorem~\ref{Compactness} 
hold true.

\noindent
\textbf{Property A: Strong bounds.}
The strong bounds (A1) and (A2) follow directly from the energy inequality for approximate solutions \eqref{EnergyNS}.

To show (A3) we need to show that
$$
\|\tau_\dt\bu_\dt-\bu_\dt\|^2_{L^2(\dt,T;L^2(\Omega_M))}=\sum_{n=1}^N\|\bu^{n}_\dt-\bu^{n-1}_\dt\|^2_{L^2(\Omega_M))}\dt\leq C\dt.
$$
Indeed, from the triangle inequality we have:
$$
\sum_{n=1}^N\|\bu^n_\dt-\bu^{n-1}_\dt\|^2_{L^2(\Omega_M)}
\leq \sum_{n=1}^N\|\bu^n_\dt-  {\bu}^{n-1}_\dt \circ A^{n-1}_\dt \|^2_{L^2(\Omega_M)}
+\sum_{n=1}^N\| {\bu}^{n-1}_\dt \circ A^{n-1}_\dt -\bu^{n-1}_\dt\|^2_{L^2(\Omega_M)}.
$$ 
Now, the first term on the right hand side can be bounded by the constant from the energy estimate \eqref{EnergyNS}. 
The second term will be estimated using the following lemma:
\begin{lemma}\label{NSDiffLemma}
There exists a constant $C>0$, independent of $\dt$ and $n$,  such that
\begin{equation}\label{NsTimeDiff}
\| {\bu}^{n}_\dt \circ A^{n}_\dt -\bu^{n}_\dt\|_{L^2(\Omega_M)}\leq C\|\nabla\bu^n_\dt\|_{L^2(\Omega_M)}\dt.
\end{equation}
\end{lemma}
\proof
By the mean value theorem we get:
$$
|\int_{\Omega_M}|({\bu}^{n}_\dt \circ A^n_{\dt})-\bu^n_{\dt}|^2 d\bx
\leq \int_{\Omega_M}\int_0^1|\nabla\bu^n(sA^n_{\dt}(\bx)-(1-s)\bx)|^2|A^n_{\dt}(\bx)-{\rm id}|^2dsd\bx
$$
$$
\leq C\dt^2\|\partial_t\bdeta\|^2_{L^{\infty}}\|\nabla\bu^n\|^2_{L^2}.
$$
Here ${\rm id}$ denotes the identity mapping. In the last inequality
 we have used the Lipschitz property in time of $\bdeta$, i.e., of the ALE mappings $A^n_\dt$, to obtain 
 the desired estimate.
\qed

\noindent
{\bf Property B: Uniform time-derivative bound.}
We want to estimate the following norm:
$$
\|P^{n+1}_\dt\frac{\bu^{n+1}_\dt-\bu^n_\dt}{\dt}\|_{(Q^{n+1}_\dt)'} = \sup_{\|\bq \|_{Q^{n+1}_\dt}=1}
\Big |\int_{\Omega^{n+1}_\dt}\frac{{\bf u}^{n+1}_\dt- {\bf u}^{n}_\dt}{\Delta t}\cdot{\bq} d\bx \Big |.
$$

We rewrite the integral on the right hand-side as follows:
$$
\int_{\Omega^{n+1}_\dt}\frac{\bu^{n+1}_\dt-\bu^n_\dt}{\dt}\cdot\bq=
\int_{\Omega^{n+1}_\dt}\frac{\bu^{n+1}_\dt -{\bu}^{n}_\dt \circ A^{n}_\dt}{\dt}\cdot\bq
+\int_{\Omega^{n+1}_\dt}\frac{{\bu}^{n}_\dt \circ A^{n}_\dt -\bu^n_\dt }{\dt}\cdot\bq,\quad \bq\in Q^{n+1}_\dt,
$$
and estimate each of the two integrals on the right hand-side. 
The first one is estimated by using the weak form \eqref{Ex1AppSol} of the problem:
$$
|\int_{\Omega^{n+1}_\dt}\frac{\bu^{n+1}-{\bu}^{n}_\dt \circ A^{n}_\dt}{\dt}\cdot\bq|\leq |\int_{\Omega^{n+1}_\dt}
({\bu}^{n}_\dt  \cdot\nabla)\bu^{n+1}_\dt)\cdot\bq
+\int_{\Omega^{n+1}_\dt}\nabla\bu^{n+1}_\dt:\nabla\bq|
$$
$$
\leq C\|\bu^n_\dt\|_H\|\bu_\dt^{n+1}\|_{V^{n+1}_\dt}\|\bq\|_{L^{\infty}}+\|\bu_\dt^{n+1}\|_{V^{n+1}_\dt}\|\bq\|_{V}\leq C\|\bu_\dt^{n+1}\|_{V^{n+1}_\dt}\|\bq\|_{H^s(\Omega^{n+1}_\dt)}.
$$
The second one is estimated by using Lemma \ref{NSDiffLemma} and Cauchy-Schwarz inequality.

\qed


\noindent
\textbf{Property C: Smooth dependence of (approximate) function spaces on time.}

\noindent
\textbf{Property C1: The existence of a common test space.} 
We set $Q^{n,l}_\dt$ to be the common test space.
By Corollary \ref{corollary:shifts}  $Q^{n,l}_\dt\subset Q^{n+i}_\dt$, $i=0,\dots,l$, and so
$J^j\bq=\bq$,
where we consider the functions defined outside the minimal domain using extensions by $0$.

\noindent 
\textbf{Property C2: Approximation property of solution spaces. } 
We define $V^{n,l}_\dt$ to be the closure of $Q^{n,l}_\dt$ in $V = H_0^1(\Omega_M)^d$, and construct the operators
$$
I^i_{\dt,l,n} : V_\dt^{n+i} \to V_\dt^{n,l}
$$
for which we need to show that \eqref{C21} and \eqref{C22} hold. 

First notice that
\begin{equation}\label{Ex1B2V}
V^{n,l}_\dt=\overline{Q^{n,l}_\dt}^{H^1}=\{\bu\in H^1_0(\Omega^{n,l}_\dt):\nabla\cdot\bu=0\},
\end{equation}
and so the functions in $V^{n,l}_\dt$ are in $H_0^1$, defined on the common domain $\Omega^{n,l}_\dt$, and are divergence free.

To construct an operator $I^i_{\dt,l,n}$, for each fixed $i\in \{0,\dots,l\}$, we take a function $\bu \in V_\dt^{n+i}$ and associate to it a function $\tilde{\bu} \in V^{n,l}_\dt$,
which is constructed in three steps, as follows. First,  we consider the restriction of $\bu$ on $V^{n,l}_\dt$.
Such a function is not necessarily zero on the boundary $\partial\Omega^{n,l}_\dt$. 
This is why we ``move'' a bit inside $\Omega^{n,l}_\dt$, i.e., we construct a subset $\tilde\Omega^i_{\dt,l,n} \subset \Omega^{n,l}_\dt$, 
and consider $\bu$ outside of $\tilde\Omega^i_{\dt,l,n}$, i.e., on $\Omega^{n+i}_\dt\setminus \tilde\Omega^i_{\dt,l,n}$, extend it to the rest of $\tilde\Omega^i_{\dt,l,n}$ via an extension operator, to get a function $\bv$. 
By subtracting $\bv$ from $\bu$ we obtain a function which is zero on $\partial\Omega^{n,l}_\dt$, but is no longer divergence free inside $\Omega^{n,l}_\dt$.
We ``correct'' this by adding a function $\bf{w}$ which is such that $\tilde\bu := \bu - \bv + {\bf{w}}$ is divergence free.
By paying attention to details in this construction, as we explain below,  $\tilde{\bu}$ can be constructed such that it is in $V^{n,l}_\dt$, 
and is such that $\tilde\bu$ is a good approximation of $\bu$ in the sense of Property C2.

More precisely, fix $i$, $n$ and $l$, and recall that $h = l \dt$. 
Define $\Omega_{\tilde\gamma}$, which is contained in all the domains $\Omega^{n+i}_\dt$, as well as $\Omega^{n,l}_\dt$, by taking 
 $\tilde\gamma = {4\frac{c_L}{C_L}h}$. Consider the domain 
$$\tilde\Omega^i_{\dt,l,n} = \bdeta((n+i)\dt,\Omega_{\tilde\gamma}), \ \tilde\gamma = {4\frac{c_L}{C_L}h},$$
 which is
obtained by taking $\bdeta$ at time $t = (n+i)\dt$, but defined on the smallest domain $\Omega_{\tilde\gamma}$.
By Lemma \ref{LemGeom1} we have $\tilde\Omega^i_{\dt,l,n} \subset\subset\Omega^{n,l}_\dt\subset\subset\Omega^{n+i}_\dt$. 

Define $S_h=\Omega^{n+i}_\dt\setminus \tilde\Omega^i_{\dt,l,n}$. 
Because $\tilde\gamma = {4\frac{c_L}{C_L}h}$, and $\bdeta$ are Lipschitz continuous, it is easy to see $|S_h|\leq Cl\dt$.
Let $\bv=E_h\bu|_{{S}_h}$, where $E_h$ is an extension operator from domain ${S}_h$
onto the entire domain $\Omega^{n+i}_\dt$. 
By the well-known theorem for extension operators, see Thm 5.24 in \cite{ADA},  we have:
\begin{align}
&\|{\bv}\|_{W^{1,p}(\Omega_M)}\leq C\|{\bu}\|_{W^{1,p}({S}_h)}
\leq C\|{\bu}\|_{H^{1}(\Omega_M)}|{S}_h|^{(2-p)/p},\; p< 2,
\label{in1}
\\
&\|{\bv}\|_{H^1(\Omega_M)}\leq C\|{\bu}\|_{H^1(\Omega_M)}.
\label{in2}
\end{align}
Here, the first inequality in \eqref{in1} follows from \cite{ADA}. 
The second inequality in \eqref{in1} is a consequence of H\"{o}lder inequality.
We note that the constant above depends on the domain (Lipschitz domain), therefore, 
it depends on $i,n,l$. However, for small enough $h$, we argue that 
we can take the constant independent of $h$ by using the same partition of unity for all domains $S_h$,
since these domains are close to each other. 

By this construction,  the trace of $\bu-\bv$ is zero on $\partial\Omega^{n,l}_\dt$. However, the difference $\bu-\bv$ is not divergence free.
To rectify this,  we introduce a function 
${\bf w}\in W^{1,p}_0(\tilde\Omega^i_{\dt,l,n})$, $p\leq 2$, which is a solution to the following problem: 
\begin{equation}\label{DivCorrector}
\nabla\cdot{\bf w}=\nabla\cdot{\bf v}\quad {\rm in}\; \tilde\Omega^i_{\dt,l,n}.
\end{equation}
Solutions of this problem are given by construction in \cite{galdi2011introduction}, Thm. III.3.1,  
provided that the right hand-side satisfies the
following compatibility condition:
$$
\int_{{\tilde\Omega^i_{\dt,l,n}}}\nabla\cdot{ \bf v}
=\int_{\partial \tilde\Omega^i_{\dt,l,n}}{ \bv}\cdot\bn
=\int_{\partial\tilde\Omega^i_{\dt,l,n}}{ \bu}\cdot\bn
=\int_{\tilde\Omega^i_{\dt,l,n}}\nabla\cdot\bu
=0.
$$

Using this ``correction'' $\bf{w}$ we can now define a function ${\tilde\bu}:={\bu}-{\bv}+{\bf w}$, 
and the mapping $I^i_{\dt,n,l}$ so that
$$
I^i_{\dt,n,l} \bu = \tilde\bu,
$$
which satisfies the following properties:
(1)  ${\tilde\bu}\in H^1_0(\tilde\Omega^i_{\dt,l,n})\subset V^{n,l}_\dt $, and
(2) $\nabla\cdot {\tilde\bu}=0$.

We now show that inequalities \eqref{C21} and \eqref{C22} are satisfied. 
To show \eqref{C21} we need to prove the existence of a universal constant $C > 0$ such that
$$
\|\tilde\bu\|_{V^{n,l}_\dt} \le C \|\bu \|_{V^{n+i}_\dt}.
$$
First, from \cite{galdi2011introduction}, Thm. III.3.1 we have 
\begin{equation}\label{in3}
\|{\bf{w}}\|_{W^{1,p}_0(\tilde\Omega^i_{\dt,l,n})}\leq C(t,h)\|{\bv}\|_{W^{1,p}(\tilde\Omega^i_{\dt,l,n})}, \ p\leq 2.
\end{equation}
However, the constant $C$ depends on $t$ and $h$. We claim that there exists an absolute constant $C>0$ such that $C(t,h) < C$,
and \eqref{C21} is satisfied. Suppose that this is not true, i.e. there exists a sequence $(t_n,h_n)\to (t,h)$ such that $C(t_n,h_n)\to\infty$. 
The existence of a sequence   $(t_n,h_n)\to (t,h)$ such that $C(t_n,h_n)\to\infty$ contradicts
Thm. III.3.1 from \cite{galdi2011introduction} and the fact that $\bdeta$, which is the limiting domain,  is a Lipschitz function with 
a finite constant $C(t,h)<\infty$. 

Now we have:
$$
\|\tilde\bu\|_{V^{n,l}_\dt}  \le \|\bu\|_{V^{n,l}_\dt} + \|{\bf{w}}\|_{V^{n,l}_\dt} + \|\bv\|_{V^{n,l}_\dt} 
\le \|\bu\|_{V^{n+i}_\dt} + C \| \bv \|_{V^{n+i}_\dt}\le C \| \bu \|_{V^{n+i}_\dt},
$$
where the last inequality follows from \eqref{in2}, by recalling the $\bv$ and $\bu$ are extended by zero outside $\Omega^{n+i}_\dt$.
Therefore, \eqref{C21} holds.

To show that \eqref{C22} holds, we need to show that there exists a monotonically increasing function $g$ such that $g \to 0$ as $h\to 0$, and
\begin{equation}\label{C22_recall}
\|{\bu}-{\tilde\bu}\|_{L^2(\Omega_M)} \leq g(l\dt) \|\bu \|_{V^{n+i}_\dt}.
\end{equation}
We calculate:
\begin{equation}\label{NsDiff}
\|{\bu}-{\tilde\bu}\|_{L^2(\Omega_M)}\leq \|{\bf v}\|_{L^2(\Omega_M)}+\|{\bf w}\|_{L^2(\Omega_M)}
\leq  C(\|{\bf v}\|_{W^{1,p}}+\|{\bf w}\|_{W^{1,p}})
\leq Ch^{(2-p)/p} \|\bu \|_{H^1(\Omega_M)},
\end{equation}
with $2d/(d+2)\leq p<2$ (Sobolev embedding of $L^2$ into $W^{1,p}$), where we have used \eqref{in3}, \eqref{in1} and \eqref{in2} in the last inequality.
We conclude that \eqref{C22_recall} holds by recalling that $\bu$ is zero outside $\Omega^{n+i}_\dt$.

\begin{remark}[On the regularity assumption on $\bdeta$]{\rm{
Property (C2) holds even when $\bdeta$ is not necessarily Lipschitz. 
For example, Property (C2) can be verified for a fluid domain which is only a sub-graph of a H\"{o}lder continuous function. 
This shows that our framework is rather flexible, in the sense that it can be applied to scenarios requiring significantly less regularity of $\bdeta$
than the current results in literature.}}
\end{remark}

\noindent
\textbf{Property C3: Uniform Ehrling Property.}
We follow the main steps of the proof of the classical Ehrling's lemma and use the method of contradiction. 
To simplify notation we introduce only one index $i$ to denote
the sequences of function spaces $H^{n,l}_\dt$, $V^{n,l}_\dt$, $Q^{n,l}_\dt$, which will now be denoted by $H^i$, $V^i$, $Q^i$. 

Suppose that the statement of the Uniform Ehrling Property is false, i.e., that there exists a $\delta> 0$ such that for every $k > 0$, 
there exists a sub-sequence, $\bu_{i_k} \in {V^{i_k}}$, indexed by $i_k$, such that
\begin{equation}\label{Ex1Ehrl}
\|\bu_{i_k}\|_{H^{i_k}}>\delta\|\bu_{i_k}\|_{V^{i_k}}+k\|\bu_{i_k}\|_{(Q^{i_k})'}.
\end{equation}
Again, for simplicity, in the rest of the proof we will write $\bu_k$ instead $\bu_{i_k}$. 

Without the loss of generality we assume $\|\bu_k\|_{H^k}=1$ (otherwise we can divide \eqref{Ex1Ehrl} by $\|\bu_k\|_{H^k}$). 
Inequality \eqref{Ex1Ehrl} implies that 
\begin{equation}\label{step1}
\|\bu_k\|_{V^k}\leq \frac{1}{\delta},\ {\rm and} \  \|\bu_k\|_{(Q^k)'}\to 0, \  {\rm as} \ k \to \infty.
\end{equation}
We will show that \eqref{step1} contradicts $\|\bu_k\|_{H^k}=1$.

For this purpose, we first notice that due to the uniform boundedness $\|\bu_k\|_{V^k}\leq \frac{1}{\delta}$,
and the fact that $V_k\subset V$ for each $k$, and $V\subset\subset H$, 
there exists a subsequence, still denoted by $\bu_k$, such that $\bu_k$ converges 
to some $\bu$ in $H$, i.e., $\bu_k\to \bu$ in $H$.

On the other hand, each of those spaces $V_k$, has an associated function $\bdeta^k$, which defines the domain $\Omega^k=\bdeta^k(\Omega_{\gamma_k})$,
on which the spaces $H^k$, $V^k$, $Q^k$  are defined. 
Since sequence $\bdeta_k$ is bounded in $W^{1,\infty}$, and $\gamma_k$ in $\R$, there exists a subsequence (again denoted by $\bdeta_k$, and by $\gamma_k$) such that 
$\bdeta_k$ converges to some $\bdeta$, i.e., $\bdeta_k\to\bdeta$ in $C(\overline{\Omega_M})$, and $\gamma_k\to\gamma$, for some $\gamma\geq 0$. 

Let $\Omega_\eta=\bdeta(\Omega_{\gamma})$.
We will show that $\bu = 0$ on $\Omega_\eta$, and on $\Omega_M\setminus\Omega_\eta$, which implies $\bu = 0$ in $\Omega_M$,
which contradicts the assumption that $\|\bu_k\|_{H^k} = 1$.

Let $K$ be a compact subset of $\Omega_\eta$, and $\bq\in C_c^{\infty}(\Omega_M)$ such that supp $\bq\subset K$. 
Then, because of the uniform convergence of $\bdeta_k$, there exists a $k_0>0$ such that 
$K$ is contained compactly in infinitely many sets $\Omega^k$, i.e., $K\subset\subset\Omega^k$, $k\geq k_0$. 
Therefore, $\bq\in Q^k$, for all $k\geq k_0$, and hence,
by using the second statement in \eqref{step1}, we obtain:
\begin{equation}\label{ehrling_zero}
\int_{\Omega_\eta}\bu\cdot\bq=\lim_{k\to\infty}\int_{\Omega^k}\bu_k\cdot\bq=0.
\end{equation}

In other words, \eqref{ehrling_zero} implies that $\bu$ is orthogonal to the divergence free fields $\bq$
(recall both $\bu$ and $\bq$ are zero outside $\Omega_\eta$). 
But, since $\bu$ is divergence free itself, it follows from the Helmholtz decomposition
that $\bu$ is both in the orthogonal complement of divergence free fields, and 
in the space of divergence free fields, and so $\bu$ must be equal to zero. 

Let us now take a compact set $S\subset\Omega_M\setminus\overline{\Omega_\eta}$. 
Again, because of the uniform convergence of $\bdeta_k$, there exists a $k_1>0$ such that $\bu_k=0$ on $S$, for all $k\geq k_1$. 
Therefore $\bu=0$ on $S$. Since, $S$ was an arbitrary compact subset of $\Omega_M\setminus\overline{\Omega_\eta}$ we conclude that $\bu=0$ on $\Omega_M\setminus\Omega_\eta$.

Therefore $\bu=0$ in $\Omega_M$, which contradicts the assumption $\|\bu_k\|_{H^k}=1$, and concludes the proof of the uniform Ehrling property.
\qed

\subsection{Fluid-structure interaction with the no-slip condition}\label{sec:FSI}
We show here how Theorem~\ref{Compactness} can be applied to prove existence of a weak solution to a fluid-structure interaction 
problem between an incompressible, viscous fluid and an elastic structure modeled by the linear Koiter shell equations.
Existence of weak solutions to this FSI problem was first studied in \cite{BorSun}. 
We show here how the compactness theorem discussed in the present manuscript
can be directly applied to the problem studied in \cite{BorSun}. 
We present here the definition of the FSI problem, and the main energy estimates obtained in \cite{BorSun}, 
which are important for the compactness proof.
Then we verify assumptions (A)-(C),  and state the compactness results of Theorem~\ref{Compactness}
as it applies to the approximate sequences of functions, constructed in \cite{BorSun}, approximating weak solutions of the underlying FSI problem.

\subsubsection{Problem definition}

\begin{figure}[htp]
\center
\includegraphics[scale=0.27]{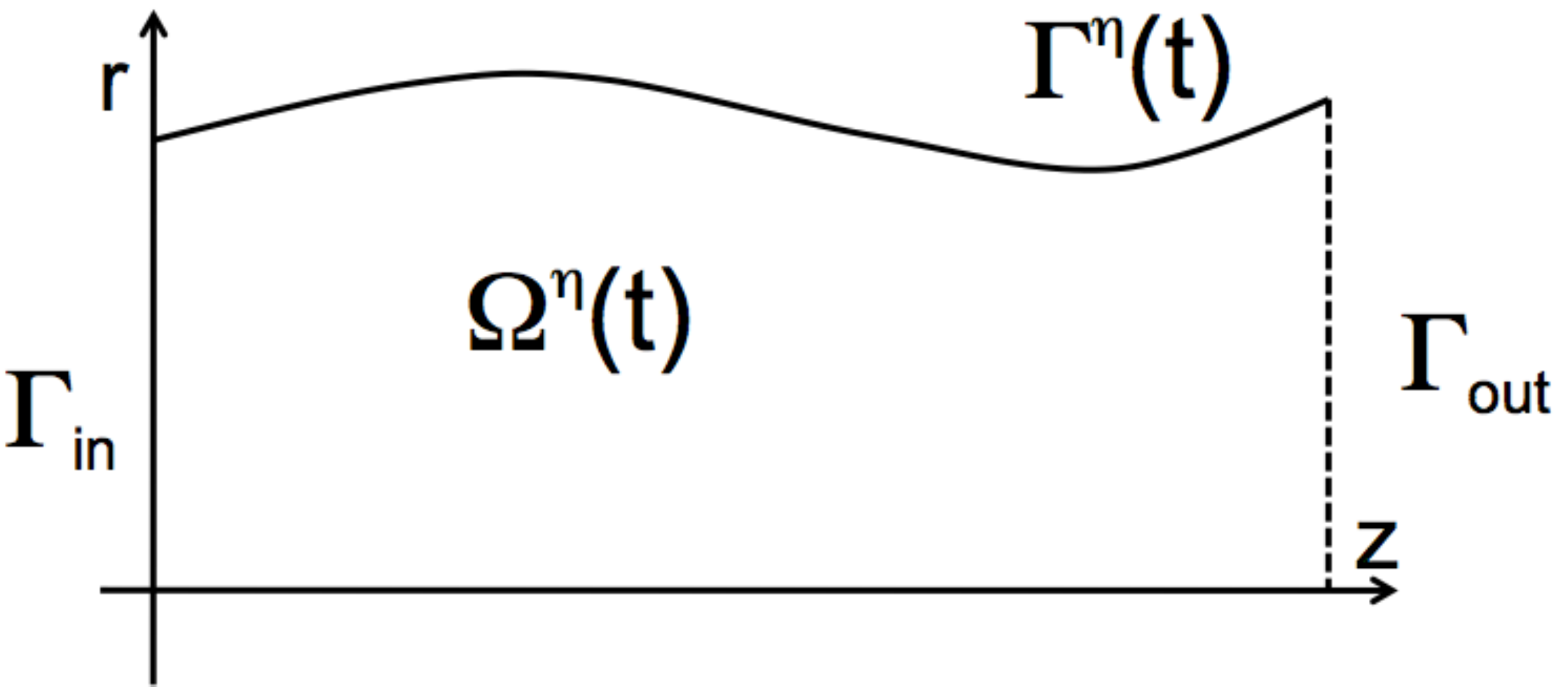}
\caption{A sketch of the FSI domain. }
\label{domain}
\end{figure}

We consider 
 a nonlinear moving boundary problem modeling the flow of an incompressible, viscous fluid through 
an elastic cylinder. The fluid flow is modeled by the Navier-Stokes equations for an incompressible, viscous fluid:
\begin{equation}\label{NS}
\left.\begin{array}{rcl}
\rho_f\big ({\partial_t{\bu}}+({\bu}\cdot\nabla){\bu}\big)&=&\nabla\cdot \boldsymbol\sigma    \\
\nabla\cdot{\bu}&=&0 
\end{array}
\right\}\ \textrm{in}\; \Omega^{\bdeta}(t),\ t\in (0,T), 
\end{equation}
where $\boldsymbol\sigma$ is the Cauchy stress tensor, $\rho_f$ is the fluid density, and $\bu=\bu(\boldsymbol{x},t)$ is the fluid velocity. 
The equations are defined on the domain $\Omega^{\bdeta}(t)$ which depends on time, and also on the position of the
elastic structure whose displacement is denoted by $\bdeta$. 
Displacement $\bdeta(\boldsymbol{x},t)$ is determined from the elastodynamics of the structure modeled by, e.g.,
the elastic shell equations, which are defined on a reference configuration $\Gamma$ (in Lagrangian formulation), and can be written in general form as:
\begin{eqnarray}
\rho_{S} h  \partial_{tt} \bdeta + {\cal L}_e \bdeta = {\bf f}, \quad & {\rm on} \  \Gamma,\ t\in (0,T),
\label{Koiter}
\end{eqnarray}
where $\rho_S$ is the structure density, $h$ the elastic shell thickness, ${\bf f}$ is linear force density acting on the shell, 
and $\bdeta$ is  shell displacement from the reference configuration $\Gamma$. See Figure~\ref{domain}.
Operator ${\cal L}_e$ is a linear, continuous, positive-definite, coercive operator on a function space $\chi$.\
For example, for the linearly elastic Koiter shell considered in \cite{BorSun}, 
where only radial component of displacement was considered to be different from zero,
 the space $\chi = H_0^2(\Gamma)$.
Equations \eqref{NS} and \eqref{Koiter} are supplemented with initial and boundary conditions.

The fluid and structure interact through a two-way coupling, which takes place at the fluid-structure interface, which we denote by 
 $\Gamma^{\bdeta}(t)$.
The fluid flow influences the motion of the structure through traction forces 
exerted by the fluid onto the structure at $\Gamma^{\bdeta}(t)$, while the motion of
the structure influences the fluid through the kinetic and potential energy associated with the elastodynamics of the shell. 
Additionally, at the fluid-structure interface $\Gamma^{\bdeta}(t)$ the velocities of the fluid and structure are coupled via a kinematic coupling condition
such as, e.g., the no-slip condition. In the next section we consider the kinematic coupling condition described by the Navier slip condition. 

Thus, the following two coupling conditions evaluated at the fluid-structure interface $\Gamma^{\bdeta}(t)$ describe the fluid-structure coupling:
\begin{align}
\partial_t\bdeta &= \bu|_{\Gamma^{\bdeta}(t)},
\label{kinematic}\\
\rho_{S} h  \partial_{tt} \bdeta + {\cal L}_e \bdeta &= \boldsymbol\sigma\boldsymbol{n}|_{\Gamma^{\bdeta}(t)}.
\label{dynamic}
\end{align}
Here $\boldsymbol{n}$ is the (outward) unit normal to the fluid-structure interface $\Gamma^{\bdeta}(t)$, and
notation $\bu|_{\Gamma^{\bdeta}(t)}$ means that the fluid velocity $\bu$ is evaluated at the current location of the interface $\Gamma^{\bdeta}(t)$.
Condition \eqref{kinematic} is called the kinematic, and \eqref{dynamic} is called the dynamic coupling condition.
Condition \eqref{dynamic} assumes zero traction from the external environment. Otherwise, the right hand side of \eqref{dynamic}
is given by the jump between the traction exerted by the fluid onto the thin structure and the exterior. 

Equations \eqref{NS}, \eqref{Koiter}, \eqref{kinematic}, and \eqref{dynamic} define a nonlinear moving-boundary problem for 
the unknown functions $\bu$ and $\bdeta$.
This problem can be written as a first-order system \eqref{ProblemA} by introducing a new variable $\boldsymbol{v} = \partial_t \bdeta$,
which denotes the structure velocity.
In that case problem \eqref{NS} - \eqref{dynamic} reads:
\begin{align}
&\qquad \rho_f {\partial_t{\bu}} = \rho_f  ({\bu} \cdot\nabla){\bf u} + \nabla\cdot \boldsymbol\sigma ,\ {\rm with \ contraint} \   \nabla\cdot{\bu}=0
\quad {\rm in} \ \Omega^{\bdeta} (t),\  t \in (0,T),
\nonumber \\
&\left.
\begin{array}{rcl}
\partial_t\bdeta &=& \bu|_{\Gamma^{\bdeta}(t)}
\\
\rho_{S} h  \partial_{t} \boldsymbol{v} &=& -  {\cal L}_e \bdeta + \boldsymbol\sigma\boldsymbol{n}|_{\Gamma^{\bdeta}(t)}\\
\partial_t \bdeta &=&  \boldsymbol{v}
\end{array}
\right\}\ {\rm in} \ \Gamma, \ t\in(0,T).
\label{Monolithic}
\end{align}
Recall that the problem studied in \cite{BorSun} considered  only the radial component of structure displacement to be different from zero,
denoting it by $\eta$. The corresponding structure velocity $\bv$ becomes $v{\bf e}_r$.

To write a weak formulation of problem \eqref{Monolithic} as considered in \cite{BorSun}, the following function space describes
the behavior of the fluid velocity:
$$
V_F(\Omega^{\eta})=\{\bu\in H^1(\Omega^\eta):\nabla\cdot\bu=0,\; \bu=v{\bf e}_r\; {\rm on}\;\Gamma^\eta,\; \bu\cdot{\bf e}_r=0\; {\rm on}\;\partial\Omega^\eta\setminus\Gamma^\eta\}.
$$
The test functions $\bq$ for the fluid velocity  and the test functions $\psi$ for the structure are required to belong to
$$
Q^\eta=\{(\bq,\psi)\in \big ( V_F(\Omega^{\eta})\cap H^4(\Omega^{\eta})\big )\times H^{2}(\Gamma):\bq_{|\Gamma^{\eta}}=\psi\}.
$$
The weak formulation is obtained by multiplying the momentum equation in the Navier-Stokes 
equations by $\bq$, and the structure equation \eqref{Koiter} by $\psi$
and integrating by parts. After taking into account the divergence-free condition and
the dynamic coupling condition \eqref{dynamic}, it was shown in \cite{BorSun} that the following weak formulation holds:
find $(\bu,\eta)$ such that 
\begin{align}
&\displaystyle{\rho_f \int_{\Omega^{\eta}(t)}\big ( \partial_t{\bu} \cdot {\bq} + \big ({\bu}\cdot\nabla){\bu}\cdot{\bq}\big )
 +  2\mu\int_{\Omega^{\eta}(t)}{\bf D}({\bu}):{\bf D}({\bq})}
\displaystyle{ +\rho_S h \int_{\Gamma} \partial_t v \psi } 
\displaystyle{+  a_e(\eta,\psi)  = \langle {\bf F},{\bq}\rangle},\; 
\nonumber \\
&\int_{\Gamma} \partial_t\eta   \psi = \int_{\Gamma} {v} \psi,\quad
\forall (\bq,\psi) \in Q^\eta,
\label{FSIweak}
\end{align}
with the prescribed initial conditions for the fluid velocity, the structure position and structure velocity.
Here, the source term functional ${\bf F}$ is used to account for the boundary 
conditions at the inlet and outlet of the fluid domain, and $\langle \phantom{x},\phantom{x}\rangle$ denotes the duality pairing 
in the corresponding spaces.
The term $a_e(\eta,\psi)$ is the bi-linear form that results from operator ${\cal{L}}_e$ after integration by parts.

To make sense of the acceleration term involving the integral of $\partial_t \bu$ over $\Omega^{\eta}(t)$,  it is convenient
to introduce an Arbitrary Lagrangian-Eulerian (ALE) mapping which maps the moving domain onto a fixed reference domain $\Omega$.
To follow the approach in \cite{BorSun}, one can introduce an explicit 
formula for the ALE mapping:
\begin{equation}
A_{\eta}(t):\Omega\rightarrow\Omega^{\eta}(t),\quad
A_{\eta}(t)(\tilde{z},\tilde{r}):=\left (\begin{array}{c}\tilde{z}\\(R+\eta(t,\tilde z))\tilde{r}\end{array}\right ),\quad (\tilde{z},\tilde{r})\in\Omega,
\label{RefTrans}
\end{equation} 
where $(\tilde z,\tilde r)$ denote the coordinates in the reference domain $\Omega=(0,L)\times (0,1)$.
The mapping $A_\eta(t)$ is a bijection, 
and its Jacobian $J^\eta$ and the related ALE domain velocity ${\bf w}^{\eta}$ are given by
\begin{equation}\label{ALE_Jacobian}
J^\eta = |{\rm det} \nabla A_\eta(t)| = |R + \eta(t,\tilde z)|,\quad
\displaystyle{{\bf w}^{\eta}=\partial_t\eta{\tilde r}{\bf e}_r}.
\end{equation}
Note that the transformed gradient, $\nabla^{\eta}$ is given by
$
\nabla^{\eta} = \left (
\displaystyle{\partial_{\tilde{z}}-\tilde{r}\frac{\partial_z\eta}{R+\eta}\partial_{\tilde{r}}}
,
\displaystyle{
\frac{1}{R+\eta}\partial_{\tilde{r}}}
\right )^T,
$
and so
$
\nabla^{\eta}{\bf v}=\nabla{\bf v}(\nabla A_{\eta})^{-1}.
$
We denote
$ \sigma^{\eta}=-p^{\eta}{\bf I}+2\mu{\bf D}^{\eta}({\bf u}^{\eta})$, where 
${\bf D}^{\eta}({\bf u}^{\eta})=\frac 1 2(\nabla^{\eta}{\bf u}^{\eta}+(\nabla^{\eta})^{\tau}{\bf u}^{\eta}).$

The weak form \eqref{FSIweak} can now be written in ALE form on $\Omega$ (with the symmetrized advection term), as follows \cite{BorSun}:
\begin{gather}
\displaystyle{
\rho_F\int_{\Omega}J^{\eta}\partial_t {\bf u}\cdot{\bf q} 
+ \frac{\rho_F}{2}\int_{\Omega}J^{\eta}\Big ((({\bf u}-{\bf w}^{\eta})\cdot\nabla^{\eta}){\bf u}\cdot {\bf q}-(({\bf u}-{\bf w}^{\eta})\cdot\nabla^{\eta}){\bf q}\cdot {\bf u}}
-(\nabla^{\eta}\cdot{\bf w}^{\eta}){\bf q}\cdot {\bf u}\Big)
\nonumber
\\ 
\displaystyle{
+2\mu\int_{\Omega} J^\eta {\bf D}^{\eta}({\bf u}):{\bf D}^{\eta}({\bf q})
+\rho_sh\int_{\Gamma} \partial_t v \psi
+a_S(\eta,\psi)
= \langle {\bf F},{\bq}\rangle},\\
\int_{\Gamma} \partial_t\eta \cdot  \psi = \int_{\Gamma} {v}\cdot \psi,\quad
\forall (\bq,\psi) \in Q^\eta,
\label{FSIweak2}
\end{gather}
where $Q^\eta$ is the test space defined for each $t$ by the ALE mapping $A_\eta$.
The problem studied in \cite{BorSun}  is driven by the inlet and outlet dynamic pressure data $P_{in/out}$, 
prescribed at $\Gamma_{in}$ and $\Gamma_{out}$ respectively, and so
\begin{equation}\label{in_out_data}
\langle {\bf F},{\bq}\rangle := 
R\big (
P_{in}(t) \int_{\Gamma_{in}}(q_z)|_{\Gamma_{in}}-P_{out}(t) \int_{\Gamma_{out}} (q_z)|_{\Gamma_{out}}
\big ).
\end{equation}

The work in \cite{BorSun} addressed the existence of a weak solution to problem \eqref{Monolithic},
providing the first constructive existence proof for FSI problems. It involved 
a semi-discretization approach, similar to the one described in Section~\ref{sec:introduction}.

\subsubsection{The semi-discretized problem}\label{sec:FSIsemi}
We consider the time interval $(0,T)$  and divide it into $N$ sub-intervals, each of width $\dt = T/N$.
The following 
sequences of approximations 
$\{{\bf u}^n_\dt\}$, $\{\eta^n_\dt\}$, and $\{ v^{n+i/2}_\dt \}$, $i = 0,1$ are constructed using the operator splitting strategy
described in \cite{BorSun}.
The resulting functions $\{{\bf u}^n_\dt\}$, $\{\eta^n_\dt\}$, and $\{ v^{n+i/2}_\dt \}$, $i = 0,1$
satisfy the following semi-discretized version of the weak form \eqref{FSIweak2},
which can be obtained after taking into account  the precise form of the fluid domain velocity ${\bf w}^{\eta}$ and of $\nabla^\eta \cdot {\bf w}^{\eta}$ 
as specified above:
\begin{equation}
\begin{array}{c}
\displaystyle{\rho_f\int_{\Omega}(R+\eta^{n}_\dt) \frac{{\bf u}^{n+1}_\dt-{\bf u}^{n}_\dt}{\Delta t}\cdot{\bf q}+\frac{\rho_f}{2} \int_{\Omega}{v^{n+1/2}_\dt}{\bf u}^{n+1}_\dt\cdot{\bf q}}
+2\mu\int_{\Omega}(R+\eta^{n}_\dt){\bf D}^{\eta^{n}_\dt}({\bf u}^{n+1}_\dt):{\bf D}^{\eta^{n}_\dt}({\bf q})
\\ \\
\displaystyle{+\frac{\varrho_f}{2} \int_{\Omega} {(R+\eta^{n}_\dt)}
\big ( 
\left[({\bf u}^{n}_\dt-v^{n+1/2}_\dt r{\bf e}_r)\cdot\nabla^{\eta^{n}_\dt}\right ]{\bf u}^{n+1}_\dt\cdot{\bf q}- \left [({\bf u}^{n}_\dt-v^{n+1/2}_\dt r{\bf e}_r)\cdot\nabla^{\eta^{n}_\dt} \right] {\bf q}\cdot{\bf u}^{n+1}_\dt}\big )
\\ \\
\displaystyle{+\rho_sh\int_\Gamma \frac{v^{n+1}_\dt-v^{n}_\dt}{\Delta t}\psi
+a_e(\eta^{n+1}_\dt,\psi)}
= F^n(\bq),
\\ \\
\displaystyle{\int_\Gamma \frac{\eta^{n+1}_\dt - \eta^n_\dt}{\dt}\psi = \int_\Gamma v^{n+1/2}_\dt \psi,}\quad
\forall (\bq,\psi) \in Q^{\eta^n_\dt},
\end{array}
\label{ProbUk}
\end{equation}
${\rm with}\ \nabla^{\eta^{n}_\dt}\cdot{\bf u}^{n+1}_\dt=0,\quad {\bf u}^{n+1}_\dt{|_\Gamma}=v^{n+1}_\dt {\bf e}_r,$
and
$
F^n(\bq) \displaystyle{:=R\big (
P^{n}_{in}\int_{\Gamma_{in}} (q_z)|_{\Gamma_{in}}-P^{n}_{out}\int_{\Gamma_{out}} (q_z)|_{\Gamma_{out}}\big ),}
$
where $P^{n}_{in/out}$ denote the piecewise constant approximations of $P_{in/out}(t)$.

Notice that  ${\bf u}^{n+1}_\dt$ is defined via the ALE mapping associated with the ``previous'' (known) domain $\Omega^{\eta^n_\dt}$, i.e.,
${\bf u}^{n+1}_\dt=\bu^{n+1}_\dt\circ A^{-1}_{\eta^{n}_\dt}$. 
Thus, approximations ${\bf u}^{n+1}_\dt$  with the superscript $n+1$ correspond to the functions defined on 
the "previous" fluid domain $\Omega^{\eta^n}_\dt$, not the ``current'' one $\Omega^{\eta^{n+1}_\dt}$. Notice that this was our choice, which was based on 
typical approaches used in numerical simulations of this class of FSI problems. As $\dt \to 0$, 
the solution/limit will be independent of this choice.

It was shown in \cite{BorSun} that the following energy estimate holds.
\begin{proposition}
For each fixed $\Delta t > 0$, functions $\{{\bf u}^n_\dt\}$, $\{\eta^n_\dt \}$, and $\{ v^{n+i/2}_\dt \}$, $i = 0,1$ satisfy the following discrete energy inequality:
\begin{equation}
\begin{array}{c}
E_\dt^{n+1}+\displaystyle{\frac{\rho_f}{2}\int_{\Omega}(R+\eta^n_\dt)|{\bf u}^{n+1}_\dt-{\bf u}^n_\dt|^2+\frac{\rho_s h}{2}\|v^{n+1}_\dt-v^{n+\frac 1 2}_\dt\|^2_{L^2(\Gamma)}}\\
+D^{n+1}_N\leq E_\dt^{n}+C\Delta t((P_{in}^n)^2+(P_{out}^n)^2),
\end{array}
\label{DEE}
\end{equation}
where the total (kinetic and elastic) energy $E_N^{n+1}$ and viscous dissipation $D^{n+1}_\dt$ are defined by
 \begin{equation}
\displaystyle{E_\dt^{n+1}=\frac 1 2\Big (\rho_f\int_{\Omega}(R+\eta^{n})|{\bf u}^{n+1}_\dt|^2+
\rho_s h\|v^{n+1}_\dt\|^2_{L^2(\Gamma)} + \| \eta^{n+1}_\dt \|^2_{{\cal{L}}_e}}\Big ),
\label{kenergija}
\end{equation}
\begin{equation}
\displaystyle{D_\dt^{n+1}=\Delta t \mu\int_{\Omega}(R+\eta^n_\dt)|D^{\eta^n_\dt}({\bf u}_\dt^{n+1})|^2}.
\label{kdisipacija}
\end{equation}
The constant $C$ in the above estimate depends only on the parameters in the problem, and not on $\Delta t$ (or $N$).
\end{proposition} 
Here $\| \cdot \|_{{\cal{L}}_e}$ denotes the norm associated with the bilinear functional $a_e$, i.e., with the 
elasticity operator ${{\cal{L}}}_e$. For a linearly elastic cylindrical Koiter shell allowing displacement only in the radial direction 
$\| \cdot \|_{{\cal{L}}_e}$ is equivalent to the ${H^2(\Gamma)}$-norm. 

As a consequence, the following uniform energy bounds hold \cite{BorSun}.
\begin{lemma}\label{EnergyEstimates}
There exists a constant $C>0$ independent of $\Delta t$ (and $N$), 
such that the following estimates hold:
\begin{enumerate}
\item
$E_\dt^{n+1}\leq C$, for all $ n = 0,...,N-1$, where $\dt = T/N$,
\item $ \sum_{j=1}^N D_\dt ^j\leq C,$
\item
$\displaystyle{\sum_{n=0}^{N-1}\left(\int_{\Omega}(R+\eta^n_\dt)|{\bf u}^{n+1}_\dt-{\bf u}^n_\dt |^2+\|v^{n+1}_\dt-v^{n+\frac 1 2}_\dt \|^2_{L^2(\Gamma)}
 +\|v^{n+\frac 1 2}_\dt-v^{n}\|^2_{L^2(\Gamma)}\right)\leq C}$,
\item
$\displaystyle{\sum_{n=0}^{N-1}\| \eta^{n+1}_\dt - \eta^n_\dt \|_{{\cal{L}}_e} \le C.}$
\end{enumerate}
The constant $C = E_0 + \tilde{C} \left(\|P_{in}\|_{L^2(0,T)}^2 + \|P_{out}\|_{L^2(0,T)}^2\right)$, where $\tilde{C}$ is the constant from 
\eqref{DEE},
which depends only on the parameters in the problem.
\end{lemma}

For each $\dt > 0$, approximations $({\bf u}^n_\dt,v^n_\dt,\eta^n_\dt), n = 1, \dots, N,$ define $({\bf u}_\dt,v_\dt,\eta_\dt)$, which is piecewise constant
on each sub-interval $(n \dt, (n+1) \dt) \subset (0,T)$, as in \eqref{FuncDef2}. As $\dt \to 0$, the goal in \cite{BorSun} was to show that 
there exists a sub-sequence of $\{({\bf u}_\dt,v_\dt,\eta_\dt)\}$ which converges
to a weak solution of problem \eqref{Monolithic}. A crucial component in the proof is to show that $\{({\bf u}_\dt,v_\dt)\}$ is precompact in $L^2(0,T;H)$, 
where $H$ is a Hilbert space to be specified below. We show here how Theorem~\ref{Compactness} can be used to show precompactness of 
$\{({\bf u}_\dt,v_\dt)\}$ in $L^2(0,T;H)$.

\subsubsection{Compactness result}

To show the compactness result, in the rest of this section we will work in the physical space, namely, on moving fluid domains $\Omega^\eta(t)$, instead of the fixed, reference 
domain $\Omega$. 
To state the weak formulation \eqref{ProbUk} of the semi-discretized problem on moving domains, we introduce the following spaces:
\begin{equation}\label{Qn_dt}
Q^n_\dt=\{(\bq,\psi)\in \big ( V_F(\Omega^{\eta^{n}_\dt})\cap H^4(\Omega^{\eta^{n}_\dt}\big )\times H^{2}(\Gamma):\bq_{|\Gamma^{\eta_\dt^{n}}}=\psi\},
\end{equation}
which will serve as the test spaces for each fixed $\dt$, and the solution spaces  
\begin{equation}\label{Vn_dt}
V^{n+1}_\dt=\{(\bu,v)\in V_F(\Omega^{\eta^{n}_\dt})\times H^{1/2}(\Gamma):\bu_{|\Gamma^{\eta_\dt^{n}}}=v\}.
\end{equation}

The weak formulation \eqref{ProbUk} is then transformed to:
\begin{equation}
\begin{array}{c}
\displaystyle{\varrho_f\int_{\Omega^{\eta^{n}_\dt}}\frac{{\bf u}^{n+1}_\dt-\tilde{\bf u}^{n}_\dt}{\Delta t}\cdot{\bf q}
+\frac{\rho_f}{2} \int_{\Omega^{\eta^{n}_\dt}}\frac{v^{n+1/2}_\dt}{R+\eta^{n}_\dt}{\bf u}^{n+1}_\dt\cdot{\bf q}} 
+2\mu\int_{\Omega^{\eta^{n}_\dt}}{\bf D}({\bf u}^{n+1}_\dt):{\bf D}({\bf q})
\\ \\
\displaystyle{+\varrho_f\int_{\Omega^{\eta^{n}_\dt}}\frac{1}{2}\left[(\tilde{\bf u}^{n}_\dt-\frac{v^{n+1/2}_\dt r}{R+\eta^{n}_\dt}{\bf e}_r)\cdot\nabla\right .
{\bf u}^{n+1}_\dt\cdot{\bf q}-\left . (\tilde{\bf u}^{n}_\dt-\frac{v^{n+1/2}_\dt}{R+\eta^{n}_\dt}{\bf e}_r)\cdot\nabla \right] {\bf q}\cdot {\bf u}^{n+1}_\dt}
\\ \\
\displaystyle{+\rho_sh\int_\Gamma \frac{v^{n+1}_\dt-v^{n}_\dt}{\Delta t}\psi+a_e(\eta^{n+1}_\dt,\psi)=F^{n}({\bf q}),}
\\ \\
\displaystyle{\int_\Gamma \frac{\eta^{n+1}_\dt - \eta^n_\dt}{\dt}\psi = \int_\Gamma v^{n+1/2}_\dt \psi,}\quad
\forall (\bq,\psi) \in Q^{n}_\dt,
\end{array}
\label{Eulerian}
\end{equation}
${\rm with}\ \nabla\cdot {\bf u}^{n+1}_\dt=0,\quad {\bf u}^{n+1}_\dt{|_{\Gamma^{\eta^{n}_\dt}}}=v^{n+1}_\dt {\bf e}_r$.
Here:
\begin{equation}\label{tilde_u}
{\bf \tilde{u}}^{n}_\dt ={\bf u}^{n}_\dt \circ A_{\eta^{n}_\dt }^{-1},
\end{equation}
which is a function that takes  points from $\Omega^{\eta^n_\dt}$, maps them into $\Omega$ via  $A_{\eta^{n}_\dt}^{-1}$
and assigns to those points the values
of ${\bf u}^{n}_\dt$ (which originally came from $\Omega^{\eta^{n-1}_\dt}$).

To simplify notation, without the loss of generality, in the rest of this section we will be taking all the physical constants to be equal~$1$. 

To study compactness of the corresponding sequences, 
it is useful to work with a domain $\Omega^M$ which contains all the fluid domains $\Omega^{\eta^n}$. 
 Similarly, to obtain the desired estimates for the time shifts $\tau_h {\bf u}_\dt$, it is useful to work on a common domain
 that contains all the fluid domains determined by the time shifts from $t$ to $t +h$.
The existence of such domains is guaranteed by the following lemma.

\begin{lemma}\label{DoljeGore}
There exist smooth functions $m(z)$ and $M(z)$, $z \in \overline\Gamma = [0,L]$, such that
the sequence $\{\eta^n_\dt(z)\}$ corresponding to the semi-discretization of $(0,T)$
can be uniformly bounded as follows:
\begin{itemize}
\item $m(z)\leq \eta^{n}_\dt (z)\leq M(z),\; z\in [0,L],\; n=0,\dots,N,\;N\in\N$, $\dt = T/N$, and
\item $M(z)-m(z)\leq C(T),\; z\in [0,L],$
where $C(T)\to 0$, $T\to 0$. 
\end{itemize}
Moreover, for each fixed $N$, $l\in\{1,\dots,N\}$, and $n\in\{1,\dots,N-l\}$, 
there exist smooth functions $m_\dt^{n,l}(z)$, $M_\dt^{n,l}(z)$ such that 
the functions $\{\eta^n_\dt(z)\}$ corresponding to the semi-discretization of $(t,t+h)=(n\dt,(n+l)\dt)$ can be bounded as follows:
\begin{itemize}
\item $m_\dt^{n,l}(z)\leq \eta^{n+i}_\dt(z)\leq M_\dt^{n,l}(z),\; z\in [0,L],\; i=0,\dots,l,$ and
\item $M_\dt^{n,l}(z)-m_\dt^{n,l}(z)\leq C\sqrt{l\dt},\; z\in [0,L],\; {\rm and} \ \|M_\dt^{n,l}-m_\dt^{n,l}\|_{L^2(0,L)}\leq Cl\dt.$
\end{itemize}
\end{lemma}
\proof
The existence and properties of the functions $M$ and $m$ are a direct consequence of Proposition 5 in \cite{BorSun} (p. 29). 

To prove the second statement of Lemma~\ref{DoljeGore}, we fix $N \in \N$ and 
consider the finitely many functions $\eta^{n+i}_\dt(z)$ for $i = 0,\dots,l$, which are defined on the time interval 
from $t = n\dt$ to $t + h = n\dt + l\dt$. 
Define $m_\dt^{n,l}(z)$ and $M_\dt^{n,l}(z)$, $z\in [0,L]$ to be the functions obtained by considering the minimum and maximum of the finitely
many functions $\eta^{n+i}_\dt(z)$ for $i = 0,\dots,l$, mollified if necessary to get the smooth functions.
The properties of the functions $m_\dt^{n,l}$, $M_\dt^{n,l}$ then  follow from the proof of the same Proposition 5  in \cite{BorSun}. 
Namely, from the proof of Proposition 5 \cite{BorSun} it follows that the upper bound on $\|\eta_\dt^{n+i}-\eta^n_\dt\|_{H^1}$
only depends on the width of the time interval, which is $h$, namely $\|\eta_\dt^{n+i}-\eta^n_\dt\|_{H^1}\leq C\sqrt{h},$ for all $i = 1,\dots,l$.
Since $\eta^{n}_\dt$ are defined on $\overline\Gamma = [0,L]$, then the estimate also holds point-wise (because of the Sobolev embedding of $H^1(0,L)$ into $C[0,L]$). 
Thus the first inequality is proven since $h=l\dt$.
The second inequality, namely the $L^2$-bound on the difference between $M_\dt^{n,l}$ and $m_\dt^{n,l}$, 
follows from 
$$
\|\eta^n_\dt-\eta_0\|_{L^2(0,L)}\leq \sum_{i=0}^{n-1}\|\eta^{i+1}_\dt-\eta^{i}_\dt\|_{L^2(0,L)}
=\Delta t\sum_{i=0}^{n-1}\|v^{i+\frac 1 2}_\dt \|_{L^2(0,L)} \le \| v_\dt \|_{L^2(0,T,L^2(0,L))},
$$
where $\eta^{0}_\dt=\eta_0$, 
and from the uniform energy bound of the structure
velocity $v_\dt$ in $L^2(0,T;L^2(0,L))$, which follows from Lemma~\ref{EnergyEstimates}.
\qed

\noindent
{\bf Set up for the compactness result.} We are now in a position to define  the overarching function spaces for the compactness result.
Let  $\Omega^{M}$ be the maximal domain determined by the function $M(z)$ from Lemma~\ref{DoljeGore}
containing all the physical domains  $\Omega^{\eta^n_\dt}$. 
We define the spaces $H$ and $V$ as follows:
\begin{equation}\label{spaces}
H=L^2(\Omega^M)\times L^2(\Gamma), \quad
V=H^s(\Omega^M)\times H^s(\Gamma), 0<s<1/2.
\end{equation}
Notice that $H$ and $V$ are Hilbert spaces and that $V \subset\subset H$. 

\noindent
{\bf Extensions of fluid velocity approximations to $\Omega^M$.}
We consider the functions  $\bu^n_\dt$ defined on $\Omega^{\eta^{n-1}_\dt}$,
and the structure velocity approximations $v^n$ defined on $\Gamma$.
We {\bf{extend}} the functions $\bu^n_\dt$
to the maximal domain $\Omega^M$ by zero. 
\begin{lemma}\label{extensions}
The extensions by $0$ to $\Omega^M$ of $\bu^n_\dt$ are such that 
$$
(\bu^n_\dt,v^n_\dt)\in V = H^s(\Omega^M)\times H^s(\Gamma), 0<s<1/2.
$$
\end{lemma}
\proof
The proof is a consequence of Lemma~2.7 in \cite{Michailov}.
More precisely, to apply Lemma~2.7 in \cite{Michailov} we 
first notice that the functions $\eta_\dt^n$ are uniformly Lipschitz on $\Gamma$. 
This follows from the uniform energy estimates presented in Lemma~\ref{EnergyEstimates} above,
which 
imply that there exists a constant $C > 0$, independent of $\dt$, such that the $H^2(\Gamma)$-norms of $\eta_\dt^n$ are uniformly bounded by $C$ (estimate 4. in Lemma~\ref{EnergyEstimates}).
Since $\Gamma 
\subset \R$, this implies that $\eta_\dt^n$ are uniformly Lipschitz, where the Lipschitz constant is independent of $\dt$.
Now we can apply Lemma~2.7 in \cite{Michailov} to the functions $\bu^n_\dt$. 
First, the uniform energy estimates from Lemma~\ref{EnergyEstimates} imply that $\bu^n_\dt$ are in $H^1(\Omega^{\eta^n_\dt})$.
Then, Lemma~2.7 in \cite{Michailov} implies that the extensions by $0$ to $\Omega^M$ belong to $H^s(\Omega^M)$, $0 < s < 1/2$. 
\qed

Moreover,
there exists an absolute constant $C>0$, which depends only on the maximum of the Lipschitz constants for $\eta^n_\dt$, such that
$$
\|\bu^n_\dt\|_{H^s(\Omega^M)}\leq C\|\nabla\bu^n_\dt\|_{L^2(\Omega^{\eta^n_\dt})}, 0 < s < 1/2.
$$
Thus, by considering extensions by $0$ outside of the fluid domain $\Omega^{\eta^{n-1}_\dt}$, we see that $V^n_\dt$ is embedded in $V$. 

Similarly, we extend the test functions $\bq \in Q^n_\dt$ to the maximal domain $\Omega^M$ by $0$. 
Notice that the test functions are no longer smooth this way.
 However, as we shall see below, we will work with ``admissible'' test functions whose smoothness
 will depend only on the smoothness of $\bq$ within the domains $\Omega^{\eta^n_\dt}$.
 
We define now the velocity functions that depend on both time and space by introducing
$\{(\bu_\dt,v_\dt)\} \subset L^2(0,T;H)$  which are piecewise constant in $t$, i.e.
\begin{equation*}
\left.
\begin{array}{rcl}
\bu_\dt&=&\bu_\dt^n\\
v_\dt &=& v^n_\dt 
\end{array}
\right\}
\quad {\rm on}\; ((n-1)\dt,n\dt],\;n=1,\dots,N,
\end{equation*}
 as well as the corresponding time-shifts, denoted by $\tau_h$, defined by
$
\tau_h \bu_\dt (t,.)=\bu_\dt(t-h,.),\  h\in\R.
$

\begin{theorem}\label{CompactnessFSI}{\rm{\bf(Compactness result for FSI problem)}}
The sequence $\{(\bu_\dt,v_\dt)\}$, introduced in Sec.~\ref{sec:FSIsemi}, satisfying the weak formulation \eqref{Eulerian}
 energy estimates from Lemma \ref{EnergyEstimates},  is relatively compact in $L^2(0,T;H)$,
where $H=L^2(\Omega^M)^2\times L^2(\Gamma)$.
\end{theorem}

\proof
We would like to show that the assumptions (A)-(C) from Theorem~\ref{Compactness} 
hold true.

\noindent
\textbf{Property A: Strong bounds.}
The strong bounds (A1) and (A2) follow directly from the uniform energy bounds in Lemma~\ref{EnergyEstimates}
and from Lemma~\ref{extensions} above.  Because of Theorem \ref{A3Bez}, Property (A3) is not needed for the proof of Theorem \ref{CompactnessFSI}. However, it can be proved in the same way as Property (A3) in Example \ref{example:NS} by using the same calculation as presented below in the proof of Property B.

\noindent
\textbf{Property B: The time derivative bound.}
We want to estimate the following norm:
$$
\|P^n_\dt\frac{\bu^{n+1}_\dt-\bu^n_\dt}{\dt}\|_{(Q^n_\dt)'} = \sup_{\|(\bq,\psi)\|_{Q^{n}_\dt}=1}
\Big |\int_{\Omega^{\eta^n_\dt}}\frac{{\bf u}^{n+1}_\dt- {\bf u}^{n}_\dt}{\Delta t}\cdot{\bf q}+\int_0^L\frac{v^{n+1}_\dt-v^{n}_\dt}{\Delta t}\psi\Big |.
$$
For this purpose we add and subtract the function $\tilde\bu^n$ which is defined in \eqref{tilde_u}:
\begin{align}
\Big |\int_{\Omega^{\eta^n_\dt}}\frac{{\bf u}^{n+1}_\dt- {\bf u}^{n}_\dt \pm \tilde\bu^n}{\Delta t}\cdot{\bf q}+\int_0^L\frac{v^{n+1}_\dt-v^{n}_\dt}{\Delta t}\psi\Big |
&\leq \Big |\int_{\Omega^{\eta^n_\dt}}\frac{{\bf u}^{n+1}_\dt- \tilde\bu^n}{\Delta t}\cdot{\bf q}+\int_0^L\frac{v^{n+1}_\dt-v^{n}_\dt}{\Delta t}\psi\Big |
\nonumber
\\
&+\Big |\int_{\Omega^{\eta^n_\dt}}\frac{\tilde {\bf u}^{n}_\dt- \bu^n}{\Delta t}\cdot{\bf q} \Big |,
\nonumber
\end{align}
and estimate the two terms on the right hand-side.

The first term is estimated by using the weak form of the problem given by equation \eqref{Eulerian}:
$$
\Big |\int_{\Omega^{\eta^n_\dt}}\frac{{\bf u}^{n+1}_\dt-\tilde{\bf u}^{n}_\dt}{\Delta t}\cdot{\bf q}+\int_0^L\frac{v^{n+1}_\dt-v^{n}_\dt}{\Delta t}\psi\Big |
$$
$$
\leq C\frac{R+M}{R+m}\|\nabla {\bf q}\|_{L^{\infty}}(\|v^{n+1/2}_\dt\|_{L^2}+\|{\bf u}^{n}_\dt\|_{L^2})\|\nabla {\bf u}^{n+1}_\dt\|_{L^2}
+ C_1 \|\nabla {\bf u}^{n+1}_\dt\|_{L^2}\|\nabla{\bf q}\|_{L^2}
$$
$$
+\|\eta^{n+1}_\dt\|_{H^2}\|\psi\|_{H^2}+C\|{\bf q}\|_{H^1}
\leq C(\|\nabla {\bf u}_\dt^{n+1}\|_{L^2}+\|\eta_\dt^{n+1}\|_{H^2}+1)\|({\bf q},\psi)\|_{(Q^{n,l}_\dt)}.
$$
Here we used the energy estimates from Lemma~\ref{EnergyEstimates} 
from where we concluded that $\|\bu^{n}_\dt\|_{L^2}$, $\|v^{n}_\dt\|_{L^2}$ are uniformly bounded by $C$. 
Notice how the choice of the space $Q^{n}_\dt$, which includes high regularity Sobolev  spaces, is useful in the last inequality 
to provide the  upper bound in terms of $\|({\bf q},\psi)\|_{(Q^{n}_\dt)}$.

To estimate the second term,
we first notice that function $\tilde{\bf u}^{n}$ is non-zero on $\Omega^{\eta^{n}}$, while function  ${\bf u}^{n}$ is non-zero on $\Omega^{\eta^{n-1}}$.  
This is why we introduce $A=\Omega^{\eta^{n-1}}\cap \Omega^{\eta^{n}}$,  $B_1=\Omega^{\eta^{n}}\setminus \Omega^{\eta^{n-1}}$, 
and $B_2=\Omega^{\eta^{n-1}}\setminus \Omega^{\eta^{n}}$,
and estimate the  integrals over $A$, $B_1$, and $B_2$ separately:
$$
|\int_{A}({\bf\tilde u}^{n}-{\bf u}^{n})\cdot {\bf q}|=|\int_A\Big ({\bf u}^{n}(z,r)-{\bf u}^{n}(z,\frac{R+\eta^{n}}{R+\eta^{n-1}}r)\Big )\cdot{\bf q}(z,r)dzdr|
\leq C \Delta t\|v^{n-1/2}\|_{L^2}\|\nabla{\bf u}^{n}\|_{L^2(A)}\|{\bf q}\|_{L^{\infty}(A)}.
$$
To estimate the integral over $B_1$, we use the fact that $\bu^n = 0$ on $B_1$ to obtain:
$$
|\int_{B_1}({\bf\tilde u}^{n}-{\bf u}^{n})\cdot {\bf q}|
\leq | \int_0^Ldz\int_{R+\eta^{n-1}(z)}^{R+\eta^{n}(z)}({\bf\tilde u}^{n}-{\bf u}^{n})(z,r)\cdot{\bf q}dr|
=  \int_0^Ldz\int_{R+\eta^{n-1}(z)}^{R+\eta^{n}(z)}  {\bf\tilde u}^{n}  (z,r)\cdot{\bf q}dr|
$$
$$
\leq \|{\bf q}\|_{L^{\infty}}\int_0^L\max_{r}(\tilde{\bu}^{n}(z,r))dz\int_{R+\eta^{n-i}(z)}^{R+\eta^{n}(z)}dr\leq C\|{\bf q}\|_{L^{\infty}}\int_0^L\|\partial_ru^{n}_r(.,z)\|_{L^2_r}|\Delta t v^{n-1/2}(z)|dz 
$$
$$
\leq C\Delta t\|{\bf q}\|_{L^{\infty}}\|\nabla{\bf u}^{n}\|_{L^2(\Omega^{n}_\dt)}.
$$
Here, we used the fact that for $f \in H^1(0,1)$ we can estimate $\|f\|_{L^\infty} \le C \|f'\|_{L^2}$, and applied this to the function ${\bf{u}}^{n}$ above, 
viewed as a function of only one variable, $r$,  satisfying the condition that the $r$-component ${{u}}_r(z,0)=0$ is equal to zero at $r = 0$,
and the z-component $u_z$ is equal to zero on the boundary $r = R+\eta^{n-1}(z)$.
This is a formal estimate which can be justified by density arguments.

The integral over $B_2$ is, in fact, equal to zero.


\noindent
\textbf{Property C: Smooth dependence of function spaces on time.}
To define a common function space that will help with estimating the time shifts by $h=l\dt$ of $\bu^n_\dt$, namely $\bu^{n,l}_\dt$, we recall that
our scheme is designed in such a way that the functions $\bu^{n+1}_\dt$ are defined on the ``previous'' domain $\Omega^{\eta^n_\dt}$.

\noindent
\textbf{Property C1: The common test space.}
Consider the maximum  and minimum domains $\Omega^{M^{n,l-1}_\dt}$ and $\Omega^{m^{n,l-1}_\dt}$ determined by the maximum and minimum functions $M^{n,l-1}_\dt$ and $m^{n,l-1}_\dt$
given by Lemma~\ref{DoljeGore}.
%
\if 1 = 0
TODO: Ja bih bez ove nove notacije. Lakse je pratiti.

In order to simplify notation we introduce the following notation:
$$
\Omega^{n,l-1}_\dt=\Omega^{M^{n,l-1}_\dt},\; \omega^{n,l-1}_\dt=\Omega^{m^{n,l-1}_\dt}.
$$
\fi
%
We define the  "common" test space required by the general property (C1) to be the space consisting of all the (smooth) functions $(\bq,\psi)$
such that $\bq$ is defined on $\Omega^{m^{n,l}_\dt}$ and then extended to the maximal domain $\Omega^{M^{n,l}_\dt}$ by the trace $\psi {\bf e}_r$
(which is constant in the  ${\bf e}_r$ direction):
$$
Q^{n,l}_\dt=\{(\bq,\psi)\in \big ( V_F(\Omega^{M^{n,l}_\dt})\cap H^3(\Omega^{M^{n,l}_\dt})\big )\times H^{3}(\Gamma):\bq=\psi{\bf e}_r\;{\rm on}\; \Omega^{M^{n,l}_\dt} \setminus\Omega^{m^{n,l}_\dt}\}.
$$

\begin{remark}{\rm{
Notice that $Q^{n,l}_\dt\subset Q^{n+i}_\dt$, for all $i=0,\dots,l$, and the functions from $Q^{n,l}_\dt$ are admissible for problem \eqref{Eulerian} 
for all domains 
$\Omega^{\eta^{n+i}_\dt}$, $i=0,\dots, l$. In particular, they are divergence-free for all $i = 0,\dots,l$.}}
\end{remark}

We can now define the common test space for all the time shifts by $i\dt$, where $i = 0,\dots,l$ in the following way.
First notice that the test functions from $Q^{n,l}_\dt$ are non-zero on $\Omega^{M^{n,l}}$.
Since $\Omega^{M^{n,l}}\supset\Omega_\dt^{\eta^{n+i}_\dt}$ we can define the test functions belonging to $Q^{n+i}_\dt$
to be the  {\bf restrictions} onto $\Omega^{\eta^{n+i}_\dt}$ of the test functions from $Q^{n,l}_\dt$. More precisely, we define
\begin{equation}\label{Ex2B1}
\bq^i_{\dt,l,n}=\left\{\begin{array}{l}
\bq,\; {\rm in}\; \Omega^{\eta^{n+i}_\dt}
\\ 0,\; {\rm elsewhere}
\end{array}\right . , 
 \ \;\bq\in Q^{n,l}_\dt.
\end{equation}
Now, the mappings $J^i_{\dt,l,n} : Q^{n,l}_\dt \to Q^{n+i}_\dt$ are defined by:
$$
J^i_{\dt,l,n}(\bq,\psi)=(\bq^i_{\dt,l,n},\psi).
$$
Indeed, from the definition of the test space $Q^{n,l}_\dt$ we can see that $J^i_{\dt,l,n}(\bq,\psi)\in Q^{n+i}_\dt$. 
Moreover, we immediately see that 
$$\|J^i_{\Delta t,l,n}\bq\|_{Q^{n+i}_\dt}\leq C\|\bq\|_{Q^{n,l}_\dt}.$$ 

To check that inequality \eqref{Ji1} holds,
we use Lemma \ref{DoljeGore} and compute:
$$
\Big ((J^{j+1}_{\dt,l,n}(\bq,\psi)-J^j_{\dt,l,n}(\bq,\psi)),(\bu^{n+j+1}_\dt,v^{n+j+1}_\dt)\Big )_H
=|\int_0^L\int_{\eta^{n+j}_\dt(z)}^{\eta^{n+j+1}_\dt(z)}\psi(z)\bu^{n+j+1}_\dt(z,r)drdz|
$$
$$
\leq \|\psi\|_{L^{\infty}}\|\nabla\bu^{n+j+1}_\dt\|_{L^2} \|\eta^{n+j+1}_\dt-\eta^{n+j}_\dt\|_{L^1(\Gamma)}
\leq C \|(\bq,\psi)\|_{Q^{n,l}_\dt} \|\bu^{n+j+1}_\dt\|_{V^{n+j+1}_\dt}\dt,\;j\in\{0,\dots,l-1\}.
$$
The penultimate inequality is proved analogously as in the proof of Property B.

To check that inequality \eqref{Ji2} holds, we recall that $J^i_{\dt,l,n} \bq$ and $\bq$ differ only in the region 
$\Omega^{M} \setminus\Omega^{\eta^{n+i}_\dt}$, and so the $H$-norm difference between the two functions 
can be bounded by the $Q^{n,l}_\dt$-norm of $\bq$ and the difference between $M^{n,l}_\dt$ and $m^{n,l}_\dt$,
which is bounded by $C \sqrt{l \dt}$ according to  Lemma~\ref{DoljeGore}.
This completes the construction of a common test space for the time shifts, as specified by Property C1.

\noindent
\textbf{Property C2: Approximation property of solution spaces.}
We would like to be able to compare and estimate the time shifts by $h$ of the fluid velocity function $\bu$, given at some time $t+h$,
with the function $\bu$ given at time $t$.
The time-shift and the function itself are defined on different physical domains, since they also depend on time. 
We would like to define a common solution space on which we can make the function comparisons,
where the common function space needs to be constructed in such a way that the functions from that space 
approximate well the original functions defined on different domains.
The common space will be defined on a domain that is contained in all the other domains (we can do this because
our approximate domains are close). The functions corresponding to the time-shifts are going to be defined on this common domain
by a ``squeezing'' procedure, also used in \cite{CDEM},
and then extended to the largest domain by the trace on the boundary in a way that keeps the divergence free condition satisfied. 
Keeping in mind that $h = l\dt$, we introduce the space
$V^{n,l}_\dt$ to be the closure of  $Q^{n,l}_\dt$ in $V$ (as in the general case), and introduce the mappings $I^i_{\dt,l,n}$:
$$
V^{n,l}_\dt=\overline{Q^{n,l}_\dt}^{V},\quad 
I^i_{\dt,l,n}:V^{n+i}_\dt \to V^{n,l}_\dt, \  i=0,1,\dots,l,
$$
which will be used to approximate the functions in $V^{n+i}_\dt$ by the corresponding approximations defined in the common space $V^{n,l}_\dt$.
For each $i$, the mapping $I^i_{\dt,l,n}$ will be defined via a ``squeezing'' operator, defined below as follows (see also  \cite{CDEM}):
\begin{definition}[Squeezing]\label{def:squeezing}
Let $\eta_m$, $\eta$, $\eta_M$ be three functions defined on $[0,L]$ such that $-R<\eta_m(z)\leq\eta(z)\leq\eta_M(z)$, $z\in [0,L]$,
so that $R+\eta_m$ defines $\Omega^{\eta_m}$ with $R+ \eta_m > 0$. 
Let $\bu$ be a divergence-free function defined on $\Omega^{\eta}$ such that $\bu=v{\bf e}_r$ on $\Gamma^{\eta}$. For a given $\sigma\geq 1$, such that $\sigma\eta_m\geq \eta$ we define $\bu_{\sigma}\in H^1(\Omega^{\eta_M})$ in the following way:
\begin{equation}\label{Squeez}
{\bf u}_{\sigma}(z,r)=\left \{\begin{array}{cl}(\sigma u_z(z,\sigma r),u_r(z,\sigma r)),\;  &\sigma r\leq R+\eta(z),\\ v{\bf e}_r,\;&{\rm elsewhere}.\end{array}\right .
\end{equation}
\end{definition}
This defines a ``squeezing'' operator, which associates to any function $\bu$ defined on $\Omega^\eta$
a function $\bu_\sigma$ defined on the large domain $\Omega^{\eta_M}$, containing all the important information about the 
original function $\bu$ squeezed within the minimal domain $\Omega^{\eta_m}$. The operator is designed by: 
(1)  
first squeezing the function $\bu$ from domain $\Omega^\eta$ 
into $\Omega^{\eta_m}$ and rescaling the function $u_z$ by $\sigma$ so that the divergence free condition
remains to be satisfied, and 
(2) extending the squeezed $\bu$ to the remainder of the maximal domain $\Omega^{\eta_M}$ by the values of its trace on $\eta$, 
where the extension is constant in the ${\bf e}_r$ direction. 
The resulting function is divergence free.

\begin{remark}{\rm{
Notice that $\nabla\cdot{\bf u}=0$ implies that also $\nabla\cdot\bu_{\sigma}=0$. Moreover if $({\bf u},v)\in H^s(\Omega^{\eta})\times H^s(0,L)$, then ${\bf u}_\sigma\in H^s(\Omega^{\eta_M})$, $0\leq s\leq 1$.}}
\end{remark}

Now that the squeezing operator has been defined, for each $i=0,1,\dots,l$ we can define the mappings 
$I^i_{\dt,l,n}$ 
by squeezing the functions defined on $\Omega^{\eta^{n+i}_\dt}$ onto the smallest, ``common''  domain $\Omega^{m^{n,l}_\dt}$, 
and extending them by the trace on the boundary to the largest domain $\Omega^{M^{n-1,l+1}_\dt}$,
as in \eqref{Squeez}:
\begin{equation}\label{mappingsI}
I^i_{\dt,l,n} : \bu^{n+i}_\dt \mapsto (\bu^{n+i}_\dt)_{\sigma^i_{n,l}} \in V^{n,l}_\dt,
\end{equation}
where $(\bu^{n+i}_\dt)_{\sigma^i_{n,l}}$ is obtained from Definition~\ref{def:squeezing} by setting:
\begin{itemize}
\item $\eta_M =M^{n-1,l+1}_\dt$, $\eta=\eta^{n+i-1}_\dt$, $\eta_m=m^{n,l}_\dt$, and
\item $\sigma^i_{n,l} > 0$ is such that $\displaystyle{\sigma^i_{n,l}\geq \max_{z\in [0,L]}\frac{R+M^{n-1,l+1}(z)}{R+\eta^{n+i-1}_\dt(z)}}$.
\end{itemize}

We need to prove that there exists a universal constant $C>0$ such that
\begin{equation}\label{C2_1}
\| I^i_{\dt,l,n}(\bu^{n+i}_\dt ) \|_{V^{n,l}_\dt} \le C \| \bu^{n+i}_\dt  \|_{V^{n+i}_\dt },\quad i = 0,\dots,l
\end{equation}
and a universal, monotonically increasing function $g$, which converges to $0$ as $h\to0$, where $h = l\dt$, such that
\begin{equation}\label{C2_2}
\| I^i_{\dt,l,n}(\bu^{n+i}_\dt ) - \bu^{n+i}_\dt  \|_{L^2(\Omega^M)} \le g(l\dt) \| \bu^{n+i}_\dt  \|_{V^{n+i}_\dt }, \quad i = 0,\dots,l.
\end{equation}
The first inequality follows from Lemma~2.7 in \cite{Michailov} and 
the fact that the squeezed functions remain in $H^1$, with the uniformly bounded norm provided by the energy estimates 
in Lemma~\ref{EnergyEstimates}.

To show \eqref{C2_2} we first prove the following lemma, which states that the difference between $\bu$ and $\bu_\sigma$ on $\Omega^M$ can be estimated by
the $L^2$-norms of the gradient of $\bu$ on $\Omega^\eta$ and the $L^2$-norm of the velocity on the boundary, where the constant in the estimate
depends only on the Lipschitz constant of the fluid domain boundary.

\begin{lemma}\label{DenTmp}
Let $\eta_m$, $\eta$, and $\eta_M$ be as in Definition~\ref{def:squeezing}, 
$\bu\in H^1(\Omega^{\eta})$, $\bu=v{\bf e}_2$ on $\Gamma^{\eta}$,
and let $\sigma$ be such that $\displaystyle{\sigma\geq \max_{z\in [0,L]}\frac{R+\eta_M(z)}{R+\eta(z)}}$.
Then the following estimate holds:
$$
\|{\bf u}-{\bf u}_\sigma\|_{L^2(\Omega^{\eta_M})}\leq C\sqrt{\sigma-1} (\|\nabla{\bf u}\|_{L^2(\Omega^\eta)}+\|v\|_{L^2(0,L)}),
$$
where $C$ depends only on $\|R+\eta(z))\|_{L^{\infty}}$.
\end{lemma}
\proof
We will prove this estimate in two steps. 
First, we obtain an $L^2(\Omega^{\eta_M})$ estimate of the difference between $\bu$ and a slightly modified function $\tilde{\bu}_\sigma$,
defined by:
\begin{equation*}
\tilde{{\bf u}}_{\sigma}(z,r)=\left \{
\begin{array}{cl}(u_z(z,\sigma r),u_r(z,\sigma r)),\; &\sigma r\leq R+\eta(z),\\ 
v{\bf e}_r, \;&{\rm elsewhere},
\end{array}
\right. ,
\end{equation*}
and then estimate the difference between $\tilde{\bu}_\sigma$ and $\bu_\sigma$.
Notice that the only difference between $\tilde{\bu}_\sigma$ and $\bu_\sigma$ is the factor $\sigma$ in front of $u_z$. 

The $L^2$-norm of the difference between $\bu$ and $\tilde{\bu}_\sigma$ can be broken into three parts:
\begin{gather}
\int_{\Omega^{\eta_M}} | \bu(z,r) - \tilde{\bu}_\sigma(z,r) |^2 dr dz=
\int_0^L \int_0^{\frac{R+\eta(z)}{\sigma}}  | \bu(z,r) - \tilde{\bu}_\sigma(z,r) |^2 dr dz
\nonumber
\\
+
\int_0^L \int_{\frac{R+\eta(z)}{\sigma}}^{R+\eta}  | \bu(z,r) - \tilde{\bu}_\sigma(z,r) |^2 dr dz
+
\int_0^L \int_{R+\eta}^{R+\eta_M}  | \bu(z,r) - \tilde{\bu}_\sigma(z,r) |^2 dr dz
= I_1 + I_2 + I_3.
\nonumber
\end{gather}
The first part, $I_1$, can be estimated as follows:
$$
I_1=\int_0^L \int_0^{\frac{R+\eta(z)}{\sigma}} |{\bf u}(z,r)-{\bf u}(z,\sigma r)|^2drdz
=
\int_0^L \int_0^{\frac{R+\eta(z)}{\sigma}}  |\partial_r{\bf u}(z,\xi)|^2|r(1-\sigma)|^2drdz
\leq
C\|\nabla{\bf u}\|_{L^2(\Omega^{\eta_m})}^2(1-\sigma)^2.
$$
Notice that the constant $C$ depends on the size of  $\Omega^{\eta_m}$, which can be estimated from above by the size of
the maximal domain. 

To estimate $I_2$ we recall that $\tilde{\bu}_\sigma$ is defined to be $v{\bf e}_r$ outside the squeezed domain bounded by $(R+\eta(z))/\sigma$,
which is exactly the trace of $\bu$ on $R+\eta(z)$. Therefore, we get:
\begin{gather}
I_2 =\int_0^L  \int_{\frac{R+\eta}{\sigma}}^{R+\eta}{|\bf u}(z,r)-{\bf u}(z,R+\eta(z))|^2dr\ dz
\leq\int_0^L\int_{\frac{R+\eta}{\sigma}}^{R+\eta}|\partial_r{\bf u}(z,\xi)|^2(R+\eta(z))^2(1-\frac{1}{\sigma})^2
\nonumber
\\
\leq C (\sigma-1)^2\|\nabla{\bf u}\|_{L^2(\Omega^{\eta_m})}^2.
\nonumber
\end{gather}
Finally, to estimate $I_3$ we recall that $\bu$ is extended by zero outside $\Omega^\eta$, and that, as before, $\tilde{\bu}_\sigma$ 
is  equal to $v{\bf e}_r$ in the same region. Therefore, we obtain:
$$
I_3=\int_0^L \int_{R+\eta}^{R+\eta_M}|v(z)|^2 dr dz \leq \int_0^L|v(z)|^2(\eta_M-\eta)
\leq \|v\|_{L^2(0,L)}^2(\sigma-1)\|R+\eta(z))\|_{L^{\infty}},
$$
where the last inequality holds due to the assumption 
$\displaystyle{\sigma\geq \max_{z\in [0,L]}\frac{R+\eta_M(z)}{R+\eta(z)}}$.

In order to complete the proof and take into account  the multiplication of $u_z$ by $\sigma$, we use the triangle inequality:
$$\| u_z- (u_\sigma)_z \|_{L^2}\leq \| u_z-(u_\sigma)_z \pm  (\tilde{u}_\sigma)_z  \|_{L^2}
= \| u_z - (\tilde{u}_\sigma)_z \|_{L^2} + \| u_z \|_{L^2} (\sigma-1) 
\leq C\sqrt{\sigma-1} (\|\nabla{\bf u}\|_{L^2}+\|v\|_{L^2}). $$
This proves the lemma.
\qed

To show \eqref{C2_2} we use this lemma to get:
$$
\|I^i_{\dt,n,l}{\bf u}^{n+i}_\dt -{\bf u}^{n+i}_\dt \|_{L^2(\Omega^{\eta_M})}\leq C\sqrt{\sigma^i_{n,l}-1} \| {\bf u}^{n+i}_\dt \|_ {V^{n+i}_\dt}
\le C\sqrt{\sigma_{n,l}-1} \| {\bf u}^{n+i}_\dt \|_ {V^{n+i}_\dt},
$$
where 
$$
\sigma_{n,l} := \max_{z\in[0,L]} \frac{R+M^{n-1,l+1}_\dt}{R+m^{n,l}_\dt}.
$$
It is now crucial to show that we can estimate $(\sigma_{n,l} - 1)$ appearing on the right hand-side of the above inequality
by a function $g(h)$ such that $g \to 0$ as $h \to 0$, where $h = l\dt$.
Indeed, we have:
$$
\sigma_{n,l} - 1 = \frac{R+M^{n-1,l+1}_\dt}{R+m^{n,l}_\dt} - \frac{R+m^{n,l}_\dt}{R+m^{n,l}_\dt} = \frac{M^{n-1,l+1}_\dt - m^{n,l}_\dt}{R+m^{n,l}_\dt} 
\leq C (M^{n-1,l+1}_\dt - m^{n,l}_\dt) \leq C \sqrt{l \dt},
$$
where the last inequality follows from Lemma~\ref{DoljeGore}. Thus, we have shown that $g(h) = C h^{1/4}$, and the following 
density results holds:

\begin{corollary}[Density]\label{density}
Let $n,l,N\in\N$ be such that $n+l\leq N$, $h=l\Delta t$ and $i\in\{0,\dots, l\}$. Let $(\bu,v)\in V^{n+i}_\dt$. 
Then there exits a constant $C$ (independent of $N$, $n$, $l$) and the operators
$I^i_{\dt,n,l}$ defined by \eqref{mappingsI} such that 
$I^i_{\dt,n,l} ({\bf u},v)\in V^{n,l}_\dt$ and
$$
\|I^i_{\dt,n,l} {\bf u}-{\bf u}\|_{L^2}\leq C h^{1/4} \| {\bf u}^{n+i}_\dt \|_ {V^{n+i}_\dt}, \ i\in \{0,\dots,l\}.
$$
\end{corollary}

\noindent
\textbf{Property C3: Uniform Ehrling Property.}
We need to prove the Uniform Ehrling Property, stated in \eqref{Ehrling}.
The main difficulty comes, again, from the fact that we have to work with moving domains,
which are parameterized by $\dt,n$ and $l$. To show that the uniform Ehrling estimate \eqref{Ehrling} holds,
independently of all three parameters, we simplify the notation (only in this proof) and replace the indices $\dt,n,l$ 
with only one index, $n$, so that, e.g., the minimum and maximum functions $m^{n,l}_\dt$, $M^{n,l}_\dt$
are re-enumerated as $M_n$, $m_n$, with the corresponding maximal fluid domains $\Omega^{M_n}$
and the function spaces $H_n$, $V_n$ and $Q_n'$ all defined on $\Omega^{M_n}$.

We prove the Uniform Ehrling Property by contradiction, and by using the 
"uniform squeezing property" (\cite{NageleDiss}, Lemma 5.7).
We start by assuming that the statement of the uniform Ehrling property \eqref{Ehrling} is false.
More precisely, we assume that there exists a $\delta_0>0$ and a sequence $({\bf f}_n,g_n)\in H_n$ such that
$$
\|({\bf f}_n,g_n)\|_{H} =  \|({\bf f}_n,g_n)\|_{H_n}> \delta_0\|({\bf f}_n,g_n)\|_{V_n}+n\|({\bf f}_n,g_n)\|_{Q_n'}.
$$
Here, as before, we have extended the functions $f_n$ onto the entire domain $\Omega^M$, which is determined
by the maximal function $M(z)$ defined in Lemma~\ref{DoljeGore},
so that $\|({\bf f}_n,g_n)\|_{H} =  \|({\bf f}_n,g_n)\|_{H_n}$.
Recall that $H$ is defined as the $L^2$ product space, so extensions by $0$ do not change the norm. 

It will be convenient to also replace the $V_n$ norm on the right hand-side by the norm on $V$. Here, however,
since $V$ is defined as the $H^s$ product space, the norm on $V_n$ is bounded from below by a constant times the norm on $V$, see Lemma \ref{extensions}.
Thus, we have:
$$
\|({\bf f}_n,g_n)\|_{H}
>\delta_0\|({\bf f}_n,g_n)\|_{V_n}+n\|({\bf f}_n,g_n)\|_{Q_n'}\geq C\delta_0\|({\bf f}_n,g_n)\|_{V}+n\|({\bf f}_n,g_n)\|_{Q_n'}.
$$

Without the loss of generality we can assume that our sequence $({\bf f}_n,g_n)$ is such that 
$\|({\bf f}_n,g_n)\|_{H}=1$. For example, we could consider $\frac{1}{\|({\bf f}_n,g_n)\|_{H}}({\bf f}_n,g_n)$ instead of $({\bf f}_n,g_n).$
Notice that we now have a sequence such that the two terms on the right hand-side are uniformly bounded in $n$,
which implies that there exists a subsequence, which we again denote by $({\bf f}_n,g_n)$, such that:
\begin{equation}\label{sequence}
\|({\bf f}_n,g_n)\|_{H}=1,\quad \|({\bf f}_n,g_n)\|_{V}
\leq \frac{1}{C\delta_0},\quad \|({\bf f}_n,g_n)\|_{Q_n'}\to 0.
\end{equation}
Since $({\bf f}_n,g_n)$ is uniformly bounded in $V$, and by the compactness of the embedding of $V$ into $H$, we conclude 
that there exists a subsequence $({\bf f}_{n_k},g_{n_k})\to ({\bf f},g)$ strongly in $H$.

Now, from the energy estimates, we recall that sequences $\{M_n\}$, $\{m_n\}$ are uniformly bounded in $H^2(\Gamma)$,
and so by the compactness of the embedding of $H^2(\Gamma)$ into $C(\bar{\Gamma})$
 there exist  subsequences $M_{n_k}\to \tilde{M}$, $m_{n_k}\to \tilde{m}$.
With a slight abuse of notation, we denote the convergent subsequences by index $n$. Therefore, we are now working with the convergent
sub-sequences
$$
M_{n}\to \tilde{M}, \  m_{n} \to \tilde{m}, ({\bf f}_{n},g_{n}) \to ({\bf f},g) \ {\rm in} \ H.
$$
These functions define the maximal domain $\Omega^{\tilde{M}}$, the minimal domain $\Omega^{\tilde{m}}$,
and the corresponding minimal and maximal domains that depend on $n$: $\Omega^{{m}_n}$ and $\Omega^{{M}_n}$.

Let $\tilde{S}$ be the strip 
between the minimal and maximal domains $\Omega^{\tilde{m}}$ and $\Omega^{\tilde{M}}$,
and $S_n$ be the strip between the minimal and maximal domains that depend on $n$:
$$
S_n = \Omega^{M_n} \setminus\Omega^{m_n}, \quad
\tilde{S} = \Omega^{\tilde M} \setminus\Omega^{\tilde m}.
$$
We now want to show that $({\bf f},g) = 0$ first outside $\Omega^{\tilde M}$, i.e., in $\Omega^M\setminus\Omega^{\tilde M}$, then in $\Omega^{\tilde M}$.
(Recall that $\Omega^M$ was defined to contain all the approximate domains.)
This will contradict the assumption that $\|({\bf f}_n, g_n)\|_H = 1$ and complete the proof.

First, we show that ${\bf f} = 0$ outside of $\Omega^{\tilde M}$, and inside $\Omega^M$. This will follow simply because ${\bf f}$ is the limit of a sequence of 
functions that are zero outside $\Omega^{M_n}$ which converge to $\Omega^{\tilde M}$. More precisely,
we introduce the characteristic functions
 $\chi_{\Omega^{M_n}}$  of the sets $\Omega^{M_n}$ and recall that  $\chi_{\Omega^{M_n}}$ converge uniformly to $\chi_{\Omega^{\tilde{M}}}$.
Therefore:
$$
{\bf f}(1-\chi_{\Omega^{\tilde{M}}})
=\lim_{n}{\bf f}_n(1-\chi_{\Omega^{M_n}})=0,
$$
since ${\bf f}_n$ are extended by $0$ to $\Omega^M$. 

Next, we show that ${\bf f} = 0$ inside $\Omega^{\tilde M}$.
We start by showing that ${\bf f}=g{\bf e}_r$ in $\tilde S$.
This follows immediately from
$$
({\bf f}-g{\bf e}_r)\chi_{\tilde{S}}=\lim_n ({\bf f}_n-g_n{\bf e}_r)\chi_{S_n}=0,
$$
because $({\bf f}_n,g_n)\in V_n=\overline{Q_n}^{V}$, which implies ${\bf f}_n=g_n{\bf e}_r$ in $S_n$.

We finish the proof by using $\|({\bf f}_n,g_n)\|_{Q_n'}\to 0$ to show $({\bf f},g)=({\bf{0}},0)$. Recall that:
$$
Q^n=\{(\bq,\psi)\in \big (V_F(\Omega^{M_n})\cap H^4(\Omega^{M_n})\big )\times H^2_0(\Gamma):\bq=\psi{\bf e}_3\;{\rm on}\; S_n\}.
$$
Let $\varepsilon>0$. By the density of $Q^n$ in $H$, and by the uniform convergence of $M_n$ and $m_n$, 
combined with the ``squeezing procedure'' from the proof of Property C2,
we can construct a test function $(\bq,\psi)$ and $n_0\in\N$ such that:
\begin{equation}\label{TestEhrlingC}
\|({\bf f},g)-(\bq,\psi)\|_H\leq \varepsilon,
\quad
(\bq,\psi)\in Q^n,\; n\geq n_0.
\end{equation}
Therefore we have:
$$
\langle({\bf f},g),(\bq,\psi) \rangle_H
=\lim_n \langle({\bf f}_n,g_n),(\bq,\psi)\rangle_H
=\lim_n {_{Q_n'}} \langle ({\bf f}_n,g_n),(\bq,\psi) \rangle_{Q_n}
\leq \|({\bf f}_n,g_n)\|_{Q_n'}\underbrace{\|(\bq,\psi)\|_{Q_n}}_{\leq C} \to 0.
$$
Here, the uniform boundedness of $\|(\bq,\psi)\|_{Q_n}$ follows from the uniform convergence of $M_n$ and $m_n$.
Hence,
$$
\|({\bf f},g)\|_H^2
=\langle({\bf f},g),({\bf f},g)\pm(\bq,\psi)\rangle_H
=\langle({\bf f},g),({\bf f}-\bq,g-\psi)\rangle_H
\leq \varepsilon\|({\bf f},g)\|_H. 
$$
Therefore, $\|({\bf f},g)\|_H<\varepsilon$. Since $\varepsilon$ is arbitrary, $({\bf f},g)=({\bf 0},0)$, and this completes the proof of the Uniform Ehrling Property.

 Since this was the last step in verifying that all the assumptions of Theorem~\ref{Compactness} hold, we conclude that the statement of 
 the compactness theorem, Theorem~\ref{CompactnessFSI}, holds true for the sequence of functions approximating a solution of the 
 moving-boundary problem \eqref{Eulerian}. 
 \qed

\subsection{Fluid-structure interaction with the Navier slip condition}\label{sec:FSI_Navier}

In this section we only outline how Theorem~\ref{Compactness} can be applied to obtain a compactness result 
for a FSI problem involving the Navier slip boundary condition, studied in \cite{muha2016existence}.
The fluid domain and the fluid equations are the same as in Example~\ref{sec:FSI}. 
In contrast with Example~\ref{sec:FSI}, however, the structure model allows both the longitudinal and radial components of displacement
from the reference configuration $\Gamma=\Gamma_0=(0,L)$. 
The elastodynamics problem is given in terms of 
 a general  {\em continuous, self-adjoint, coercive, linear operator} ${\cal L}_e$, defined on $H_0^2(0,L)$, 
for which there exists a constant $c > 0$ such that 
\begin{equation}\label{coercivity}
\left<{\cal{L}}_e \bdeta,\bdeta \right> \ge c  \| \bdeta \|^2_{H_0^2(\Gamma)}, \quad \forall \bdeta \in H_0^2(\Gamma),
\end{equation}
where $\left< \cdot,\cdot \right>$ is the duality pairing between $H_0^2$ and $H^{-2}$.
The structure elastodynamics problem is defined by a clamped shell problem:
\begin{eqnarray}
\rho_{S} h  \partial_{tt} \bdeta + {\cal L}_e \bdeta = {\bf f}, \quad & z\in \Gamma=(0,L),\ t\in (0,T),
\label{Koiter}
\\
\bdeta(t,0)=\partial_z\bdeta(t,0)=\bdeta(t,L)=\partial_z\bdeta(t,L)=0,\; & t\in (0,T),
\label{KoiterBC}
\\
\bdeta(0,z)=\bdeta_0,\; \partial_t\bdeta(0,z)={\bf v}_0,; & z\in \Gamma=(0,L),
\nonumber
\end{eqnarray}
where $\rho_S$ is the structure density, $h$ the elastic shell thickness, ${\bf f}$ is linear force density acting on the shell, $\bdeta=(\eta_z,\eta_r)$ is the shell displacement, and $\bdeta_0$ and ${\bf v}_0$ are the initial structure displacement and the initial structure velocity, respectively. 
The fluid and structure equations are coupled via two sets of coupling conditions, the kinematic and dynamic coupling condition,
where the kinematic coupling condition is the Navier slip condition:
\begin{itemize}
\item {\bf The\  kinematic\  coupling\ condition (Navier slip condition):}
\begin{equation}\label{Coupling1a}
\begin{array}{l}
{\rm Continuity \ of \ normal\ velocity\;on\;}\Gamma^{\eta}(t) {\rm (the\ impermeability\ condition)}:\\
\displaystyle{\partial_t\bdeta(t,z)\cdot\bnu^{\eta}(t,z)} = \boldsymbol{u}(\bphi(t,z))\cdot\bnu^{\eta}(t,z), 
\quad
(t,z) \in (0,T)\times\Gamma,
\\  \\
 {\rm The\  slip\ condition \ between\ the \ fluid\ and\ thin\ structure\;on\;}\Gamma^{\eta}(t):\\
\displaystyle{(\partial_t\bdeta(t,z)-{\bf u}(\bphi(t,z)))\cdot\bt^{\eta}(t,z)}
\\
\displaystyle{=\alpha\bsigma\big (\bphi(t,z)\big )\bnu^{\eta}(t,z)\cdot\bt}^{\eta}(t,z),\quad (t,z)\in (0,T)\times\Gamma.
\end{array}
\end{equation}
Here, $\alpha$ is the friction constant ($1/\alpha$ is the slip length), 
$\bphi$ is the domain deformation, 
$\bnu^\eta$ and $\bt^\eta$ denote the unit normal and tangent vectors to $\Gamma^\eta$, respectively,
and $\boldsymbol\sigma$ is the fluid Cauchy stress tensor.
\item {\bf The\  dynamic\  coupling\ condition:}
\begin{equation}
\label{Coupling1b}
\rho_{S} h \partial_{tt}\bdeta(t,z) =  -{\mathcal L}_e\bdeta(t,z) -S^{\eta}(t,z)\bsigma\big (\bphi(t,z)\big )\bnu^{\eta}(t,z),\; (t,z)\in (0,T)\times\Gamma,
\end{equation}
The term $S^\eta$ is the Jacobian of the transformation between the Eulerian 
and Lagrangian formulations
of the fluid and structure problems, respectively.
\end{itemize}
Notice that there is no pressure contribution in the slip condition \eqref{Coupling1a}.
The pressure contributes only through the dynamic coupling condition \eqref{Coupling1b}.

Therefore, the problem reads: 
find $({\bf u},\bdeta)$ such that

\begin{equation}
\left .
\begin{array}{rcl}
\rho_F (\partial_t{\bf u}+{\bf u}\cdot\nabla{\bf u})&=&\nabla\cdot\bsigma,
\\
\nabla\cdot{\bf u}&=&0,
\end{array}
\right \}\ {\rm in}\ \Omega^{\eta}(t),\ t\in (0,T);
\label{NSP}
\end{equation}

\begin{equation}\label{SlipP}
\left.
\begin{array}{rcl}
\rho_{S} h \partial_{tt}\bdeta(t,z) +{\mathcal L}_e\bdeta(t,z) &=& -S^{\eta}(t,z)\bsigma\big (\bphi(t,z)\big )\bnu^{\eta}(t,z),
\\
\displaystyle{\partial_t\bdeta(t,z)\cdot\bnu^{\eta}(t,z)} &=& \boldsymbol{u}(\bphi(t,z))\cdot\bnu^{\eta}(t,z),
\\
\displaystyle{(\partial_t\bdeta(t,z)-{\bf u}(\bphi(t,z))\cdot\bt^{\eta}(t,z)} &=& \alpha\bsigma\big (\bphi(t,z)\big )\bnu^{\eta}(t,z)\cdot\bt^{\eta}(t,z),
\end{array}
\right\}
{\rm on} \ \Gamma, t\in(0,T),
\end{equation}
with
\begin{equation*}
\bdeta(t,0)=\partial_z\bdeta(t,0)=\bdeta(t,L)=\partial_z\bdeta(t,L)=0,\; t\in (0,T).
\end{equation*}
This problem is supplemented with the inlet/outlet and initial conditions:
\begin{equation}
p+\frac{\rho_F}{2}|{\bf u}|^2 =P_{in/out}(t),\; {\bf u}\cdot\bt=0\; {\rm on}\;\Gamma_{in/out}, t\in (0,T),
\label{ProblemBC}
\end{equation}
\begin{equation}
{\bf u}(0,.)={\bf u}_0,\;\bdeta(0,.)=\bdeta_0,\; \partial_t\bdeta(0,.)={\bf v}_0.
\label{ProblemIC}
\end{equation}
More general boundary conditions on the rigid part of the boundary were considered in \cite{muha2016existence}. 

The initial data must satisfy the compatibility conditions:
$
{\bf u}_0\in L^2(\Omega_0)^2,\; \nabla\cdot{\bf u}_0=0,\;  {\rm in}\; \Omega^0,
{\bf u}_0\cdot\bnu_0={\bf v}_0\cdot\bnu_0, \  {\rm on}\; \Gamma^0,
$
where $\Omega^0=\Omega^{\eta}(0)$, $\Gamma^0=\Gamma^{\eta}(0)$, $\bnu_0=\bnu^{\eta}(0,.)$. 
The fluid equations are defined on a cylinder with initial displacement 
$\bdeta_0 = 0$. A more general initial displacement was considered in \cite{muha2016existence}.

To deal with the motion of the fluid domain we again consider the ALE mapping ${\bf A}_\eta$, given by
(20) in \cite{muha2016existence}. 
This ALE mapping defines the fluid domain velocity ${\bf w}^\eta = d{\bf A}_\eta/dt$.

The weak formulation of problem 
\eqref{NSP}-\eqref{ProblemIC},
written in ALE formulation, is then given by the following (see  (43) in \cite{muha2016existence} pp. 8566):
find $({\bf u},\bdeta)$ such that
\begin{equation}\label{WeakALEFormula}
\begin{array}{c}
\displaystyle{-\rho_F\int_0^T\int_{\Omega}J^{\eta}{\bf u}\cdot\partial_t{\bf q} 
+ \frac{\rho_F}{2}\int_0^T\int_{\Omega}J^{\eta}\Big ((({\bf u}-{\bf w}^{\eta})\cdot\nabla^{\eta}){\bf u}\cdot {\bf q}-(({\bf u}-{\bf w}^{\eta})\cdot\nabla^{\eta}){\bf q}\cdot {\bf u}
-(\nabla^{\eta}\cdot{\bf w}^{\eta}){\bf q}\cdot {\bf u}\Big)}
\\ \\
\displaystyle{
+\int_0^T2\mu\int_{\Omega} J^\eta {\bf D}^{\eta}({\bf u}):{\bf D}^{\eta}({\bf q})
+\frac{1}{\alpha}\int_0^T \int_{\Gamma}(u_{t^{\eta}}-\partial_t\eta_{t^{\eta}})q_{t^{\eta}}S^{\eta}dzdt }
\\ \\
\displaystyle{-\rho_S h\int_0^T\int_{\Gamma} \partial_t\bdeta\partial_t\bpsi dz dt+\int_0^T \langle{\mathcal L}_e\bdeta,\bpsi\rangle
+\frac{1}{\alpha}\int_0^T \int_{\Gamma}(\partial_t\eta_{t^{\eta}} -u_{t^{\eta}})\psi_{t^{\eta}}S^{\eta}dz dt}
\\ \\
\displaystyle{
=\int_0^T\langle {\bf F},{\bf q}\rangle+\int_{\Omega}J_0{\bf u}_0\cdot{\bf q}(0)+\int_{\Gamma}{\bf v}_0\cdot\bpsi},
\quad \forall ({\bf q},\bpsi)\in{\mathcal Q}^{\eta}(0,T),
\end{array}
\end{equation}
where the sub-script $t^\eta$ denotes the tangent component of a given function,
$S^\eta$ is the Jacobian of the transformation between Eulerian and Lagrangian coordinates,
$J^\eta$ is the determinant of the Jacobian of  the ALE mapping, 
and ${\bf w}^\eta$ is the ALE velocity.
The term $\langle {\bf F},{\bq}\rangle$ denotes the contribution from the inlet/outlet data, as in \eqref{in_out_data}.

The solution space and the test space for the problem with the Navier slip boundary condition are given by the following. 
The solution space incorporates the continuity of normal velocities:
\begin{equation}\label{W_eta}
{\mathcal W^\eta}(0,T)=\{({\bf u}^\eta,\bdeta)\in {\cal W}_F^\eta(0,T)\times{\cal W}_S(0,T):{\bf u}^\eta_{|\Gamma}\cdot\bnu^{\eta}=\partial_t\bdeta\cdot\bnu^{\eta}\},
\end{equation}
and similarly, the corresponding test space is defined by
\begin{equation}\label{Q_eta}
{\mathcal Q^\eta}(0,T)=\{({\bf q}^\eta,\bpsi)\in C^1_c([0,T);{\cal V}_F^\eta\times H_0^2(\Gamma)^2):{\bf q}^\eta_{|\Gamma}\cdot\bnu^{\eta}=\bpsi(t,z)\cdot\bnu^{\eta}\}.
\end{equation}
Here:
\begin{equation}\label{VF_eta}
{\mathcal V}_F^{\eta}=\{{\bf u}^\eta=(u_z^\eta,u_r^\eta)\in H^1(\Omega)^2:\nabla^{\eta}\cdot{\bf u}^\eta=0,
{\bf u}^\eta\cdot\btau =0,\; {\rm on}\ \Gamma_{in/out} \},
\end{equation}
\begin{equation}\label{FluidSpace}
{\mathcal W}_F^\eta(0,T)=L^{\infty}(0,T;L^2(\Omega))\cap L^2(0,T;{\cal V}_F^\eta),
\end{equation}
\begin{equation}\label{struc_test}
{\mathcal W}_S(0,T)=W^{1,\infty}(0,T;L^2(0,L))^2\cap L^{\infty}(0,T;H^2_0(\Gamma))^2.
\end{equation}

The existence of a weak solution to this problem was proved in \cite{muha2016existence} in a constructive way,
using Rothe's method. A sequence of approximate solutions $\{{\bf u}_{\Delta t}, \bdeta_{\Delta t}\}_{\Delta t > 0}$ was constructed,
as in Example~\ref{sec:FSI} of this paper. The definition of the approximate solution sequence
$\{{\bf u}_{\Delta t}, \bdeta_{\Delta t}\}_{\Delta t > 0}$ is presented in Section~5.1 in \cite{muha2016existence},
with the uniform energy estimates satisfied by $\{{\bf u}_{\Delta t}, \bdeta_{\Delta t}\}_{\Delta t > 0}$ given in Proposition~4 \cite{muha2016existence}.
The goal was to show that as $\Delta t \to 0$, a sub-sequence of approximate solutions converges to a weak solution of the 
coupled FSI problem. 
Crucial for the existence result were compactness arguments for the fluid and structure velocities,
for the structure displacement, and for the geometric quantities which include 
the unit normal and tangent to the current location of the fluid domain boundary $\bnu_{\Delta t}$ and $\bt_{\Delta t}$,
the gradient of the ALE mapping $\nabla {\bf A}_{\Delta t}$, the Jacobian of the ALE mapping $J_{\Delta t}$, and the 
Jacobian of the transformation between Eulerian and Lagrangian coordinates $S_{\Delta t}$. 

The following convergence results were proved in \cite{muha2016existence} (see Theorem 5 and Corollaries 2 and 3) for the geometric quantities:
\begin{enumerate}
\item $\displaystyle{\bdeta_{\Delta t}\to\bdeta\;{\rm in}\; L^{\infty}(0,T;H^{2s}(\Gamma))}$, and
   $\displaystyle{\tau_{\Delta t}\bdeta_{\Delta t}\to\bdeta\;{\rm in}\; L^{\infty}(0,T;H^{2s}(\Gamma))}$, $s < 1$,
\item $\displaystyle{\bdeta_{\Delta t}\to\bdeta\;{\rm in}\; L^{\infty}(0,T;C^1(\overline{\Gamma}))}$, and
   $\displaystyle{\tau_{\Delta t}\bdeta_{\Delta t}\to\bdeta\;{\rm in}\; L^{\infty}(0,T;C^1(\overline{\Gamma}))}$,
\item $\bnu_{\Delta t}\to \bnu^{\eta}\;{\rm in}\; L^{\infty}(0,T;C(\overline{\Gamma}))$, and
   $\bt_{\Delta t}\to \bt^{\eta}\;{\rm in}\; L^{\infty}(0,T;C(\overline{\Gamma}))$,
\item $S_{\Delta t}\to S^{\eta}\;{\rm in}\; L^{\infty}(0,T;C(\overline{\Gamma}))$, and
   $J_{\Delta t}\to J^{\eta}\;{\rm in}\; L^{\infty}(0,T;C(\overline{\Omega}))$,
\item $(\nabla{\bf A}_{\Delta t})^{-1}\to (\nabla{\bf A}_{\eta})^{-1}\;{\rm in}\; L^{\infty}(0,T;C(\overline{\Omega}))$.
\end{enumerate}
What we will show in the present manuscript is how Theorem~\ref{Compactness} applies to showing compactness of
the fluid and structure velocities.

For this purpose we first recall from \cite{muha2016existence}, Section 3.2, that due to the low regularity in time of the fluid domain motion, namely, $\bdeta \in L^\infty(0,T;H^2(\Gamma))\cap W^{1,\infty}(0,T;L^2(\Gamma))$,
the function ${\bf w}^\eta$, which is the time derivative of the ALE mapping, is not in $H^1(\Omega_\eta)$. 
This is why we need to first explain how to interpret the 
term $-\frac{\varrho_F}{2}\int_0^T\int_{\Omega}J^{\eta}(\nabla^{\eta}\cdot{\bf w}^\eta){\bf q}\cdot\bu$ 
in the weak formulation \eqref{WeakALEFormula} above, which is (43) in \cite{muha2016existence} pp. 8566.
We formally recalculate this term on the moving domain as follows:
$$
-\frac{\varrho_F}{2}\int_0^T\int_{\Omega}J^{\eta}(\nabla^{\eta}\cdot{\bf w}^\eta){\bf q}\cdot\bu
=-\frac{\varrho_F}{2}\int_0^T\int_{\Omega^{\eta}(t)}(\nabla\cdot{\bf w}){\bf q}\cdot\bu
=\frac{\varrho_F}{2}\int_0^T\int_{\Omega^{\eta}(t)}{\bf w}\cdot\nabla({\bf q}\cdot\bu)
$$
$$
-\frac{\varrho_F}{2}\int_0^T\int_{\Gamma^{\eta}(t)}{\bf w}\cdot{\bf n}({\bf q}\cdot\bu)d\Gamma^{\eta}
=\frac{\varrho_F}{2}\Big (\int_0^T\int_{\Omega^{\eta}(t)}\Big (({\bf w}\cdot\nabla)\bu\cdot\bq+({\bf w}\cdot\nabla)\bq\cdot\bu\Big ) -\int_0^T\int_{\Gamma}({\bf v}\cdot{\bf n})({\bf q}\cdot\bu)S^{\eta}\Big).
$$
Now, the right-hand side is well defined, and we take it as the definition of the integral involving $\nabla^\eta\cdot{\bf w}$.
Consequently,  the entire convective term can now be written as:
$$
\frac{\rho_F}{2}\int_0^T\int_{\Omega}J^{\eta}\Big ((({\bf u}-{\bf w}^{\eta})\cdot\nabla^{\eta}){\bf u}\cdot {\bf q}-(({\bf u}-{\bf w}^{\eta})\cdot\nabla^{\eta}){\bf q}\cdot {\bf u}
-(\nabla^{\eta}\cdot{\bf w}^{\eta}){\bf q}\cdot {\bf u}\Big)
$$
$$
=-{\varrho_F}\Big (\int_0^T\int_{\Omega^{\eta}(t)}((\frac{1}{2}\bu-{\bf w}^{\eta})\cdot\nabla)\bq\cdot\bu
+\frac{1}{2}\int_0^T\int_{\Omega^{\eta}(t)}(\bu\cdot\nabla)\bu\cdot\bq
-\frac{1}{2}\int_{\Gamma}({\bf v}\cdot{\bf n})({\bf q}\cdot\bu)S^{\eta}\Big).
$$


For the compactness, the only questionable term is the boundary integral
$\int_{\Gamma}({\bf v}\cdot{\bf n})({\bf q}\cdot\bu)S^{\eta}$. 
Namely, with $\bu_{\Delta t}$ in $H^1$, by the trace theorem
 we have that ${\bu}_{\Delta t} \cdot \bq$ is bounded in $L^2(0,T;H^{1/2}(\Gamma))$.
 Moreover, $\bn_\dt S_\dt$ is bounded in $L^{\infty}(0,T;H^1(\Gamma))$ (because $\bdeta$ is in $H^2$), 
 and therefore $(\bu_\dt\cdot\bq)\bn_{\Delta t} S_\dt$ is bounded in $L^2(0,T;H^{1/2}(\Gamma))$.
  Therefore, to be able to pass to the limit in this integral, 
we only need the convergence of
${\bf v} _{\Delta t} \cdot {\bf n}_{\Delta t}$ in $H^{-1/2}(\Gamma)$. 

In fact, in contrast with the no-slip case, even if we wanted to {\sl directly}
%
obtain compactness of 
${\bf v} _{\Delta t} \cdot {\bf n}_{\Delta t}$ in $L^2(0,T;L^2(\Gamma))$ we would not be able to do it.
The reason for this is that, in order to obtain compactness of ${\bf v} _{\Delta t} \cdot {\bf n}_{\Delta t}$ in $L^2(0,T;L^2(\Gamma))$, 
we would, heuristically speaking, have to bound
$\partial_t(\bv\cdot\bn)=\partial_t\bv\cdot\bn+\bv\cdot\partial_t\bn$ in the sense of some dual norm. 
However, we only have $\bv\in L^{\infty}(0,T;L^2(\Gamma))$ and $\bn\in L^\infty(0,T;H^1(\Gamma))$, 
and therefore we cannot even define the product $\bv\cdot\partial_t\bn$ on $\Gamma$. 
Therefore, we will have to settle for a compactness result in a weaker norm.

From the physical point of view, the reason for  a weaker compactness result lies
in the fact that the fluid viscosity regularizes the evolution of the fluid-structure interface 
in a different way from the no-slip case, where the continuity of the fluid and structure velocities
holds in both normal and tangential components. 

However, we can obtain the strong convergence of the {\it normal} component of fluid velocities in $L^2$
indirectly, by using the following arguments. First, from the uniform boundedness of ${\bf v} _{\Delta t} \cdot {\bf n}_{\Delta t}$ in $L^2$ (from the energy estimates), 
we know that there exists a sub-sequence, also denoted by ${\bf v} _{\Delta t} \cdot {\bf n}_{\Delta t}$, which converges weakly in $L^2$.
If we can show that the norms $\| {\bf v} _{\Delta t} \cdot {\bf n}_{\Delta t} \|_{L^2(0,T;L^2(\Gamma))}$ converge, we would have strong convergence of ${\bf v} _{\Delta t} \cdot {\bf n}_{\Delta t}$ in $L^2$
(this is because we already have weak convergence in $L^2$).
We can, indeed, show the convergence of the norms from:
$$
\int_0^T\| {\bf v} _{\Delta t} \cdot {\bf n}_{\Delta t} \|^2_{L^2(\Gamma)} =\int_0^T{\phantom{a}}_{H^{-s}}\langle {\bf v} _{\Delta t} \cdot {\bf n}_{\Delta t} ,  {\bf v} _{\Delta t} \cdot {\bf n}_{\Delta t} \rangle_{H^{s}}
\to \int_0^T{\phantom{a}}_{H^{-s}}\langle {\bf v}  \cdot {\bf n},  {\bf v} \cdot {\bf n} \rangle_{H^{s}} = \int_0^T\| {\bf v} \cdot {\bf n} \|^2_{L^2(\Gamma)}, 
$$
where $0<s<1/2$, and the above convergence follows from the strong convergence of ${\bf v}_\dt  \cdot {\bf n}_\dt$ in $L^2(H^{-s})$, and the weak convergence of the same sequence in $L^2(H^{s})$,  for $0<s<1/2$.
The last statement is a consequence of the regularization by the fluid viscosity.
This way we showed that the fluid viscosity still regularizes the normal component of the structure velocity  via the no-penetration condition.
The dissipation due to friction in the tangential direction is, however, not sufficient to obtain strong $L^2$-convergence of the tangential component. 

To state the main compactness theorem we introduce the following overarching spaces:
\begin{equation}\label{spaces_slip}
H=L^2(\Omega^M)^2\times H^{-s}(\Gamma)^2, V=H^s(\Omega^M)^2\times L^2(\Gamma)^2, 0<s<1/2,
\end{equation}
where $\Omega^M$ is a rectangle $(0,L)\times (0,R_M)$. Here $(0,L)\times (0,R_M)$ is chosen in such a way that  $\Omega^n_\dt$,
which is determined by the boundary $\bdeta^n_\dt$, is contained in $(0,T)\times (0,R_M-\varepsilon)$, for all $\dt>0$, and $n\in\N$.
The main compactness result for the weak solution of the FSI problem  with the Navier slip condition \eqref{WeakALEFormula} is the following.
\begin{theorem}\label{CompactnessSlip}{\rm{\bf(Compactness result for FSI problem with Navier slip condition)}}
The sequence $\{(\bu_\dt,\bv_\dt)\}$, introduced in Section~5.1 in \cite{muha2016existence}, is relatively compact in $L^2(0,T;H)$,
where $H=L^2(\Omega^M)^2\times H^{-s}(\Gamma), 0<s<1/2$.
\end{theorem}

To prove this theorem we proceed the same way as in Example~\ref{sec:FSI}, except for a few notable differences, as made precise below.
The outline of the proof is as follows.

\vskip 0.2in
\noindent
{\bf Property A: Strong bounds.} These follow directly from the uniform energy estimates in Proposition 4 in \cite{muha2016existence}.

\vskip 0.2in
\noindent
{\bf Property B: The time derivative bound.} These are obtained in exactly the same way as in Example~\ref{sec:FSI}.

\vskip 0.2in
\noindent
{\bf Property C: Smooth dependence of function spaces on time.}
Before we outline the main ideas behind verifying properties (C1)-(C3), we introduce the following test space:
\begin{equation}\label{funky_function_space}
Q^n_\dt=\{(\bq,\bpsi)\in \big ( V_F(\Omega^{n}_\dt)\cap H^4(\Omega^{n}_\dt)\big )^2\times H^{2}_0(\Gamma)^2:\bq_{|\Gamma^{\eta_\dt^n}}=\bpsi\}.
\end{equation}
Notice that this test space has a different (stronger) condition at the interface: it requires that both velocity components of test functions
are continuous at the interface (no-slip). 
This test space is smaller then the test space in which continuity of only normal components of velocity is required.
However, it can be shown, using density arguments, that this test space approximates well the original test space \eqref{Q_eta} in appropriate topology..
Furthermore, the solution space is defined as follows:
\begin{equation}\label{Ex3:solution_space}
V^n_\dt=\{(\bu,\bv)\in \mathcal{V}_F(\Omega^{n-1}_\dt)^2\times L^2(\Gamma)^2:(\bu_{|\Gamma^{\eta_\dt^{n-1}}}-\bv)\cdot{\bf n}^{\bdeta^{n-1}_\dt}=0\}.
\end{equation}

\vskip 0.1in
\noindent
{\bf Property C1: The common test space.} 
To define a common test space, which contains all the fluid domains
determined by the time shifts $\tau_h \bu_\dt$ from $t$ to $t+h$, where $t = n\dt$ and $h = l\dt$,
we use the motivation explained above, and some additional geometric considerations.
Namely, due to the fact that in this example 
the tangential displacement is different from zero, the fluid domain may degenerate in such a way that
it ceases to be a subgraph of a function.
However, since we are interested in the local-in-time existence result presented in \cite{muha2016existence},
one can prove the following lemma, which guarantees the existence of a time during which the fluid domain will stay a sub-graph,
starting from a straight cylinder with initial displacement $\bdeta = 0$.
Let
$$
\boldsymbol\varphi_\dt^n = (R,z) + \bdeta_\dt^n(z).
$$
Analogously, as in Lemma \ref{DoljeGore}, one can prove, using regularity arguments, 
that there exists a time $T > 0$ such that 
not only the displacements in our FSI problem are small, but their tangents do not change by much,
namely that the following holds:
\begin{equation}\label{Ex3DG}
\|\bphi^{n}_\dt-(R,z)\|_{C^1[0,L]}\leq C(T),\; \dt>0,\; N\in\N,
\end{equation}
where $C(T)\to 0$, as $T\to 0$. 
By taking $T>0$ small enough, one can immediately show that there exists a $T>0$ such that the tangent lines
will stay between $-\pi/2$ and $\pi/2$, i.e., the fluid domain boundary will stay a sub-graph:
\begin{lemma}\label{Ex3Slope}
Let $\alpha\in (0,\pi/2)$. There exists a $T>0$ such that the slope of the tangent to $\Gamma^n_\dt$ is smaller then $\alpha$, i.e.
$$
-\alpha<\arccos(\frac{\varphi_z'(z)}{\|\bphi'(z)\|})<\alpha,\  z\in [0,L],\;\dt>0,\; n\in\N.
$$
\end{lemma}

This will imply an analogue of Lemma \ref{DoljeGore}, and the existence of the minimal and maximal fluid domains
associated with the time shifts $\tau_h$, with $h = l\dt$:
\begin{lemma}\label{Ex3DoljeGore}
For each fixed $\dt>0$, $l\in\{1,\dots,N\}$, $n\in\{1,\dots,N-l\}$, 
there exist smooth functions $m_\dt^{n,l}(z)$, $M_\dt^{n,l}(z)$ such that 
the domains $\Omega^{n+i}_\dt$, corresponding to the semi-discretization of the time interval $(t,t+h)=(n\dt,(n+l)\dt)$, satisfy the following:
\begin{itemize}
\item 
$\Omega^{m_\dt^{n,l}}\subset \Omega^{n+i}_\dt\subset \Omega^{M_\dt^{n,l}},\; i=0,\dots,l,$ and
\item $M_\dt^{n,l}(z)-m_\dt^{n,l}(z)\leq C\sqrt{l\dt},\; z\in [0,L],\; \|M_\dt^{n,l}-m_\dt^{n,l}\|_{L^2(0,L)}\leq Cl\dt.$
\end{itemize}
\end{lemma}
A common test space is then defined in the following way:
\begin{equation}\label{Ex1BSpace}
Q^{n,l}_\dt=\{(\bq,\bpsi)\in \big ( V_F(\Omega^{n,l}_\dt)\cap H^4(\Omega^{n,l}_\dt)\big )^2\times H^{2}_0(\Gamma)^2:\bq_{|\Gamma^n_\dt}=\bpsi\},
\end{equation} 
where $\Omega^{n,l}_\dt=\Omega^{M^{n,l}_\dt}$.

Let us take $({\bf q},\bpsi)\in Q^{n,l}_\dt$. We define $J^i$ as a restriction of $\bq$ onto the corresponding domain:
$$
J^i_{\dt,l,n}({\bf q},\bpsi)=({\bf q}_{|\Omega^{n+i}_\dt},{\bf q}_{|\Gamma^{n+i}_\dt})=:({\bf q}^i,\bpsi^i).
$$
To verify properties \eqref{Ji1} and \eqref{Ji2} of this mapping, we calculate
$\int_{\Omega}({\bf q}^{j+1}-{\bf q}^j)\cdot\bu^{n+j+1}+\|\bpsi^i-\bpsi^j\|_{L^2(\Gamma)}$. 
The first term is estimated in a same way as in  Example~\ref{sec:FSI} (no-slip), 
using the fact that $\bq$ is uniformly bounded, and Lemma \ref{Ex3DoljeGore}. Namely:
$$
|\int_{\Omega}({\bf q}^{j+1}-{\bf q}^j)\cdot\bu^{n+j+1}|
\leq \int_{\Omega^{j+1}_\dt\triangle\Omega^j_\dt}|\bq\cdot\bu^{n+j+1}|
\leq C\| \bq\|_{L^{\infty}}\dt\|\bu^{n+j+1}\|_{V^{n+j+1}_\dt}
$$
$$
\leq C\dt\|(\bq,\bpsi)\|_{Q^{n,l}_\dt}\|\bu^{n+j+1}\|_{V^{n+j+1}_\dt},
$$
where $A \triangle B := A\setminus B \cup B \setminus A$.
The second term is estimated using Lemma \ref{TraceLema} below, to obtain:
$$
\|\bpsi^i-\bpsi^j\|_{L^2(\Gamma)}\leq \|\bdeta^{n+i}_\dt-\bdeta^{n+j}_\dt\|_{L^2(\Gamma)}\|\nabla{\bf q}\|_{L^{\infty}}\leq C|i-j|\dt\|({\bf q},\bpsi)\|_{Q^{n,l}_\dt}.
$$

\begin{lemma}\label{TraceLema}
Let $\Omega'\subset\Omega^M$ and $\bphi^1,\bphi^2:\Gamma\to \Omega'$ be two curves such $\bphi^1=\bphi^2=\bdeta_0$ on $\partial\Gamma$ and $\bu\in H^1(\Omega')$. Then we have the following estimates:
$$
\|\bu(\bphi^1(.))-\bu(\bphi^2(.))\|_{L^2(\Gamma)}\leq \|\bphi^1-\bphi^2\|_{L^{\infty}(\Gamma)}\|\nabla\bu\|_{L^2(\Omega')},
$$
$$
\|\bu(\bphi^1(.))-\bu(\bphi^2(.))\|_{L^2(\Gamma)}\leq \|\bphi^1-\bphi^2\|_{L^{2}(\Gamma)}\|\nabla\bu\|_{L^{\infty}(\Omega')}.
$$
\end{lemma}
\proof
Let $z\in\Gamma$.
Define $g(s):=\bu(\bphi^1(z)+s(\bphi^2(z)-\bphi^1(z)))$. Then $g'(s)=\nabla\bu(\bphi^1(z)+s(\bphi^2(z)-\bphi^1(z)))\cdot(\bphi^2(z)-\bphi(z)$. 
Then, one gets:
$$
|\bu(\bphi^1(z))-\bu(\bphi^2(z))|=|\int_0^1g'(s)ds|\leq \int_0^1|\nabla\bu(\bphi^1(z)+s(\bphi^2(z)-\bphi^1(z))||\bphi^2(z)-\bphi^1(z)|ds.
$$
By taking the square on both sides, and integrating with respect to $z$,  $\int_{\Gamma}dz$, we obtain:
$$
\|\bu(\bphi^1(z))-\bu(\bphi^2(z))\|_{L^2(\Gamma)}^2\leq \int_{\Gamma}|\bphi^2(z)-\bphi^1(z)|^2\int_0^1|\nabla\bu(\bphi^1(z)+s(\bphi^2(z)-\bphi^1(z))|^2dsdz.
$$
\qed

This completes an outline of the verification of property (C1).

\vskip 0.1in
\noindent
{\bf Property C2: Approximation property of solution spaces.}
We define $V^{n,l}_\dt={\overline{Q^{n,l}_\dt}}^V$, which gives, as in \eqref{Ex3:solution_space}:
\begin{equation}\label{Ex3CV}
V^{n,l}_\dt=\{(\bu,\bv)\in H^s(\Omega^{n,l}_\dt)^2\times L^2(\Gamma)^2:\nabla\cdot\bu=0,\;(\bu_{|\Gamma^n_\dt}-\bv)\cdot{\bf n}^{\bdeta^n_\dt}=0\},
\end{equation}
for an $s$ such that $0 < s < 1/2$.
To verify that every function in $V^{n+i}_\dt$, given by the time shift $i\dt$,  can be approximated by a function in $V^{n,l}_\dt$,
one needs to construct a mapping $I_{\dt,l,n}^i:  V^{n+i}_\dt\to V^{n,l}_\dt$ with ``good approximation properties'' \eqref{C21}  and \eqref{C22}.
In the no-slip case where only the radial component of displacement was considered, the ``squeezing operator'', given by Definition~\ref{def:squeezing},
provided the desired properties. 
In the present case, when both radial and longitudinal components of displacement are allowed to be different from zero,
an extension operator, similar to the one in Section~\ref{example:NS}, has to be used
to verify Property C2.
Namely, for a function $({\bf u}_\dt^{n+i},{\bf v}_\dt^{n+i})\in V^{n+i}_\dt$ we will use a divergence free extension $\tilde{\bu}_\dt^{n+i}$ 
 given by  the following Lemma:
\begin{lemma}\label{Ex3Ext}
Let $\bu\in V^n_\dt$. Then there exists a divergence free function $\tilde{\bu}\in V$ such that $\tilde{\bu}|_{\Omega^{n-1}_\dt}=\bu$ and
$$
\|\tilde{\bu}\|_V\leq C\|\bu\|_{V^n_\dt},
$$
where $C$ is independent of $\dt,\; n$.
\end{lemma}
\proof
Let $\bu_1$ be extension of $\bu$ to $\R^2$, and $\bu_2$ the restriction of $\bu_1$ on $\Omega^M$, $\bu_2=(\bu_1)|_{\Omega^M}$. Here we used the existence of a ``total extension operator'' (see Adams \cite{ADA} Thm 5.24), which gives the following bound:
$$
\|\bu_2\|_{V}\leq C\|\bu\|_{V^{n}_\dt}.
$$
Unfortunately, $\bu_2$ is not divergence free.
To get around this difficulty, we first introduce $F=-\int_{\Omega^M}\nabla\cdot\bu_2$, 
and define 
$$
\bu_3=\bu_2+\frac{F}{L}\theta {\bf e}_r,
$$
where $\theta\in C^{\infty}_c([0,L]\times (R_M-\varepsilon,R_M]$ is such that $\theta(z,R_M)=1$.
This function has the property that $\int_{\Omega^M} \nabla \cdot \bu_3 = 0$. Moreover, $\bu_3$ is still an extension of $\bu_2$,
since the support of $\theta$ is disjoint from $\Omega^n_\dt$. 
Now, since $|F|\leq C\|\bu\|_{V^{n}_\dt}$, we have the following estimate on $\bu_3$:
$$
\|\bu_3\|_{V}\leq C\|\bu\|_{V^{n}_\dt}.
$$
We still need to correct $\bu_3$ to get a function which is divergence free on $\Omega^M$.
For this purpose we introduce a function 
$G\in H^1_0(\Omega)^2$ such that
$$
\nabla\cdot G=-\nabla\cdot\bu_3,
$$
and define $\tilde\bu=\bu_3+G$.
The existence of $G$ is well-known (see \cite{galdi2011introduction}, Thm. III.3.1). Moreover, one can construct a $G$ such that $G=0$ on $\Omega^n_\dt$  because $\nabla\cdot\bu_3=0$ on $\Omega^n_\dt$ (see \cite{galdi2011introduction} Lemma III.3.1).
With these considerations, we now have an extension of $\bu$ to $\Omega^M$ which is divergence free, and it satisfies the desired estimate:
$$
\|\tilde{\bu}\|_V\leq C\|\bu\|_{V^n_\dt}.
$$

The constant $C$ in this lemma does not depend on $n$ and $\dt$ because we are working with the domains that are ``well-controlled'' in the sense of Lemma~\ref{Ex3Slope}.

\qed

Using this lemma we define:
\begin{equation}\label{Ex3DefI}
I^i_{\dt,l,n}({\bf u}_\dt^{n+i},{\bf v}_\dt^{n+i})=\Big ((\tilde{\bu}_\dt^{n+i})|_{\Omega^{n,l}_\dt},((\tilde{\bu}_\dt^{n+i})|_{\Gamma^{n}_\dt}\cdot{\bf n}^n_\dt){\bf n}^n_\dt+({\bf v}_\dt^{n+i}\cdot\btau^n_\dt)\btau^n_\dt\Big ).
\end{equation}
It is clear from \eqref{Ex3CV} and \eqref{Ex3DefI}, and from Lemma \ref{Ex3Ext} that $I^i_{\dt,l,n}({\bf u}_\dt^{n+i},{\bf v}_\dt^{n+i})\in V^{n,l}_\dt$.
We need to estimate $\|I^i_{\dt,l,n}({\bf u}_\dt^{n+i},{\bf v}_\dt^{n+i})-({\bf u}_\dt^{n+i},{\bf v}_\dt^{n+i})\|_H.$ 
To simplify notation, we drop the subscripts $\dt,l,n$, and introduce the following notation $(I^i{\bf u}_\dt^{n+i},I^i{\bf v}_\dt^{n+i}):=I^i({\bf u}_\dt^{n+i},{\bf v}_\dt^{n+i})$.

We will show that we can obtain the desired estimate even in the stronger norm, namely, that we can estimate
$\|I^i{\bf u}_\dt^{n+i}-{\bf u}_\dt^{n+i}\|_{L^2(\Omega^M)}+\|I^i{\bf v}_\dt^{n+i}-{\bf v}_\dt^{n+i}\|_{L^2(\Gamma)}$. The first term is estimated by Lemma \ref{Ex3DoljeGore}:
$$
\|I^i{\bf u}_\dt^{n+i}-{\bf u}_\dt^{n+i}\|^2_{L^2(\Omega^M)}
=\int_{\Omega^{n,l}_\dt\triangle\Omega^{n+i}_\dt}\tilde{\bu}^{n+i}_\dt\leq C\sqrt{(l+1)\dt}.
$$
The second integral is estimated by Lemma \ref{TraceLema}:
$$
\|I^i{\bf v}_\dt^{n+i}-{\bf v}_\dt^{n+i}\|_{L^2(\Gamma)}
\leq \|(\tilde{\bu}_\dt^{n+i})|_{\Gamma^{n}_\dt}\cdot{\bf n}^n_\dt-(\tilde{\bu}_\dt^{n+i})|_{\Gamma^{n+i-1}_\dt}\cdot{\bf n}^{n+i-1}_\dt\|_{L^2(\Gamma)}
+\|{\bf v}^{n+i}_\dt ({\bf n}^{n+i-1}_\dt-{\bf n}^n_\dt)\|_{L^2}
$$
$$
\leq \|\big ((\tilde{\bu}_\dt^{n+i})|_{\Gamma^{n}_\dt}-(\tilde{\bu}_\dt^{n+i})|_{\Gamma^{n+i-1}_\dt}\big )\cdot{\bf n}^{n}_\dt\|_{L^2(\Gamma)}
+\|(\tilde{\bu}_\dt^{n+i})|_{\Gamma^{n+i-1}_\dt}\cdot({\bf n}^n_\dt-{\bf n}^{n+i-1}_\dt)\|_{L^2(\Gamma)}
$$
$$
+C\|{\bf v}^{n+1}_\dt\|_{L^2}\|\partial_z(\bdeta^n_\dt-\bdeta^{n+i-1}_\dt)\|_{L^{\infty}}
$$
$$
\leq C(\|{\bf n}^n_\dt\|_{L^{\infty}}\|\nabla\bu^{n+i}_\dt\|_{L^2}\|\bdeta^n_\dt-\bdeta^{n+i-1}_\dt\|_{L^{\infty}}
+(\|\nabla\bu^{n+i}_\dt\|_{L^2}+\|{\bf v}^{n+1}_\dt\|_{L^2})\|\partial_z(\bdeta^n_\dt-\bdeta^{n+i-1}_\dt)\|_{L^{\infty}})
$$
$$
\leq C\sqrt{l\dt}(\|\nabla\bu^{n+i}_\dt\|_{L^2}+\|{\bf v}^{n+1}_\dt\|_{L^2}).
$$
This finishes the proof of property (C2).

\vskip 0.1in
\noindent
{\bf Property C3: Uniform Ehrling Property.}
The proof of (C3) is identical to that associated with Example~\ref{sec:FSI}.

\qed

\section{Conclusions}
This work provides an extension of the Aubin-Lions-Simon compactness result to generalized Bochner spaces $L^2(0,T;H(t))$, where
$H(t)$ is a family of Hilbert spaces, parameterized by $t$. 
An example of a class of problems where such a compactness result is needed arises in the study of the existence of weak solutions to 
nonlinear problems governed by partial differential equations defined on moving domains, i.e., domains that depend on time. 
This work identifies the conditions on the regularity of the domain motion in time, under which an extension of the  Aubin-Lions-Simon compactness result holds. 
To demonstrate the use and the applicability of the abstract theorem, several examples were presented. They include a classical problem for the incompressible,
Navier-Stokes equations defined on a {\sl given} non-cylindrical domain, and a class of fluid-structure interaction problems for the incompressible, Navier-Stokes equations
coupled to the elastodynamics of a Koiter shell, via both the no-slip condition, and the Navier slip condition. 
This work supplements the results in \cite{BorSun,muha2016existence} by providing a consistent, general compactness proof,
which can be used in each particular case, studied in \cite{BorSun,muha2016existence}, as explained in Sec. \ref{sec:FSI} and Sec. \ref{sec:FSI_Navier}, above.

\section*{Acknowledgements}
This research has been supported in part by the National Science Foundation under grants
DMS-1613757, NIGMS DMS-1263572, DMS-1318763 (\v{C}ani\'c), DMS-1311709 (\v{C}ani\'c and Muha), and by
Croatia-USA bilateral grant ``Fluid-elastic structure interaction with the Navier slip boundary condition'' associated with NSF DMS-1613757.
Support from the Cullen Chair Fund at the University of Houston is also acknowledged (\v{C}ani\'c).


\if 1 = 0
****************OLD NAVIER SLIP*****************************

\subsection{FSI with slip}
\subsubsection{Monolithic formulation in the physical domain}

*****************************************************************
Let $\kappa>0$ be such that $\Lambda:\Gamma^{\bdeta_0}\times (-\kappa,\kappa)\to \R^2$, $(\bdeta_0(z),s)\mapsto \bdeta_0(z)+s\bn(z)$ is injective, where $\bn(z)$ is the unit outer normal at point $\bdeta_0(z)$, $z\in [0,L]$ (existence of such $\kappa$ is well known, see e.g. \cite{Lee} Thm 10.19, p. 254). Let us denote $S_0=\Lambda(\Gamma^{\bdeta_0}\times (-\kappa,\kappa))$. Furthermore, let us take $T_0$ such that $\Gamma^{\bdeta}(t)\subset S_0$, $t\in [0,T_0]$. Furthermore, because of uniform convergence of $\bdeta_\dt$, we chose $\dt_0$ such that 
$$
\Gamma^{\bdeta_\dt^n}\subset S_0,\quad \dt<\dt_0,\; n=0,\dots,N.
$$

We will repeatedly use the following simple lemma:
\begin{lemma}\label{TraceLema}
Let $\Omega\subset\R^2$ and $\bdeta^1,\bdeta^2:\Gamma\to S_0$ be two curves such $\bdeta^1=\bdeta^2=\bdeta_0$ on $\partial\Gamma$ and $\bu\in H^1(\Omega)$. Then we have the following estimates:
$$
\|\bu(\eta^1(x))-\bu(\eta^2(x))\|_{L^2(\Gamma)}\leq \|\bdeta^1-\bdeta^2\|_{L^{\infty}(\Gamma)}\|\nabla\bu\|_{L^2(\Omega)},
$$
$$
\|\bu(\eta^1(x))-\bu(\eta^2(x))\|_{L^2(\Gamma)}\leq \|\bdeta^1-\bdeta^2\|_{L^{2}(\Gamma)}\|\nabla\bu\|_{L^{\infty}(\Omega)},
$$
\end{lemma}
\proof
Let $x\in\Omega$.
We define $g(s):=\bu(\eta^1(x)+s(\eta^2(x)-\eta^1(x)))$. Then $g'(s)=\nabla\bu(\eta^1(x)+s(\eta^2(x)-\eta^1(x)))\cdot(\eta^2(x)-\eta^1(x)$. We compute:
$$
|\bu(\eta^1(x))-\bu(\eta^2(x))|=|\int_0^1g'(s)ds|\leq \int_0^1|\nabla\bu(\eta^1(x)+s(\eta^2(x)-\eta^1(x)))||(\eta^2(x)-\eta^1(x)|ds.
$$
We square and integrate $\int_{\Gamma}dx$ to obtain:
$$
\|\bu(\eta^1(x))-\bu(\eta^2(x))\|_{L^2(\Gamma)}^2\leq \int_{\Gamma}|\eta^2(x)-\eta^1(x)|^2\int_0^1|\nabla\bu(\eta^1(x)+s(\eta^2(x)-\eta^1(x)))|^2dsdx.
$$
\qed

TODO: Budi malo precizniji kakve su krivulje, tj. segmenti po kojima integriras (mogu se sjeci). U slucaju kada imamo samo radijalni pomak se ne mogu, ali kod slipa treba jos dodatno argumentirati (iako je sve OK) koristeci cinjenicu iz diferencijalne geometrije da se okolina glatke krivulje moze parametrizirati za normalama (onda normale sluze kao segmenti).

***********************************************************************

\subsubsection{Kratki opis - heuristika}
\begin{itemize}
\item Nije dobro raditi konvergenciju normalnih kmponenta $\bv\cdot\bn$. Naime ne mozemo ograniciti derivaciju $\partial_t (\bv\cdot\bn)$ jer nam probleme stvara clan $\bv\cdot\partial_t\bn$ ($\partial_t \bn$ jer samo $H^{-1}$, $\bv$ je samo $L^2$ pa ih ne mozemo mnoziti). Kod abstraktnog teorema to se manifestira tako da uvjet B1 nije zadovoljen. Naime $\|\bn^{n+1}-\bn^n\|\leq C \sqrt{\dt}$, a a treba nam ocjena $C\dt$ koju ne mozemo dobiti!!
\item Ako dodamo viskoelasticnost, tada sve funkcionira - dovoljno $\varepsilon\partial_t\triangle\bv$.
\item Iz energetskih ocjena znamo $\partial_t\bv\cdot\btau$ je ogranicena (cak u $L^2$). Poanta je dokazati da je i $\partial_t\bu+\partial_t \bv\cdot\bn$ takodjer ogranicena (u dualnoj normi), pa onda i $\partial_t\bu+\partial_t\bv$.
\item Sada koristimo gornje zapazanje i cinjenicu da je $\partial_t \bv$ ogranicena u $L^2$ pa zbog toga kompaktna u $H^{-s}$ da dobijemo kompaktnost $(\bu,\bv)$ u $L^2\times H^{-s}$. 
\item Cilj nam je preci na limes u clanu:
$$
\int_0^T\int_{\Gamma^{\eta}(t)}({\bf v}\cdot\bn^{\eta}){\bf u}\cdot{\bf q}d\Gamma^{\eta}(t)
$$
Medjutim zbog teorema o tragu ${\bf u}\cdot{\bf q}$ je ogranicen u $H^{1/2}$ pa nam je dovoljno konvergencija od $\bv$ u $H^{-s}$.
\item Abstraktni teorem se primjenjuje na sljedeci nacin:
\begin{itemize}
\item Svaka test funkcija (sa evaluacijom na $\Gamma^0$) je oblika $(\bq,\bq_{|\Gamma}+\psi\btau^0)$.
\item Preslikavanje $J^i$ definiramo kao evaluaciju na $\Gamma^i$:
$$
(\bq,\bq_{|\Gamma}+\psi\btau^0)\mapsto (\bq,\bq_{|\Gamma^i}+\psi\btau^0)=(\bq+\bq_{|\Gamma^i}+\tilde{\psi}^i\bn^i+\psi\btau^i)
$$
Ovo nije dobro zbog slicnih problema kao u prethodnoj tocku (pogreska je reda velicini $\sqrt{\dt}$ zbog geometrijskih velicina).
\item Nacin za izvuci se je sljedeci. Promatramo test prostore koji se sastoje od parove $(\bq,\bq_{|\Gamma^n})$. Tada je $J^j(\bq,\bq_{|\Gamma^n})=(\bq,\bq_{|\Gamma^{n+j}})$. Naravno, taj test prostor je manji od pravog test prostora koji zahtjeva neprekidnost samo normalne komponente. Medjutim, taj test prostor je ipak gust u $H^s$ normi za $s<1/2$ sto nam je dovoljno da primjenimo abstraktni teorem i zakljucimo konvergenciju (TODO: oprez, treba o ovome jos dobro razmisliti i provjeriti sve detalje).

\end{itemize}

\end{itemize}

Dodatno treba paziti na geometriju. Dakle na sljedece stvari:
\begin{itemize}
\item Kako definirati traku oko funkcija $\eta^i_\dt$. Mozda najbolje napraviti samo uniju domena.
\item Prosirenje se sada ne moze definirati formulom, vec moramo napraviti divergence free prosirenje (dokazati da ne ovisi o funkcijama $\eta^i_\dt$). Poanta je da se napravi na velikoj domeni, pa se onda samo gledaju restrikcije. Treba se igrati sa cone condition i zvjezdastim domenama, te teoremima o Bogovski operatoru iz Galdija i prosirenjima iz Adamsa. Na kraju se sve svodi na to kako nasa deformacija deformira konus, tj. da li se u deformirani konus moze opisati manji konus. Nasa deformacija je oblika $I+B^n$, gdje je $B^n$ mali u $C^1$ normi. To bi trebalo iskoristiti. Takodjer, vidjeti kakve to veze ima sa kvazikonformnim preslikavanjima.
\item Zadnji korak kod dualne ocjene derivacije u vremenu treba biti oprezan. Naime vise ne mozemo tako jednostavno razbiti integral. Treba koristiti geometriju - mozda sve napisati u bazi koja se sastoji od normale i tangente na limes $\eta$. To je vezano i za prvo pitanje. 

Treba nekako iskoristiti da su sve krivulje blizu u $C^1$ normi. Mozda se to moze iskorisiti tako da pokazemo da kut izmedju tangenta ne moze biti velik, dakle nesto ovog tipa:
$$
\|\bdeta_1-\bdeta_2\|\leq \varepsilon (\|\bdeta_1\|+\|\bdeta_2\|),
$$
gdje je $\varepsilon$ mali, ali fiksiran.
\item Ako petljanje sa geometrijom previse komplicira stvari, mozemo uzeti da je $\bdeta_0$ identiteta i to ce nam olaksati dokaze. U originalnom slip clanku to i jest pretpostavka.

\end{itemize}

******************************************************************************************************

Treba preraditi u skladu sa gornjim komentarima

*********************************************************************************************

Let $\Omega_m$, $\Omega_M$ be such that $\Omega_m\subset\Omega^{\bdeta_\dt^n}\subset\Omega_M$, $\dt>0$, $n=1,\dots, N$. Notice that by taking $T$ small we can make $\Omega_M\setminus\Omega_m$ arbitrary small.

Moreover we define  (TODO + mozda izgladenje ako treba).
$$\Omega^{n,l}_\dt=\bigcup_{i=0}^l\Omega^{\bdeta_\dt^{n+i}}.$$
Let us define appropriate function spaces:
\begin{itemize}
\item $H=L^2(\Omega)^2\times H^{-s}(\Gamma)^2$, $V=H^s(\Omega)^2\times L^2(\Gamma)^2$, $0<s<1/2$. Then we have $(\bu^i_\dt,v^i_\dt)\in V$.
\item For given displacement $\theta$ We define fluid space $V_F^{\theta}$ same as in \cite{BorSun}. More precisely,
$$
V_F^{\theta}=\{\bv\in H^1(\Omega^\theta):\nabla\cdot\bv=0,\;  \bv\cdot{\bf e}_r=0\; {\rm on}\;\partial\Omega^\theta\setminus\Gamma^\theta\}
$$
Now we define solution space.
$$
V^i_\dt=\{(\bu,v)\in V_F^{\bdeta^{i-1}_\dt}\times H^{1/2}(\Gamma):\bu_{|\Gamma^{\eta_\dt^{i-1}}}\cdot{\bf n}^{\bdeta^{i-1}_\dt}=v\},
$$
Using extension by $0$ we see that $V^i_\dt$ is embedded in $V$. Let us now define test space:
$$
Q^i_\dt=\{(\bq,\psi)\in \big ( V_F(\Omega^{i-1,1}_\dt)\cap H^4(\Omega^{i-1,1}_\dt)\big )\times H^{2}(\Gamma):\bq_{|\Gamma^{\eta_\dt^i}}\cdot{\bf n}^{\bdeta^{i}_\dt}=\psi\}.
$$
\item We define "common" test space:
$$\big ( V_F(\Omega^{n-1,l}_\dt)\cap H^4(\Omega^{n-1,l}_\dt)\big )\times H^{2}(\Gamma):\bq_{|\Gamma^{\eta_\dt^i}}\cdot{\bf n}^{\bdeta^{i}_\dt}=\psi\},\; V^{n,l}_\dt=\overline{Q^{n,l}_\dt}^{V}.
$$
\end{itemize}

Let us check conditions A-C:
\begin{itemize}
\item[A)] We use the fact $H^s_0=H^s$, so by taking extension by $0$ we have $(\bu^i_\dt,v^i_\dt)\in V$. Moreover, there exists absolute constant $C$ such that $\|\bu^i_\dt\|_{H^s(\Omega)}\leq C\|\nabla\bu^i_\dt\|_{L^2(\Omega^{\eta^i_\dt})}$. Namely, such $C$ depend on maximum of Lipschitz constants for $\eta^i_\dt$ which can be bounded by uniform constant.
\item[B1)] Let us take $({\bf q},\phi)\in Q^{n,l}_\dt$. We define $J^i$ in the following way:
$$
J^i({\bf q},\phi)=({\bf q}_{|\Omega^{n+i-1,1}_\dt},{\bf q}_{|\bdeta^{n+i}_\dt}\cdot\bn^{\bdeta^{n+i}_\dt})=:({\bf q}^i,\phi^i)
$$

Let us calculate $\|{\bf q}^i-{\bf q}^j\|_{L^2(\Omega)}+\|\phi^i-\phi_j\|_{L^2(\Omega)}$. The first term is zero because both $\bq^i$ and $\bq^j$ are restrictions of the same function $\bq$ (here we extend $\bq^i$, $\bq^j$ by $\bq$). We estimate the second term in the following way.

(normale su parcijalne derivacije pa ih mozemo parcijalnom integracijom prebaciti na test funkciju!!)

The second term is estimated by Lemma \ref{TraceLema}:
$$
\|\phi^i-\phi^j\|_{L^2(\Gamma)}\leq \|\eta^{n+i}_\dt-\eta^{n+j}_\dt\|_{L^2(\Gamma)}\|\nabla{\bf q}\|_{L^{\infty}}\leq C|i-j|\dt\|({\bf q},\phi)\|_{Q^{n,l}_\dt}.
$$

\item[B2)]
Let us now verify property  B2. We extend $({\bf u}_\dt^{n+i},v_\dt^{n+i})$ in the following way:
$$
\tilde{\bf u}_\dt^{n+i}=\left \{\begin{array}{c}
{\bf u}_\dt^{n+i}\; {\rm in}\; \Omega^{\eta^{n+i}_\dt},\\
v_\dt^{n+i}{\bf e}_r,\; {\rm in}\; \Omega^{M_\dt^{n-1,l}}\setminus\Omega^{\eta^{n+i}_\dt}.
\end{array}\right .
$$
TODO: Ista napomena kao kod test funkcija, tu moze squeezing jer ne treba glatkoca.

Obviously $(\tilde{\bf u}_\dt^{n+i},v_\dt^{n+i})\in V^{n,l}_\dt$ and by Lemma \ref{DoljeGore}
$$
\|\tilde{\bf u}_\dt^{n+i}-{\bf u}_\dt^{n+i}\|_{L^2}^2=\int_0^L(v_\dt^{n+i})^2(M^{n-1,l}_\dt-\eta^{n+i-1}_\dt)^2\leq Cl\dt.
$$
\item[B3)] Let us prove uniform Ehrling lemma.

We will modify the usual proof of the Ehrling lemma and use so-called "uniform squeezing property" (\cite{Nagel}). Let $M$ be the maximal functions from Lemma \ref{DoljeGore}, and let $Q_M$, $V_M$, $H_M$ be the corresponding function spaces, i.e. $V_M$, $H_M$ are closure of $Q_M$ in $\|.\|_V$ and $\|.\|_H$ respectively. 

For the simplicity of the notation let us re-enumerate functions $M^{n,l}_N$ with just one parameter and work with $M_n$ and corresponding domain $\Omega_n$ and function spaces $H_n$, $V_n$ and $Q_n'$.

We will prove the lemma by contradiction. Let us assume that Lemma is false, i.e. there exists $\varepsilon_0$ and sequence $({\bf f}_n,g_n)\in H^s_n$ such that
$$
\|({\bf f}_n,g_n)\|_{H_n}> \varepsilon_0\|({\bf f}_n,g_n)\|_{V_n}+n\|({\bf f}_n,g_n)\|_{Q_n'}.
$$

We can extend each ${\bf f}_n$ by $g_n{\bf e}_r$ on $\Omega_M\setminus \Omega_n$, and denote obtained function by $\tilde{\bf f}_n$. Notice that $(\tilde{\bf f}_n,g_n)\in V_M$. Furthermore,
$$
\|(\tilde{\bf f}_n,g_n)\|_{H_M}>\|({\bf f}_n,g_n)\|_{H_n}>\varepsilon_0\|({\bf f}_n,g_n)\|_{V_n}+n\|({\bf f}_n,g_n)\|_{Q_n'}\geq C\varepsilon_0\|(\tilde{\bf f}_n,g_n)\|_{V_M}+n\|({\bf f}_n,g_n)\|_{Q_n'}.
$$
Without lose of generality we cane take $\|(\tilde{\bf f}_n,g_n)\|_{H_M}=1$ (by considering $\frac{1}{\|(\tilde{\bf f}_n,g_n)\|_{H_M}>}(\tilde{\bf f}_n,g_n)$ instead of $(\tilde{\bf f}_n,g_n)).$

To summarize we constructed sequence $(\tilde{\bf f}_n,g_n)\in H_M$ such that
$$
\|(\tilde{\bf f}_n,g_n)\|_{H_M}=1,\quad \|(\tilde{\bf f}_n,g_n)\|_{V_M}\leq \frac{1}{C\varepsilon_0},\quad \|({\bf f}_n,g_n)\|_{Q_n'}\to 0.
$$
By compactness of embedding of embedding $V_M$ to $H_M$ we can extract subsequence $(\tilde{\bf f}_{n_k},g_{n_k})\to (\tilde{\bf f},g)$ strongly in $H_M$. Moreover, Since sequence $(M_n)$ is bounded in $H^2(\Gamma)$, by compactness of embedding $H^2(\Gamma)$ in $C(\bar{\Gamma})$ we can extract subsequence $M_{n_k}\to \eta$.

Let $({\bf q},0)\in Q_M$ such that supp(${\bf q}$) is compact subset of $\Omega_\eta$. Because of uniform convergence of $\eta_{n_k}$, there exists $k_0$ such that $({\bf q},0)\in Q_{n_k},\; k\geq k_0$. Therefore we have 
$$
\int_{\Omega_M}{\bf f}\cdot{\bf q}=\lim_k\int_{\Omega_M}\tilde{\bf f}_{n_k}\cdot{\bf q}\leq \lim_k\|({\bf f}_{n_k},g_{n_k})\|_{Q_n'}\|{\bf q}\|_{H^3}=0.
$$
Since ${\bf q}$ is arbitrary divergence free function with compact support in $\Omega_\eta$ and ${\bf f}$ is divergence free, by Helmholtz decomposition we conclude ${\bf f}=0$ in $\Omega_{\eta}$. Let now $({\bf q},\phi)\in Q_M$
$$
\int_{\Gamma}g\phi=\lim_k\int_{\Gamma}g_{n_k}\phi_k\pm\int_{\Omega_{n_k}}{\bf f}_{n_k}\cdot{\bf q}\leq \lim_k\int_{\Omega_{n_k}}|{\bf f}_{n_k}\cdot{\bf q}|=\int_{\Omega_{n_k}\setminus\Omega_{\eta}}|{\bf f}_{n_k}\cdot{\bf q}|,\quad ({\bf q},\phi_k)\in Q_{n_k}. 
$$
(TODO: Ovdje smo mogli cak uzeti $q$ takvi da je $\phi=\phi_k$. U slipu to nece ici i trebat cemo uzeti vecu regularnost od ${\bf q}$, recimo $H^4$ (da je $W^{2,\infty})$). The second inequality is consequence of  $\|(\bq,\phi_{n_k})\|_{Q_{n_k}}\leq \|\bq\|_{H^4}\leq C$ and  $\|({\bf f}_n,g_n)\|_{Q_n'}\to 0$.
Later integrals converges to $0$ because $|\Omega_{n_k}\setminus\Omega_{\eta}|\to 0$.
Therefore we conclude $g=0$ and thus $\tilde{\bf f}=0$. This contradicts $\|(\tilde{\bf f}_n,g_n)\|_{L^2_M}=1$ and therefore the Lemma is proven.

\item[C)] First notice that $({\bf u}^{i+1},0),(\bu^i,0)\in (Q_\dt^{i})'$. We compute:
$$
\|(\frac{\bu^{i+1}_\dt-\bu^i_\dt}{\dt},\frac{v^{i+1}_\dt-v^i_\dt}{dt})\|_{(Q^{i}_\dt)'}\leq \|(\frac{\bu^{i+1}_\dt-\tilde{\bu}^i_\dt}{\dt},\frac{v^{i+1}_\dt-v^i_\dt}{dt})\|_{(Q^{i}_\dt)'}+\|(\frac{\tilde{\bu}^{i}_\dt-\bu^i_\dt}{\dt},0)\|_{(Q^{i}_\dt)'}.
$$
The first one is estimated by using the equation \eqref{Eulerian}:
$$
\Big |\int_{\Omega^{\eta^n_\dt}}\frac{{\bf u}^{n+1}_\dt-\tilde{\bf u}^{n}_\dt}{\Delta t}\cdot{\bf q}+\int_0^L\frac{v^{n+1}_\dt-v^{n}_\dt}{\Delta t}\psi\Big |
$$
$$
\leq C\frac{R+M}{R+m}\|\nabla {\bf q}\|_{L^{\infty}}(\|v^{n+1/2}_\dt\|_{L^2}+\|{\bf u}^{n}_\dt\|_{L^2})\|\nabla {\bf u}^{n+1}_\dt\|_{L^2}
+ C_1 \|\nabla {\bf u}^{n+1}_\dt\|_{L^2}\|\nabla{\bf q}\|_{L^2}
$$
$$
+\|\eta^{n+1}_\dt\|_{H^2}\|\psi\|_{H^2}+C\|{\bf q}\|_{H^1}
\leq C(\|\nabla {\bf u}_\dt^{n+1}\|_{L^2}+\|\eta_\dt^{n+1}\|_{H^2}+1) \|({\bf q},\psi)\|_{(Q^{n}_\dt)}.
$$
Here we used the energy estimates (Prop. 4, p. 20 in \cite{BorSun}) from where we concluded that $\|\bu^{n}_\dt\|_{L^2}$, $\|v^{n}_\dt\|_{L^2}$ are uniformly bounded by $C$. 
Notice how the choice of the space $Q^{n}_\dt$, which includes high regularity Sobolev  spaces, is useful in the last inequality above
to provide the  upper bound in terms of $\|({\bf q},\psi)\|_{(Q^{n}_\dt)}$.
\end{itemize}
Let us estimate the second term.
Notice that function $\tilde{\bf u}^{n}$ is non-zero on $\Omega^{\eta^{n}}$, while function  $\tilde{\bf u}^{n}$ is non-zero on $\Omega^{\eta^{n-1}}$.  
Let us denote $A=\Omega^{\eta^{n-1}}\cap \Omega^{\eta^{n}}$, $B=\Omega^{\eta^{n}}\triangle \Omega^{\eta^{n-1i}}$. Now we can estimate  integrals over $A$ and $B$ separately.
$$
|\int_{A}({\bf\tilde u}^{n}-{\bf u}^{n})\cdot {\bf q}|=|\int_A\Big ({\bf u}^{i}(z,r)-{\bf u}^{i}(z,\frac{R+\eta^{n}}{R+\eta^{n-1}}r)\Big )\cdot{\bf q}(z,r)dzdr|
$$
$$
\leq \Delta t\|v^{n-1/2}\|_{L^2}\|\nabla{\bf u}^{n}\|_{L^2(A)}\|{\bf q}\|_{L^{\infty}(A)}
$$

The next estimate is on $B$. 
$$
|\int_{B}({\bf\tilde u}^{n}-{\bf u}^{n})\cdot {\bf q}|=| \int_0^Ldz\int_{R+\eta^{n-1}(z)}^{R+\eta^{n}(z)}({\bf\tilde u}^{n}-{\bf u}^{n})(z,r)\cdot{\bf q}dr|
$$
$$
\leq \|{\bf q}\|_{L^{\infty}}\int_0^L\max_{r}(u_r^{n}(z,r)|)dz\int_{R+\eta^{n-i}(n)}^{R+\eta^{n}(z)}dr\leq C\|{\bf q}\|_{L^{\infty}}\int_0^L\|\partial_ru^{n}_r(.,z)\|_{L^2_r}|\Delta t v^{n-1/2}(z)|dz 
$$
$$
\leq C\Delta t\|{\bf q}\|_{L^{\infty}}\|\nabla{\bf u}^{n}\|_{L^2(\Omega^{i}_\dt)}
$$
Here we used that if $f \in H^1(0,1)$, then $\|f\|_{L^\infty} \le C \|f\|_{H^1}$, where the function of one variable is $u^{n}$ above, viewed as a function of $r$ only and $u_r(z,0)=0$.

\subsection{Multilayered}
Sve isto samo paziti da li nam je negdje trebala Lipshitzovost granice.
\fi

\bibliographystyle{plain} 
\bibliography{../../myrefs}
\end{document}